\documentclass[11pt]{amsart}

\usepackage[utf8]{inputenc}

\usepackage{esvect}
\usepackage{amscd}
\usepackage{amsmath}
\usepackage{amssymb}
\usepackage{amsfonts}
\usepackage{amsthm}
\usepackage{bm}
\usepackage{bbm}
\usepackage{mathrsfs}
\usepackage{graphicx}
\usepackage{epsfig}
\usepackage{enumerate}
\usepackage[all]{xy}

\usepackage[usenames,dvipsnames]{xcolor}
\usepackage[colorlinks=true,linkcolor=NavyBlue,urlcolor=RoyalBlue,citecolor=PineGreen,bookmarks=true,bookmarksopen=true,bookmarksopenlevel=2,unicode=true,linktocpage]{hyperref}

\usepackage[top=25mm, bottom=25mm, left=25mm, right=25mm]{geometry}

\newtheorem{theorem}{Theorem}[section]

\newtheorem{proposition}[theorem]{Proposition}
\newtheorem{lemma}[theorem]{Lemma}
\newtheorem{corollary}[theorem]{Corollary}

\renewcommand{\1}{\mathbbm{1}}

\newcommand{\rA}{\mathsf{A}}

\newcommand{\B}{\mathcal{B}}
\newcommand{\cB}{\mathscr{B}}

\newcommand{\build}[3]{\mathrel{\mathop{\kern 0pt#1}\limits_{#2}^{#3}}}

\newcommand{\C}{\mathbb{C}}		

\newcommand{\calc}{{\sf c}} 
\newcommand{\cC}{\mathcal{C}}

\newcommand{\diag}{\mathrm{diag}}

\newcommand{\E}{\mathbb{E}}
\newcommand{\rE}{{\sf E}}
\newcommand{\dE}{{\sf E}^{+}}

\newcommand{\End}{\mathrm{End}}
\renewcommand{\epsilon}{\varepsilon}

\newcommand{\ext}{\mathchoice
	    {\raisebox{1pt}{$\bigwedge$}}
             {\raisebox{1pt}{$\bigwedge$}}
             {\raisebox{0.5pt}{$\scriptstyle\bigwedge$}}
             {\raisebox{0.2pt}{$\scriptscriptstyle\bigwedge$}}}
\newcommand{\F}{{\sf F}}
\newcommand{\cF}{\mathcal{F}}

\newcommand{\rG}{\mathsf{G}}
\newcommand{\Gr}{\mathsf{Gr}}

\newcommand{\hol}{\mathrm{hol}}

\newcommand{\Jx}{{\sf J}^{1}_{\ul{\x}}}
\newcommand{\Jxm}{{\sf J}_{\ul{\x}}}

\newcommand{\lanx}{\langle} 
\newcommand{\ranx}{\rangle}
\newcommand{\blanx}{\big\langle} 
\newcommand{\branx}{\big\rangle}

\newcommand{\cI}{\mathcal{I}}

\newcommand{\id}{\mathrm{id}}

\newcommand{\im}{\mathop{\mathrm{im}}}

\newcommand{\kZ}{Z}
\newcommand{\K}{\mathbb{K}}

\renewcommand{\k}{\mathsf{k}}

\newcommand{\oL}{\mathop{\ell}}

\newcommand{\rL}{{\sf L}}
\newcommand{\cL}{\mathcal{L}}

\newcommand{\cM}{\mathcal{M}}

\newcommand{\ol}{\overline}
\newcommand\op[1]{\mathsf{#1}}

\newcommand{\proj}[1]{{\sf\Pi}^{#1}}

\newcommand{\Q}{{\mathsf Q}}

\newcommand{\R}{\mathbb{R}}

\newcommand{\rank}{\mathrm{rank}}

\newcommand{\SU}{\mathrm{SU}}

\newcommand{\rB}{\mathsf{B}}
\newcommand{\rC}{\mathsf{C}}

\newcommand{\oR}{\partial}
\newcommand{\oD}{d}

\newcommand{\rT}{\mathsf{T}}
\newcommand{\cT}{\mathcal{T}}
\newcommand{\rP}{\mathsf{P}}

\newcommand{\T}{\mathsf{T}}

\newcommand{\cU}{\mathcal{U}}
\newcommand{\rU}{\mathsf{U}}
\newcommand{\ul}{\underline}

\newcommand{\rV}{{\sf V}}
\newcommand{\cV}{\mathcal{V}}

\newcommand{\Vect}{{\mathrm{Vect}}}

\newcommand{\X}{{\mathsf X}}
\newcommand{\Y}{{\mathsf Y}}

\renewcommand{\P}{\mathbb{P}}

\newcommand{\x}{x}
\newcommand{\y}{y}

\newcommand{\Z}{\mathbb{Z}}
\newcommand{\rZ}{{\sf Z}}
\newcommand{\cZ}{\mathscr{Z}}

\newcommand{\cH}{\mathcal{H}}

\renewcommand{\le}{\leqslant}
\renewcommand{\leq}{\leqslant}
\renewcommand{\ge}{\geqslant}
\renewcommand{\geq}{\geqslant}

\title{Determinantal random subgraphs}

\author{Adrien Kassel}
\address{Adrien Kassel -- CNRS -- UMPA, ENS de Lyon}
\email{adrien.kassel@ens-lyon.fr}

\author{Thierry L\'evy}
\address{Thierry L\'evy -- LPSM, Sorbonne Universit\'e, Paris}
\email{thierry.levy@sorbonne-universite.fr}

\date{\today}



\begin{document}

\maketitle

\begin{abstract} We define two families of determinantal random spanning subgraphs of a finite connected graph, one supported by acyclic spanning subgraphs (spanning forests) with fixed number of connected components, the other by connected spanning subgraphs with fixed number of independent cycles. Each family generalizes the uniform spanning tree and the generating functions of these probability measures generalize the  classical Kirchhoff and Symanzik polynomials. 

We call Symanzik spanning forests the elements of the acyclic spanning subgraphs family, and single out a particular determinantal mixture of these, having as kernel a normalized Laplacian on $1$-forms, which we call the Laplacian spanning forest.

Our proofs rely on a set of integral and real or complex (which we call geometric) multilinear identies involving cycles, coboundaries, and forests on graphs. We prove these identities using classical pieces of the algebraic topology of graphs and the exterior calculus applied to finite determinantal point processes, both of which we treat in a self-contained way.

We emphasize the matroidal nature of our constructions, thereby showing how the above two families of random spanning subgraphs are dual to one another, as well as possible generalisations. 
\end{abstract}

{\small
\setcounter{tocdepth}{1}
\tableofcontents
}


\section{Introduction}

\subsection{Determinantal spanning trees and cycle-rooted spanning forests} 

Since the seminal work of Kirchhoff~\cite{Kirchhoff}, the algebraic properties of spanning trees on finite connected graphs have kept fascinating and have continually been rediscovered in various guises using a variety of techniques. Spanning trees are ubiquitous in large areas of the literature in combinatorics, mathematical physics, probability, and linear algebra. Important foundational works of Whitney \cite{Whitney} and Tutte \cite{Tutte}, followed by many others, have shown how fundamental these objects are in combinatorics, and also how they can be seen in a broader context, notably that of matroids. This point of view percolated in probability theory, notably through the work of Lyons \cite{Lyons-DPP}.

It is known since the work of Burton and Pemantle \cite[transfer current theorem]{Burton-Pemantle} that the uniform probability measure on the set of spanning trees of a finite connected graph is a determinantal point process, a class of processes first introduced by Macchi \cite{Macchi} and named that way by Borodin and Olshanski at the turn of the century, see \cite{Borodin}. This measure had been studied earlier, in particular in relation to the Markov chain tree theorem (see \cite{Aldous} and references therein), and extended to infinite graphs in \cite{Pemantle,BLPS}; see the textbook \cite{Lyons-Peres}. In the planar case, the study of its scaling limit led Schramm to the discovery of SLE, see \cite{Schramm,LSW}. Analogs of uniform spanning trees on higher dimensional simplicial complexes were defined by Lyons \cite{Lyons-Betti}, who also highlighted why the support of a determinantal probability measure is the same thing as the set of bases of a linear matroid. Later, Kenyon \cite{Kenyon} defined a determinantal probability measure on the set of cycle-rooted spanning forests of a graph, determined by a $1$-form on the graph. There are quaternion determinant analogues of these probability measures~\cite{Kenyon,Kassel-ESAIM, KL3}.

On a graph, spanning trees and cycle-rooted spanning forests are the set of bases of the circular and bicircular matroids, respectively \cite{Oxley}. The circular case is the one from which the theory of matroids arose in the first place, whereas the bicircular case was only discovered later in~\cite{Simoes-Pereira} and further studied in \cite{Matthews}. These are moreover the only matroids on the set of edges of a graph for which the set of circuits consists in all subgraphs homeomorphic to a given family of connected graphs \cite{Simoes-Pereira}.

\subsection{Further determinantal spanning subgraphs}

\begin{figure}[!ht]
\centering
\includegraphics[width=6cm]{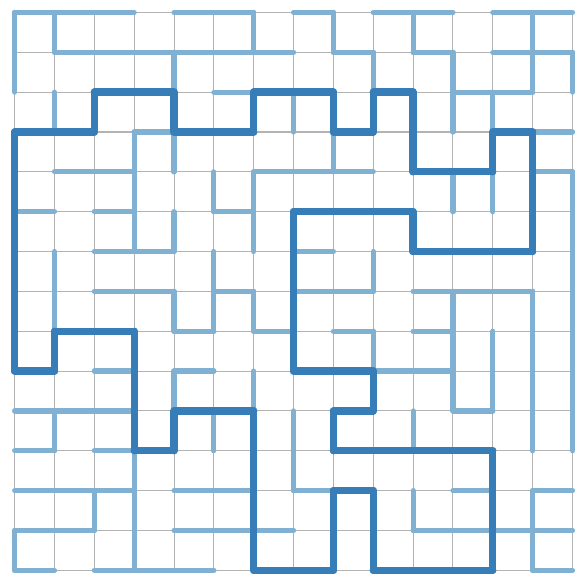} \hspace{1cm} \includegraphics[width=6cm]{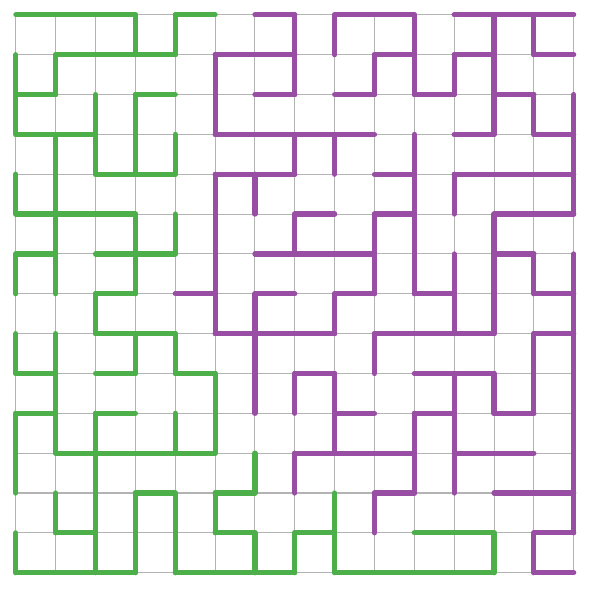}
\caption{\small A spanning unicycle (left) and a two-component spanning forest (right) of a $15\times 15$ square grid graph.}\label{fig:intro}
\end{figure}

In this paper, we define determinantal probability measures on spanning subgraphs with more complicated topology than trees (and cycle-rooted spanning forests), and give an explicit geometric formula for the weight of each graph appearing in the associated partition function. This is the content of Theorems~\ref{thm:DPP1} and~\ref{thm:DPP2}. These random subgraphs are respectively \emph{spanning forests} with fixed number of components and \emph{connected spanning subgraphs} with fixed number of independent cycles, see Figure \ref{fig:intro}.\footnote{We note that all figures presented in this paper have been obtained using the well-known algorithm for sampling determinantal point processes (see \cite{HKPV} or \cite[Proposition~3.13]{KL3}). Unless otherwise specified, we use uniformly random chosen parameters for the samples in this paper. Our basic sampling code is available on the repository
\url{https://plmlab.math.cnrs.fr/thierrylevy/determinantal-random-subgraphs}. As pointed out to us by one of the referees, it is possible to use faster sampling strategies since the determinantal point processes considered here have fixed cardinality, such as the one of \cite{Barthelme}. Letting $e$ and $v$ be the number of edges and vertices, this would yield a complexity of $O(e(v-1-k)^2)$ rather than $O(e^3)$ which in our examples would not change much since $e$ is of the order of $v$ and~$k$ is small. Moreover, we are unaware of a Wilson style algorithm~\cite{Wilson} for the probability measures we consider, but having one would naturally be of great use.} As will be apparent from \eqref{eq:symanzik-intro} and \eqref{eq:kirchhoff-intro}, let us note right away that these probability measures are \emph{not} `uniform'\footnote{Uniform in the sense of being proportional to the product of edge weights of each configuration.} unless we are considering the intersection of these families, namely the case of spanning trees.

The partition functions of these probability measures generalize the classical Symanzik and Kirchhoff polynomials (see Section~\ref{sec:polynomials} and in particular \eqref{eq:sym-pol-k}). Moreover, we generalize these results to the case of linear matroids in Theorem \ref{thm:matroid} and suggest some possible further specializations of this theorem in Section~\ref{sec:examples}.

Defining a determinantal random subgraph is equivalent to specifying a self-adjoint positive contraction on the space of $1$-forms of the graph. The families discussed above correspond to the orthogonal projections on subspaces contained in, or containing, the space of exact forms (see Section~\ref{sec:dss}).\footnote{The determinantal random spanning tree corresponds to the space where this subspace is the space of exact $1$-forms. Surprisingly, it also is the `uniform' measure, by the fundamental result of Burton and Pemantle. In that sense, these families generalize the `uniform' spanning tree case, although, as already said, the other measures we consider in this paper are generically not `uniform'.}
In the acyclic case (corresponding to subspaces contained in the space of exact forms), we further consider a certain operator canonically associated to the graph (a normalized Laplacian on $1$-forms) which is a convex combination of these (see Section \ref{sec:Laplacian-forest}) and we call the random spanning forest it induces, the Laplacian spanning forest (see Theorem~\ref{thm:laplacian-sf}).

\subsection{Main results}

We now provide a more detailed account of the content of Theorem~\ref{thm:DPP1} and Theorem~\ref{thm:DPP2}.

\subsubsection{The connected case}

Let us explain the content of Theorem~\ref{thm:DPP1}. On a weighted graph~$(\rG,\ul\x)$, given an integer $k\geq 0$, we consider the set $\cC_{k}(\rG)$ of connected spanning subgraphs with exactly~$k$ linearly independent cycles. Let us choose linearly independent $1$-forms $\theta_{1},\ldots,\theta_{k}$ which span a subspace that does not contain any non-zero exact $1$-form. For every $K\in \cC_{k}(\rG)$, we choose an integral basis $(c_{1},\ldots,c_{k})$ of the free abelian group of cycles of $K$ and assign to $K$ the weight
\begin{equation}\label{eq:symanzik-intro}
w(K)\, \ul{x}^{K}=\big|\det\big(\theta_{i}(c_{j})\big)_{1\leq i,j\leq k}\big|^{2} \ \ul\x^{K},\, 
\end{equation}
where $\theta(c)=\sum_{e\in c} \theta_{e}$ and $\ul\x^{K}$ is the product of the weights of the edges of $K$. We prove that the corresponding probability measure is determinantal, associated with the orthogonal projection on the direct sum of the space of exact $1$-forms and the span of $\vartheta=(\theta_{1},\ldots,\theta_{k})$:
\begin{equation}
H_{\vartheta}=\im d \oplus \Vect(\theta_1, \ldots, \theta_k)\,,
\end{equation}
where $d$ is the discrete derivative on the graph ($df(e)=f(\ol{e})-f(\ul{e})$ for all edges $e=(\ul{e},\ol{e})$ and $f\in \R^\rV$).

In particular, in the simplest case $k=1$, writing $\theta$ for $\theta_1$, the determinantal measure thus considered assigns to any spanning unicycle $K$ a weight proportional to $\theta(c)^2\ul{x}^K$. This probability measure was first considered by Kenyon in \cite{Kenyon} and arose as the limit of a sequence of determinantal probability measures supported on cycle-rooted spanning forests; we generalize this limit construction in~\cite{KL6}.

We can also compute the partition function of the measure \eqref{eq:symanzik-intro}. By Proposition \ref{prop:egalitesconnexes}, it satisfies
\begin{equation}\label{eq:partition-c}
\sum_{K\in \cC_k(\rG)} w(K)\,  \ul{x}^K = \rT(\ul{x}) \; \lVert \proj{\ker d^*}(\theta_1)\wedge\ldots\wedge\proj{\ker d^*}(\theta_k) \rVert^2\,.
\end{equation}
where $\T(\ul{x})=\sum_{T\in \cT(\rG)}  \ul{x}^T$, often called the Kirchhoff polynomial, is the generating function of spanning trees of $\rG$, the set of which we write $\cT(\rG)$.
In the case $k=1$, Equation \eqref{eq:partition-c} already appeared in \cite[Lemma 2]{KaWu} and \cite[Lemma 4]{Kassel-Kenyon}.

\subsubsection{The acyclic case}\label{sec:kirchhoff-intro}

Let us now explain the content of Theorem~\ref{thm:DPP2}, but in the reformulation given in Theorem \ref{thm:sym-forests}. 
Let $k\ge 1$ be an integer and consider a collection ${\sf q}=(q_1, \ldots, q_k)$ of $k$ linearly independent functions over vertices such that for each $i\in\{1,\ldots, k\}$, we have $\sum_{v\in \rV} q_i(v)=0$. Consider the set $\cF_{k}(\rG)$ of spanning forests with~$k+1$ connected components, and for each such forest $F=\{T_0,\ldots, T_{k}\}$, define the weight
\begin{equation}\label{eq:kirchhoff-intro}
w(F)\, \ul{x}^{F}=\big|\det\big(q_{i}(V_{j})\big)_{1\leq i,j\leq k}\big|^{2}\ \ul\x^{F}\,,
\end{equation}
where $V_j$ is the set of vertices of the tree $T_j$, and $q_i(V_j)=\sum_{v\in V_j} q_i(v)$. 
This weight does not depend on the choice of labelling of the $k+1$ trees constituting $F$, and the measure on $\cF_{k}(\rG)$ which assigns to each $F$ the weight $w(F)\, \ul{x}^{F}$ is non-zero, and when normalized, is a determinantal probability measure associated to a subspace of the space of~$1$-forms described explicitly by 
\begin{equation}
H_{\sf{q}}=\Big\{df: f\in \R^{\rV}, \sum_{v\in \rV} f(v) q_i(v)=0, \forall i\in \{1,\ldots, k\} \Big\}=d\Big( \big(\bigoplus_{i=1}^{d}\R q_i\big)^{\perp}\Big)\,.
\end{equation}

Let us emphasize that this random subgraph, which has a fixed total number of edges, is not the same thing as the often considered determinantal probability measure on spanning forests, without restriction on the number of components, called the {\em massive spanning forest}, which assigns to each spanning forest a weight proportional to the product of the total mass of each connected component.

In the simplest case where $k=1$, writing $q$ for $q_1$, the determinantal measure that we consider assigns to any two-component spanning forest (of which we denote by $T$ the component containing an arbitrary fixed vertex) a weight proportional to $q(\rV(T))^2$. The partition function of this measure is well known to physicists and is known as the first Symanzik polynomial; see~\cite{Bogner-Weinzierl,ABBGF} and the discussion in Section~\ref{sec:symanzik-kirchhoff-pols} below. For this reason, we call the weighted random spanning forests associated with the weight \eqref{eq:kirchhoff-intro} \emph{Symanzik spanning forests}. 

The partition functions of the Symanzik spanning forests feature some interesting polynomial identities (see Proposition \ref{prop:amini-k}). Given positive edge-weights $\underline{x}=(x_e)$, we have: 
\begin{equation}\label{eq:sym-green-intro}
\sum_{F\in \cF_{k}(\rG)}\big|\det\big(q_{i}(V_{j})\big)_{1\leq i,j\leq k}\big|^{2} \ \ul x^{F}= T_{\rG}(\ul{x}) \;\det\left[\big(\lanx q_i, G_{\ul{x}} q_j\ranx\big)_{1\le i,j\le k}\right]\,,
\end{equation}
where 
$T_{\rG}(\ul{x})=\sum_{T\in \cF_{1}(\rG)}  \ul x^{F}$ is the classical Kirchhoff polynomial of $\rG$, and where $G_{\ul{x}}$ is a certain operator on $\Omega^0(\rG)$ (a Green function, the inverse of the discrete Laplacian on the orthogonal of its kernel, see Section~\ref{sec:ratio-omid}).

All the subgraphs considered so far are determinantal with projection kernels. A normalized Laplacian on $1$-forms provides us with the kernel of a non-trivial mixture of Symanzik forests that we call the \emph{Laplacian spanning forest} (see Theorem \ref{thm:laplacian-sf}).

\subsubsection{The unifying matroidal case}

As the similarity between the above two displayed equations~\eqref{eq:symanzik-intro} and~\eqref{eq:kirchhoff-intro} shows, the two cases (acyclic and connected) are very similar. In fact, they are dual from the point of view of matroid theory. We prove in Theorem~\ref{thm:matroid} a common generalization of these two results to the case of a representable (or linear) matroid, whose representing map has a Euclidean target space, which induces a determinantal probability measure on its set of bases as shown by Lyons \cite{Lyons-DPP} (and reproven in a self-contained way in this paper).

We refer to Section \ref{sec:matroid-setting} for definitions and precise statements. The readers unfamiliar with matroid theory may wish to skip this part of the introduction at first.

\medskip

In a few words, the situation is the following. Let $S$ be a finite set and $E$ a linear space with basis $(e_i)_{i\in S}$. Let $\partial:E\to F$ is a linear map. This data defines a \emph{matroid} $\cM$ on $S$ whose bases are the maximal subsets $T$ of $S$ such that $\{\partial(e_{i}) : i\in T\}$ is linearly independent. Let $Z=\ker \partial$ be the subspace of $E$ encoding the `\emph{relations of dependency}' in this matroid. 

To define a determinantal probability measure on $S$ based on the above data, we assume that $E^{*}$ and $F^{*}$, the dual spaces of $E$ and $F$, are inner product spaces, in such a way that the basis $(e^{\star}_{i})_{i\in S}$ of $E^{*}$, dual to the basis $(e_{i})_{i\in S}$ of $E$, is orthogonal. In other words, we assume that there exists a collection of positive {\em Euclidean weights} $\ul{x}=(\x_i:i\in S)$ such that for all $i,j\in S$,
\begin{equation}
\lanx e_{i}^{\star}, e_{j}^{\star} \ranx = x_{i} \delta_{ij}.
\end{equation}

We let $d:F^*\to E^*$ be the transpose map of $\partial$. The determinantal probability measure on $2^{S}$ associated with the subspace $\im d$ has support equal to the set of bases of $\cM$. This is our \emph{reference determinantal probability measure} (akin to the uniform spanning tree measure in the above cases). We assume that we have some understanding of this reference distribution and in particular, we assign to each basis $T$ a weight $w(T)$, which we assume we know how to compute in an efficient way, and such that the probability of $T$ is proportional to~$w(T)\ul{x}^T$.


\smallskip

Now, for any subspace $H$ of $E$ containing $\im d$, and setting $k=\dim H - \dim \im d$, we define a new determinantal probability measure on $2^B$, whose support is included in the set $\B_{k}$ of subsets of $S$ that can be obtained by adding $k$ elements to a basis of $\cM$. 
The set of possible choices for the space $H$ is isomorphic to the Grassmannian of the quotient space $E/\im d$.

The gist of Theorem \ref{thm:matroid} is that we can describe explicitly the weight $K\mapsto w(K) \, \ul{x}^K$ of each configuration $K\in \B_k$.
Let $\theta_1,\ldots, \theta_k$ be linearly independent vectors in $E^*$ such that $H=\im d \oplus \Vect(\theta_1,\ldots, \theta_k)$. For each $K\in \B_k$ and each basis $T$ of $\cM$ included in~$K$, we have
\begin{equation}\label{eq:matroid-intro}
w(K)\, \ul{x}^K=w(T) \left\vert\det\left((\theta_i(z^K_{T,j})\right)_{1\le i,j\le k}\right\vert^2 \ul{x}^K\,,
\end{equation}
where $(z^K_{T,j}:1\le j \le k)$ is the \emph{fundamental basis} of $Z\cap E_K$ induced by $T$ 
(see Section \ref{sec:fund-bases-mat}).
Up to the computation of $w(T)$, the weight $w(K)$ only depends on a `small' $k\times k$ determinant (irrespective of the size of the ground set $S$) which is thus easy to compute.

Note the similarity of \eqref{eq:matroid-intro} with \eqref{eq:symanzik-intro}, as well as the difference, which comes from the fact that our reference measure is not uniform in general (in the sense that $w(T)$ depends on $T$).

As in \eqref{eq:partition-c} and \eqref{eq:sym-green-intro}, we can compute the partition function of the measure \eqref{eq:matroid-intro}. By Proposition~\ref{prop:egalitesmat}, it satisfies
\begin{equation}
\sum_{K\in \B_k} w(K)\, \ul{x}^K = \rB(\ul{x}) \; \lVert \proj{\ker d^*}(\theta_1)\wedge\ldots\wedge\proj{\ker d^*}(\theta_k) \rVert^2\,.
\end{equation}
where $\rB(\ul{x})=\sum_{T} w(T)\, \ul{x}^T$ is the weighted generating function of bases of $\cM$, that is, the partition function of the reference determinantal measure.

\subsection{Method of proof}

We briefly describe the methods used to prove Theorems \ref{thm:DPP1} and \ref{thm:DPP2}, and their generalization in Theorem \ref{thm:matroid}. We focus on Theorem \ref{thm:DPP1}, as it is simpler and illustrates the overall approach.

This result comes with a key companion result which is Proposition \ref{prop:egalitesconnexes} (or Proposition \ref{prop:egalitesacycliques} or Proposition \ref{prop:egalitesmat}).

We present two strategies, each offering complementary perspectives. Neither one is clearly preferable to the other, so we include both. One approach is based on computing the \emph{partition function}, the other on computing the \emph{weight of individual configurations}. 

In both cases we need to compute the square of a volume, but in the first method we compute a \emph{global} volume depending on the whole graph $\rG$ (and then partition it, to see the probability distribution on spanning subgraphs $K\in \cC_k(\rG)$ appear), and in the second method we compute a \emph{local} volume depending only on a spanning subgraph (and then piece these `local' volumes together to recover the partition function for the whole graph). Section \ref{sec:spanningtrees} previews both methods by rederiving the determinantal structure of the uniform spanning tree in two concise ways (see Proposition~\ref{prop:lawUST}).

Both approaches to proving Theorem \ref{thm:DPP1} also rely on an additional common idea: the new probability measures can be viewed as perturbations of a reference measure—specifically, the uniform spanning tree distribution. At this point, the Schur formula, in the guise of Lemma~\ref{lem:Schurcomp}, becomes instrumental.

Recall the statement of Theorem \ref{thm:DPP1} we want to prove: the determinantal measure associated with $H=H_0\oplus \Vect(\theta_1, \ldots, \theta_k)$ where $H_0=\im d$ is supported on $K\in C_k(\rG)$ assigning each $K$ a weight proportional to $|\det(\theta_i(c_j))|^2$, where $(c_1,\ldots, c_k)$ is an integral basis of the group of cycles of~$K$. 

This result comes with a companion result, Proposition~\ref{prop:egalitesconnexes} which states that the ratio of the new partition function to the reference one, has a nice expression as the squared norm of the projection, in an exterior power of the space of $1$-forms, of $\vartheta=\theta_1\wedge\ldots\wedge\theta_k$ on $(H_0^\perp)^{\wedge k}$.

In the first method, we prove Theorem \ref{thm:DPP1} as a corollary of Proposition \ref{prop:egalitesconnexes}, and in the second method, we prove Proposition \ref{prop:egalitesconnexes} as a corollary of Theorem \ref{thm:DPP1}.

\subsubsection{Global-to-local method: perturbation of the partition function}

This method interprets the new partition function as a perturbation of the classical Laplacian determinant $\det'(d^*d)$.

\begin{enumerate}
\item We define an operator $d_\vartheta=d+\omega_{\vartheta}$ whose image is $H=H_0 \oplus \Vect(\theta_1,\ldots, \theta_k)$. By writing the squared volume of the image of an orthonormal basis under $d_\vartheta$ as the product of two terms (using Lemma \ref{lem:Schurcomp}), one of which involves only $d$, we obtain the ratio of this volume to the reference volume associated with $d$ as the squared norm of the orthogonal projection of ~$\vartheta$ on $(H_0^{\perp})^{\wedge k}$ . This proves Proposition~\ref{prop:egalitesconnexes}.
\item Combining this computation with the mean projection theorem (Theorem \ref{thm:projection}), we re-express this projection, and thus the partition function, as a sum over pairs of a spanning tree (basis of our reference matroid) and a set of $k$ edges. The multilinear identity (Proposition \ref{prop:sum2-sym}) then simplifies this to a sum over configurations $K \in \cC_k(\rG)$, with weights matching those expected via an application of Lemma~\ref{lem:Schurcomp}. This proves Theorem~\ref{thm:DPP1}.
\end{enumerate}

\subsubsection{Local-to-global method: perturbation of configuration weights}

This method starts with a configuration $K \in \cC_k(\rG)$ and computes its weight directly. From the general theory of determinantal probability measures, we know that this weight is the squared cosine of the angle between the subspace spanned by $1$-forms supported on $K$ and the space $H = H_0 \oplus \mathrm{Vect}(\theta_1, \ldots, \theta_k)$

\begin{enumerate}
\item By picking appropriate bases of the spaces involved, and performing elementary row and column operations, we show that the probability of $K\in \cC_k(\rG)$ associated with $H=H_0\oplus\Vect(\theta_1,\ldots, \theta_k)$ is proportional to the weight $\vert\det(\theta_i(c_j))_{1\le i,j\le k}\vert^2$. This proves Theorem~\ref{thm:DPP1}.
\item The sum of all these configuration weights, that is, the partition function, is then shown to be the determinant of $(d + \omega)^*(d + \omega)$, recovering the same projection norm of $\vartheta$ as in the first method for the ratio of this partition function to the reference partition function. This proves Proposition~\ref{prop:egalitesconnexes}.
\end{enumerate}

In the graphical case, simplifications occur due to the fact that invertible integer matrices have determinant equal to $\pm 1$. Without these simplifications in the general matroidal case, we need to carefully take care of new determinantal quantities.

\subsection{Toolbox}

Our main tools for proving these results are an exterior algebra version of the matrix-tree theorem (Propositions~\ref{lem:omid},~\ref{lem:dimo}, and~\ref{lem:fundamental-global}), and consequences of it, and the mean projection theorem for determinantal point processes (Theorem \ref{thm:projection}). The latter theorem was proven in~\cite[Theorem~5.9]{KL3} and another proof was given in \cite{KL4}; earlier instances of special cases of this statement appeared in \cite{Nerode-Shank}, \cite[Theorem 1]{Maurer}, \cite[Proposition~7.3]{Biggs-graph}, \cite[Proposition~6.8]{Lyons-DPP}, and \cite[Theorem A]{CCK-line, CCK-complexes}. As pointed out to us by one of the referees, a form of this theorem also appears in \cite[Theorems 2 and 3]{Derezinski}.

Our approach allows us to unify the presentation of several statements concerning spanning trees (Sections~\ref{sec:graphs} and~\ref{sec:trees}), the Jacobian torus (Section \ref{sec:sandpile}), duality (Section \ref{sec:planar-duality}) and complexes (Section \ref{sec:cellulation}), the Kirchhoff and Symanzik polynomials (Section \ref{sec:polynomials}), determinantal probability measures (Sections \ref{sec:det-toolbox} and \ref{sec:matroid-setting}), cycle-rooted spanning forests (Section \ref{sec:examples}) and matroids (Section~\ref{sec:matroid-setting}). Incidentally, it also yields a new formula for the probability density of a determinantal process (Proposition \ref{prop:proba-det}) and its restrictions (Corollary~\ref{cor:PTXK}), in addition to the description of the above-mentioned families of examples of determinantal random subgraphs. 

In that respect, this paper also provides yet another, almost self-contained, presentation of these classical topics, expressed in the unifying language of the exterior algebra. 

\subsection{Organization of the paper}

The paper is organized as follows. 
In Section~\ref{sec:graphs} we introduce some basic definitions about graphs and associated objects (chains, cochains, cycles, coboundaries, and integral bases determined by spanning trees). 
In Section~\ref{sec:trees}, we prove combinatorial multilinear identities involving spanning trees. 
Section \ref{sec:metricstructures} is devoted to the introduction of metric structures on cochains.
In Section~\ref{sec:det-toolbox}, we introduce four independent tools from the theory of determinantal probability measures and determinant computations. 
Section \ref{sec:spanningtrees} reviews, in a very short and self-contained way, the fundamental example of the uniform spanning tree. 
In Section \ref{sec:gmi} we prove geometric multilinear identities involving spanning trees. 
In Section~\ref{sec:dss} we show our main result, namely the existence of noteworthy 
determinantal probability measures on two families of subgraphs (spanning forests and spanning connected graphs with constrained Betti numbers), 
which generalize the uniform measure on spanning trees. These two families are dual in a matroidal sense, and in the case where the graph is embedded in a two-dimensional surface, this duality becomes topological, as explained briefly in Section \ref{sec:two-dimensional-case}.
In Section~\ref{sec:polynomials}, we make the connection to stable multivariate homogeneous polynomials from analytic combinatorics and theoretical physics and derive a few consequences. In particular, we define a generalization of Symanzik polynomials. 
In Section \ref{sec:symanzik-forests} we reparametrize our family of determinantal spanning forests, and show how they are related to our generalization of Symanzik polynomials, explain the link to the classical model of massive spanning forests, and single out a special model of determinantal random forest, which we call the Laplacian random forest. 
In Section \ref{sec:connected-subgraphs} we discuss a few special cases of the determinantal connected spanning subgraphs. 
In Section~\ref{sec:matroid-setting} we explain how the results presented extend from the circular matroid case to the general case of linear matroids, which in particular encompasses the case of the bicircular matroid, or `circular' and `bicircular' matroids on the set of cells of higher dimensional simplicial complexes.
The paper closes with further examples in Section \ref{sec:examples} and concluding remarks  in Section~\ref{sec:conclusion}. 

\subsection*{Acknowledgements}
We thank Omid Amini for inspiring discussions during work visits in Paris and Lyon on the topic of Symanzik polynomials and related structures. In particular, we realized while completing this work, which was motivated by different considerations in \cite{KL6}, that the natural generalization of Symanzik polynomials we encounter here (see \eqref{eq:sym-pol-k} in Section~\ref{sec:polynomials}) had already been imagined by him several years ago, in the guise of the determinantal expression in the right-hand side of Proposition \ref{prop:egalitesacycliques}, based on the abstract construction of these polynomials in \cite[Section 2.1]{ABBGF}. This paper thus also provides an answer to the question of Omid Amini of providing a concrete description and some properties of these polynomials. We also thank Javier Fres\'an for conversations at ETH Zurich which introduced us to \cite{ABBGF} back in 2015.

In addition, we thank the four anonymous referees for their feedback and in particular Referee X for a very thorough report which helped us improve the quality of the paper.


\section{Spanning trees, cycles, and coboundaries}\label{sec:graphs}

In this section, we present the notions from the algebraic topology of graphs that we will need for our study of determinantal subgraphs. Most of the material presented in the next pages is well known, or even classical (see for instance \cite{Biggs-book, Biggs-graph,Nagnibeda}). However, the detail of definitions and terminology vary from one author to the other, and to keep our exposition self-contained, we take the time to introduce ours.

\subsection{Graphs, subgraphs, orientations, ordering}

The geometric setting of this paper is a {\em finite connected graph} $\rG=(\rV,\rE)$, that we choose now and will keep fixed (unless explicitly stated otherwise) until the end of Section \ref{sec:polynomials}.
The set of {\em vertices} $\rV$ and the set of {\em edges} $\rE$ are finite. Edges are oriented and there are two maps $e\mapsto {\sf{s}}(e)$ and $e\mapsto {\sf{t}}(e)$ from $\rE$ to $\rV$, which to an edge associate its source and target, that is, its starting and ending vertices. Moreover, each edge comes with both possible orientations: there is an orientation-reversing fixed-point-free involution $e\mapsto e^{-1}$ of $\rE$, which exchanges starting and ending vertices. 

We will deal with many kinds of subgraphs of $\rG$, that will all be {\em spanning subgraphs}, in the sense that their vertex set is $\rV$. For us, a subgraph of $\rG$ will be a subset $S\subseteq \rE$ that is stable under orientation reversal, which we regard as the graph $(\rV,S)$. A subgraph of $\rG$ needs not be connected, and can have isolated vertices, which we count as connected components. The complement of a subgraph $S$ is the subgraph $S^{c}=\rE\setminus S$.

A prominent role will be played by the class of {\em spanning trees} of $\rG$, which are the minimal connected subgraphs of $\rG$, and also the maximal acyclic subgraphs of $\rG$. We denote by $\cT(\rG)$ the set of all spanning trees of $\rG$. Two families of classes of subgraphs, which depend on a non-negative integer $k$, will also play an important role in this paper, namely the class $\cC_{k}(\rG)$ of subgraphs obtained by adding $k$ edges to a spanning tree, and the class $\cF_{k}(\rG)$ of subgraphs obtained by removing $k$ edges from a spanning tree. Of course, $\cF_{0}(\rG)=\cT(\rG)=\cC_{0}(\rG)$. We will give alternative descriptions of these classes in due time (see Section \ref{sec:betti}).

Although nothing we are going to do ultimately depends on it, we need to choose an orientation of each edge of $\rG$. Thus, we fix once and for all a subset $\dE$ of $\rE$ containing exactly one element of each pair $\{e,e^{-1}\}$. Given a subgraph $S$ of $\rG$, we set $S^{+}=S\cap \rE^{+}$.

We also pick once and for all a total ordering of~$\dE$. This will allow us to write matrices indexed by subsets of $\dE$ without ambiguity on the ordering of rows or columns. Also, anticipating a bit, exterior products over sets of oriented edges will always be taken in this order. 

\subsection{Chains and cochains}\label{sec:chain-modules}

Let us now recall the classical notions of chains and cochains on a graph, their natural pairing, and the boundary and coboundary maps induced by the graph structure. Loosely speaking, chains are linear combinations of vertices, or edges, and cochains are functions on vertices, or edges. The least classical part of what follows is our discussion of relative coboundaries of the pair $(\rG,S)$  formed by our graph and a subgraph of it.

\subsubsection{Chains and cycles} Let us start with chains.
Let us denote by $C_{0}(\rG)$ the free $\Z$-module over the set $\rV$ of vertices of the graph $\rG$, and by $C_{1}(\rG)$ the quotient of the free $\Z$-module over~$\rE$ by the submodule generated by $\{e+e^{-1}:e\in \rE\}$. The classes of the elements of $\dE$ form a basis of $C_{1}(\rG)$, that we call the canonical basis, and identify with the set $\dE$.

The \emph{boundary} operator $\partial:C_{1}(\rG)\to C_{0}(\rG)$ is defined by setting, for each edge $e$,
\[\partial e={\sf{t}}(e)-{\sf{s}}(e).\]
The range of $\partial$ is contained in the kernel of the augmentation morphism $\epsilon : C_{0}(\rG)\to \Z$, which is defined by setting $\epsilon(v)=1$ for every vertex $v$. Since $\rG$ is connected, the range of $\partial$ is in fact equal to the kernel of $\epsilon$.

We define the group of {\em cycles} as
\[Z_{1}(\rG)=\ker \partial\]
and we have the exact sequence
\[0 \longrightarrow Z_{1}(\rG) \longrightarrow C_{1}(\rG) \build{\longrightarrow}{}{\partial} C_{0}(\rG)\build{\longrightarrow}{}{\epsilon} \Z \longrightarrow 0 .\]

\subsubsection{Cochains and coboundaries}\label{sec:cochains-coboundaries}
We now consider a dual construction. The groups of cochains are defined by 
\[C^{0}(\rG)={\rm Hom}(C_{0}(\rG),\Z) \ \ \text{ and } \ \ C^{1}(\rG)={\rm Hom}(C_{1}(\rG),\Z).\]
We denote by $(\1_{v})_{v\in \rV}$ the canonical basis of $C^{0}(\rG)$ and by $(e^{\star})_{e\in \dE}$ the canonical basis of~$C^{1}(\rG)$. For example, we have $e^{\star}(e^{-1})=-1$. 

The pairing between chains and cochains is denoted by round brackets: given a chain $a$ and a cochain $\alpha$ of the same degree, we set $(\alpha,a)=\alpha(a)$.

The \emph{coboundary} operator $\delta:C^{0}(\rG)\to C^{1}(\rG)$ is defined as the adjoint of the boundary operator. Thus, for every vertex $v\in \rV$, we have
\[\delta \1_{v} = \sum_{e\in \rE : {\sf{t}}(e)=v} e^{\star}.\]
The kernel of $\delta$ is equal to the range of the map $\Z\to C^{0}(\rG)$, to which we do not give a name, which sends $1$ to the augmentation morphism $\epsilon$.

We define the group of {\em coboundaries}\footnote{Some authors call the group $B^{1}(\rG)$ the {\em lattice of integral cuts} of the graph $\rG$, see for example \cite{Nagnibeda}.} as
\[B^{1}(\rG)=\im \delta.\]
The situation is now described by the exact sequence 
\[0 \longrightarrow \Z \longrightarrow C^{0}(\rG)  \build{\longrightarrow}{}{\delta} C^{1}(\rG)  \longrightarrow C^{1}(\rG)/B^{1}(\rG) \longrightarrow 0.\]
Let us notice that $B^{1}(\rG)$, the range of $\delta$, is the annihilator of the kernel of $\partial$, that is, the annihilator of $Z_{1}(\rG)$. We will  write this duality as $B^{1}(\rG)=Z_{1}(\rG)^{\circ}=\{\alpha\in C^{1}(\rG): \forall z\in Z_{1}(\rG), (\alpha,z)=0\}$.
\medskip

\subsubsection{Subgraphs} Let us discuss cycles and coboundaries on subgraphs of $\rG$, keeping in mind that, in this paper, all graphs are spanning, in the sense that their vertex set is $\rV$. Thus, let us pick a subgraph $(\rV,S)$ of $\rG$, that we simply refer to as $S$. There is no need for a special discussion of chains and cochains of degree $0$ on $S$, since they are exactly the same as on $\rG$. In degree~$1$ however, the group of chains on $S$ identifies naturally with the submodule $C_{1}(S)\subseteq C_{1}(\rG)$ generated by the edges of $S$, and the group of cochains on $S$ with the quotient $C^{1}(\rG)/C_{1}(S)^{\circ}$, which in turn is isomorphic with the submodule $C^{1}(S)\subseteq C^{1}(\rG)$ generated by $\{e^{\star} :e\in S\}$.

There are natural splittings $C_{1}(\rG)=C_{1}(S)\oplus C_{1}(S^{c})$ and $C^{1}(\rG)=C^{1}(S)\oplus C^{1}(S^{c})$, and we will denote by
$\pi_{S}:C_{1}(\rG)\to C_{1}(S)$ and $\pi_{S}:C^{1}(\rG)\to C ^{1}(S)$ the  corresponding projections.

The space $Z_{1}(S)$ of cycles on $S$ is the kernel of the boundary operator restricted to $C_{1}(S)$:
\begin{equation}\label{eq:defZ1S}
Z_{1}(S)=Z_{1}(\rG)\cap C_{1}(S) \subseteq C_{1}(\rG).
\end{equation}
The space $B^{1}(S)$ of coboundaries on $S$ is the range of the coboundary operator `co-restricted' to~$C^{1}(S)$, that is, the range of the operator $\pi_{S}\circ \delta$:
\[B^{1}(S)=\pi_{S}\big(B^{1}(\rG)\big) \subseteq C^{1}(\rG).\] 
Let us emphasize that $Z_{1}(S)$ is a submodule of $Z_{1}(\rG)$, but $B^{1}(S)$ is {\em not} a submodule of $B^{1}(\rG)$ in general. 
Instead,  it is the quotient of $B^{1}(\rG)$ by its intersection with the kernel of the map $\pi_{S}:C^{1}(\rG)\to C^{1}(S)$.
This intersection will play an important role, and we give it a name, setting 
\begin{equation}\label{eq:defB1S}
B^{1}(\rG,S)=B^{1}(\rG)\cap C^{1}(S^{c})\subseteq C^{1}(\rG),
\end{equation}
by close analogy with \eqref{eq:defZ1S}. The submodule $B^{1}(\rG,S)$ of $B^{1}(\rG)$ is the set of coboundaries of $0$-chains that are locally constant on $S$, or more simply the set of coboundaries on $\rG$ that vanish on $S$. In the classical terminology of algebraic topology, it would be called the group of {\em relative coboundaries} of the pair $(\rG,S)$. It will play for us the role of the group of coboundaries on the graph $\rG_{/S}$ obtained by contracting all the edges of $S$, but we prefer not to consider this quotient graph and to keep working on $\rG$. 

\subsubsection{Betti numbers}\label{sec:betti}
Let us conclude by mentioning Betti numbers and the Euler relation. To the subgraph~$S$ are associated four non-negative integers, namely 
\begin{itemize}
\item $c_{0}$, the number of vertices, equal to $|\rV|$,
\item $c_{1}$, the number of edges, equal to $|S^{+}|$,
\item $b_{0}$, the number of connected components of $S$,
\item $b_{1}$, the rank of $Z_{1}(S)$.
\end{itemize}
The numbers $b_{0}$ and $b_{1}$ are called the Betti numbers of $S$, and the fundamental Euler relation
\begin{equation}\label{eq:Euler}
c_{0}-c_{1}=b_{0}-b_{1}
\end{equation}
holds. The rank of $B^{1}(S)$ is $c_{0}-b_{0}$ and that of $B^{1}(\rG,S)$ is $b_{0}-1$. 

Spanning trees are characterized by the fact that they are connected and acyclic, that is, by the Betti numbers $b_{0}=1$ and $b_{1}=0$. The other classes of subgraphs that we will consider are
\[\cC_{k}(\rG)=\{\text{subgraphs with } b_{0}=1, c_{0}=k\} \text{ and } \cF_{k}(\rG)=\{\text{subgraphs with } b_{0}=k+1, c_{0}=0\},\]
respectively the classes of connected spanning subgraphs and acyclic spanning subgraphs. 

\subsection{Splittings and integral bases induced by spanning trees}

In this section, we explain how each spanning tree of $\rG$ induces on the one hand a basis of $Z_{1}(\rG)$ and a splitting of $C_{1}(\rG)$, and on the other hand a basis of $B^{1}(\rG)$ and a splitting of $C^{1}(\rG)$. We also discuss the case of subgraphs. The bases we refer to are bases of~$\Z$-modules and we call them integral bases to distinguish them from bases of the vector spaces that we will consider later. 

\subsubsection{Chains on $\rG$}
Let $T$ be a spanning tree of $\rG$. Let $C_{1}(T)$ be the submodule of $C_{1}(\rG)$ generated by the edges of $T$. A fundamental fact is that the restriction $\partial : C_{1}(T)\to C_{0}(\rG)$ is 
injective (because~$T$ is acyclic) and has the same range as $\partial : C_{1}(\rG)\to C_{0}(\rG)$, namely the kernel of the augmentation morphism (because $T$ is spanning and connected). Therefore, any element of $C_{1}(\rG)$ has the same boundary as a unique element of $C_{1}(T)$, from which it differs by a boundaryless $1$-chain, that is, a cycle. In other words, there is a splitting
\begin{equation}\label{eq:decC1Z}
C_{1}(\rG)=C_{1}(T)\oplus Z_{1}(\rG)
\end{equation}
and we denote by 
\begin{equation}\label{eq:defZT}
\rZ_{T}:C_{1}(\rG)\to Z_{1}(\rG)
\end{equation}
the corresponding projection. 

To give a more concrete expression of $\rZ_{T}$, given two vertices $v,w$ of $\rG$, let us denote by $[v,w]_{T}$ the unique simple path going from $v$ to $w$ in $T$, seen as an element of $C_{1}(\rG)$. This is the unique element of $C_{1}(T)$ such that $\partial( [v,w]_{T})=w-v$. Then for each edge $e$, we have $\partial e=-\partial ([{\sf{t}}(e),{\sf{s}}(e)]_{T})$, so that the abovementioned splitting reads, in this case, $e=-[{\sf{t}}(e),{\sf{s}}(e)]_{T}+(e+[{\sf{t}}(e),{\sf{s}}(e)]_{T})$. Therefore,
\begin{equation}\label{eq:defZT2}
\rZ_{T}(e)=e+[{\sf{t}}(e),{\sf{s}}(e)]_{T}.
\end{equation}

Let us now use $\rZ_{T}$ to construct a basis of $Z_{1}(\rG)$. Since $C_{1}(\rG)=C_{1}(T)\oplus C_{1}(T^{c})$, there is an induced isomorphism $\rZ_{T} :C_{1}(T^{c})\simeq C_{1}(\rG)/C_{1}(T)\to Z_{1}(\rG)$, and we find that the family 
\[\mathscr Z_{T}=\big\{\rZ_{T}(e) : e\in  (T^{c})^{+}\big\}\]
is a basis of $Z_{1}(\rG)$ (see also, for instance, \cite[Theorem 5.2]{Biggs-book}).

From the definition of $\rZ_{T}$ (or equivalently from \eqref{eq:defZT2}), we observe that $\rZ_{T}(e)=0$ whenever $e \in T$. Thus, we have $\rZ_{T}\circ \pi_{T}=0$, so that on $Z_{1}(\rG)$, and with a slight abuse of notation, the following equality holds:
\begin{equation}
{\rZ}_{T}\circ \pi_{T^{c}}={\rm id}_{Z_{1}(\rG)}. \label{eq:zpiid}
\end{equation}
To be more precise, the composed map
\[Z_{1}(\rG)\hookrightarrow C_{1}(\rG) \build{\longrightarrow}{}{\pi_{T^{c}}} C_{1}(T^{c}) \build{\longrightarrow}{}{\rZ_{T}} Z_{1}(\rG)\]
is the identity.

\subsubsection{Chains on connected spanning subgraphs} Let $K$ be a connected spanning subgraph of $\rG$, thus an element of $\cC_{k}(\rG)$ for $k$ equal to the rank of $Z_{1}(K)$. The discussion above can be repeated in the submodule $C_{1}(K)$ of $C_{1}(\rG)$. Spanning trees of~$K$ are exactly the spanning trees of $\rG$ that are contained in~$K$. If $T$ is such a spanning tree, then the map $\rZ_{T}$ sends the submodule $C_{1}(K)$ of $C_{1}(\rG)$ in $Z_{1}(K)$. The map $\rZ_{T}:C_{1}(K)\to Z_{1}(K)$ is onto, because it is the identity on the submodule $Z_{1}(K)$ of $C_{1}(K)$. The kernel of this map is $C^{1}(T)$. Therefore, this map induces an isomorphism $C_{1}(K)/C_{1}(T)\simeq C_{1}(K\setminus T)\to Z_{1}(K)$ and the family
\begin{equation}\label{eq:defZKT}
\mathscr Z_{T}^{K}=\big\{\rZ_{T}(e) : e\in  (K\setminus T)^{+}\big\}
\end{equation}
is a basis of $Z_{1}(G)$. Moreover, the composed map
\[Z_{1}(K)\hookrightarrow C_{1}(K) \build{\longrightarrow}{}{\pi_{K\setminus T}} C_{1}(K\setminus T) \build{\longrightarrow}{}{\rZ_{T}} Z_{1}(K)\]
is the identity of $Z_{1}(K)$, which we summarize by saying that 
\begin{equation}\label{eq:zpiidlocal}
\rZ_{T}\circ \pi_{K\setminus T} = {\rm id}\ \text{  on the submodule } \ Z_{1}(K) \text{ of } C_{1}(\rG).
\end{equation}

\subsubsection{Cochains on $\rG$} The situation for cochains is analogous to the situation for chains that we just described, with small differences. Taking \eqref{eq:decC1Z} and denoting annihilators with a circle, we find the splitting
\begin{equation}\label{eq:decC1B}
C^{1}(\rG)=Z_{1}(\rG)^{\circ} \oplus C_{1}(T)^{\circ} = B^{1}(\rG) \oplus C^{1}(T^{c}).
\end{equation}
We denote by 
\begin{equation}\label{eq:defBT}
\rB_{T}:C^{1}(\rG)\to B^{1}(\rG)
\end{equation}
the corresponding projection on the first summand. In English, every element of $C^{1}(\rG)$ coincides on~$T$ with (that is, differs by an element of $C^{1}(T^{c})$ from) a unique coboundary on $\rG$. Indeed, consider $\alpha\in C^{1}(\rG)$. Choose a vertex $v_{0}\in \rV$ and, for each vertex $v$, define $\phi(v)=(\alpha,[v_{0},v]_{T})$. Then $\phi$ is an element of $C^{0}(\rG)$ such that $\delta \phi$ and $\alpha$ coincide on $T$, and $\rB_{T}(\alpha)=\delta \phi$.

In the case where $\alpha=e^{\star}$ for some edge $e\in \rE$, then after choosing $v_{0}={\sf{s}}(e)$,  the function $\phi$ can be described as the indicator of the set $U\subseteq \rV$ of vertices of $\rG$ that are connected to ${\sf{t}}(e)$ in $T\setminus e$, and $\rB_{T}(e^{\star})=\delta \1_{U}$. In particular, if $e$ does not belong to $T$, then $\rB_{T}(e^{\star})=0$. 

The projection $\rB_{T}$ induces an isomorphism $\rB_{T}:C^{1}(\rG)/C^{1}(T^{c})\simeq C^{1}(T)\to B^{1}(\rG)$, so that the family 
\[\mathscr B_{T}=\big\{\rB_{T}(e^{\star}) : e\in T^{+}\big\}\]
is a basis of $B^{1}(\rG)$ (see also \cite{Nagnibeda}).

Moreover, in the same vein as \eqref{eq:zpiid}, we have the equalities $\rB_{T}\circ \pi_{T^{c}}=0$ and, on $B^{1}(\rG)$,
\begin{equation}
{\rB}_{T}\circ \pi_{T}={\rm id}_{B^{1}(\rG)}\,, \label{eq:bpiid}
\end{equation}
the precise meaning of this statement being that the composed map
\[B^{1}(\rG)\hookrightarrow C^{1}(\rG) \build{\longrightarrow}{}{\pi_{T}} C^{1}(T) \build{\longrightarrow}{}{\rB_{T}} B^{1}(\rG)\]
is the identity.

\subsubsection{Cochains on acyclic spanning subgraphs} As for chains, but in a dual way, this discussion can be localised to a spanning subgraph that we now assume to be acyclic, that is, a spanning forest. Let us choose $F\in \cF_{k}(\rG)$ for some integer $k\geq 0$. Let $T$ be a spanning tree of $\rG$ {\em containing} $F$. It follows from its definition that the map $\rB_{T}$ sends the submodule $C^{1}(F^{c})$ of $C^{1}(\rG)$ in $B^{1}(\rG,F)$, and the restricted map $\rB_{T}:C^{1}(F^{c})\to B^{1}(\rG,F)$ is surjective, because it is the identity on the submodule $B^{1}(\rG,F)$ of $C^{1}(F^{c})$. The kernel of this restricted map is $C^{1}(T^{c})$, and $\rB_{T}$ induces an isomorphism $\rB_{T}:C^{1}(F^{c})/C^{1}(T^{c})\simeq C^{1}(T\setminus F)\to B^{1}(\rG,F)$. Therefore, the family 
\begin{equation}\label{eq:defBFT}
\mathscr B_{T}^{F}=\big\{\rB_{T}(e^{\star}) : e\in (T\setminus F)^{+}\big\}
\end{equation}
is a basis of $B^{1}(\rG,F)$. Moreover, the map $\rB_{T}\circ \pi_{T\setminus F}$ restricts to the identity on the submodule $B^{1}(\rG,F)$ of $C^{1}(\rG)$, in the sense that
\begin{equation}\label{eq:bpiidlocal}
\rB_{T}\circ \pi_{T\setminus F} = {\rm id}\ \text{  on the submodule } \ B^{1}(\rG,F) \text{ of } C^{1}(\rG).
\end{equation}

\subsubsection{A simple identity} To conclude this section, let us state and prove a simple identity.

\begin{proposition} \label{prop:decBZ} Let $T$ be a spanning tree of $\rG$. For all $\alpha\in C^{1}(\rG)$ and $c\in C_{1}(\rG)$, one has
\[(\alpha,c)=(\rB_{T} \alpha,c)+(\alpha,\rZ_{T}c).\]
\end{proposition}

\begin{proof} The difference between the two sides is equal to
\[(\rB_{T} \alpha,c)+(\alpha,\rZ_{T}c)-(\alpha,c)=(\rB_{T}\alpha,\rZ_{T}c)-\big((\alpha-\rB_{T}\alpha),(c-\rZ_{T}c)\big)\]
and the two terms on the right-hand side are zero, the first as the pairing of a coboundary with a cycle, the second as the pairing of an element of $C^{1}(T^{c})$ with an element of $C_{1}(T)$.
\end{proof}


\section{Integral multilinear identities}\label{sec:trees}\label{sec:integralbases}

In this section, we will prove several identities involving on the one hand trees and cycles and on the other hand trees and coboundaries in a graph. These identities take place in exterior powers of groups of integral chains or cochains. For general definitions and notations about exterior calculus, we refer to \cite[Section 5]{KL3}. 

In Section \ref{sec:gmi}, we will prove identities that are superficially similar, but depend on the choice of inner products on spaces of cochains, that we will introduce in Section \ref{sec:metricstructures}. By contrast, the identities that we prove in the present section are of purely topological nature.  

Given an integral basis $\mathscr B=(u_{1},\ldots,u_{k})$ of a free $\Z$-module $M$ of rank $k$, we will denote by $\det \mathscr B$ the element $u_{1}\wedge \ldots \wedge u_{k}$ of the top exterior power $M^{\wedge k}$ of $M$. If $M$ is a submodule of another module $N$, then we see $\det \mathscr B$ as an element of $N^{\wedge k}$.

For any two bases $\mathscr B_{1}$ and $\mathscr B_{2}$ of $M$, the elements $\det \mathscr B_{1}$ and $\det \mathscr B_{2}$ are equal up to a sign, that we denote by $\det(\mathscr B_{1}/\mathscr B_{2})=\pm 1$, and which is the determinant of the matrix that expresses the elements of $\mathscr B_{1}$ in the basis $\mathscr B_{2}$.

In order for determinants to have a definite sign, it is necessary that bases be ordered. The bases $\mathscr Z_{T}^{K}$ of $Z_{1}(K)$ and $\mathscr B_{T}^{F}$ of $B^{1}(\rG,F)$ defined in the previous section (see \eqref{eq:defZKT} and \eqref{eq:defBFT}) will always be endowed with the total order inherited from that of $\dE$, fixed once and for all at the beginning.

We denote the exterior algebra of a $\Z$-module $M$ by $\ext M$. In $\ext \, C_{1}(\rG)$ and $\ext \, C^{1}(\rG)$, we will use the following notation : for every subset $S$ of $\rE$, with $S\cap \dE=\{e_{1},\ldots,e_{r}\}$  enumerated in the order of $\dE$, we write
\begin{equation}\label{eq:defeS}
e_{S}=e_{1}\wedge \ldots \wedge e_{r} \ \text{ and } \ e^{\star}_{S}=e^{\star}_{1}\wedge \ldots \wedge e^{\star}_{r}.
\end{equation}

\subsection{Connected spanning subgraphs}\label{sec:symanzik-tree}

The content of this section is somewhat related to a result of \cite{Baker}, but our statement and proof is different, and elementary.

Let us choose an integer $k\geq 0$ and a subgraph $K\in \cC_{k}(\rG)$ of $\rG$, that is, a connected spanning subgraph of $\rG$ with $k$ linearly independent cycles, in the sense that the $\Z$-module $Z_{1}(K)$ has rank~$k$. Let us fix an arbitrary integral basis $\mathscr Z^{K}$ of $Z_{1}(K)$. Our two main identities give an expression of the element $\det \mathscr Z^{K}$ of $\ext^{k}C_{1}(\rG)$, and an expression of its tensor square.

\begin{proposition}[First cycle-tree identity]\label{lem:omid}
In $\ext^{k} C_{1}(\rG)$, the following identity holds:
\begin{equation}\label{eq:omid}
\det \mathscr Z^{K}=\sum_{T\in\cT(\rG), T\subseteq K} \det(\mathscr Z^{K}/\mathscr Z^{K}_{T}) \, e_{K\setminus T}.
\end{equation}
\end{proposition}

We will use several times the following consequence of the Euler relation \eqref{eq:Euler}: if a subgraph of $\rG$  has the same number of edges as a spanning tree, that is, $|\rV|-1$, then its Betti numbers satisfy $b_{0}-b_{1}=1$, so that if it is not a spanning tree itself, it has at least two connected components, and at least one non-trivial cycle.

\begin{proof} Let us write $\mathscr Z^{K}=(z_{1},\ldots,z_{k})$ and  decompose the element $z_{1}\wedge \ldots \wedge z_{k}$ of $\ext^{k}C_{1}(\rG)$ on the basis $\{e_{S}: |S|=k\}$:
\[z_{1}\wedge \ldots \wedge z_{k}=\sum_{S\subset \rE^{+} : |S|=k} a_{S} e_{S},\]
for some integral coefficients $a_{S}$. Since $z_{1},\ldots,z_{k}$ belong to the submodule $C^{1}(K)$ of $C^1(\rG)$, the sum involves only basis elements~$e_{S}$ with $S\subseteq K$. Consider now a subgraph $S$ of $K$ with $k$ edges and assume that $K\setminus S$ is not a spanning tree of $\rG$. Then~$K\setminus S$ is a subgraph with $|\rV|-1$ edges that is not a spanning tree, so that it contains a non-trivial cycle. Let us call $y$ such a cycle and decompose it in the integral basis $(z_{1},\ldots,z_{k})$ of $Z_{1}(K)$ as $y=n_{1}z_{1}+\ldots+n_{k}z_{k}$. By reordering $(z_{1},\ldots,z_{k})$ if needed, let us make sure that~$n_{1}\neq 0$. Then $0=(\pi_{S})^{\wedge k}(y\wedge z_{2}\wedge \ldots \wedge z_{k})=n_{1}(\pi_{S})^{\wedge k}(z_{1}\wedge \ldots \wedge z_{k})=n_{1} a_{S} e_{S}$, so that~$a_{S}=0$.

Consider now a spanning tree $T$ of $\rG$ such that $T\subseteq K$. Using \eqref{eq:zpiidlocal}, we find that 
\[z_{1}\wedge \ldots \wedge z_{k}=(\rZ_{T}\circ \pi_{K\setminus T})^{\wedge k}(z_{1}\wedge \ldots \wedge z_{k})=\sum_{S\subset \rE^{+} : |S|=k} a_{S}\, (\rZ_{T}) ^{\wedge k } \big((\pi_{K\setminus T})^{\wedge k}(e_{S})\big)\]
and the only non-zero term of the last sum is that corresponding to $S=K\setminus T$, so that
\[\det \mathcal Z^{K}=z_{1}\wedge \ldots \wedge z_{k}=a_{K\setminus T}\, (\rZ_{T})^{\wedge k}e_{K\setminus T}.\]
By the definition  of $\mathcal Z^{K}_{T}$ (see \eqref{eq:defZKT}), $(\rZ_{T})^{\wedge k}e_{K\setminus T}=\det \mathcal Z^{K}_{T}$. This identifies the coefficient $a_{K\setminus T}$ as $\det(\cZ^{K}/\cZ^{K}_{T})$ and concludes the proof.
\end{proof}

For every spanning tree $T$ of $\rG$ contained in $K$, applying $(\pi_{K\setminus T})^{\wedge k}$ to \eqref{eq:omid} yields
\begin{equation}\label{eq:pitcgamma}
(\pi_{K\setminus T})^{\wedge k}(\det \mathscr Z^{K})=\det(\cZ^{K}/\cZ^{K}_{T}) \, e_{K\setminus T}.
\end{equation}
Moreover, since $\ext^{k}Z_{1}(K)$ is a module of rank~$1$ generated by $\det \mathscr Z^{K}$, the last equation combined with \eqref{eq:omid} implies that on the submodule $\ext^{k} Z_{1}(K)$ of $\ext^{k}C_{1}(\rG)$,
\begin{equation}\label{eq:equalityextZ}
\sum_{T\in \cT(\rG), T\subseteq K} (\pi_{K\setminus T})^{\wedge k}={\rm id}. 
\end{equation}

\begin{proposition}[Second cycle-tree identity] \label{prop:sum2-sym}
In $\big(\ext^{k}C_{1}(\rG)\big)^{\otimes 2}$, the following equality holds:
\[(\det \mathscr Z^{K})^{\otimes 2}=\sum_{\substack{T\in \cT(\rG),T\subseteq K}} ({\rZ}_{T})^{\wedge k} (e_{K\setminus T}) \otimes e_{K\setminus T}.\]
\end{proposition}

\begin{proof} 
Let us apply \eqref{eq:equalityextZ} to the second factor of the left-hand side, and then the $k$-th exterior power of \eqref{eq:zpiidlocal} to the first factor of each term of the sum, to find
\begin{align*}
(\det \mathscr Z)^{\otimes 2}&=\sum_{T\in \cT(\rG):T\subseteq K} \det \mathscr Z \otimes (\pi_{K\setminus T})^{\wedge k}(\det \mathscr Z)\\
&=\sum_{T\in \cT(\rG):T\subseteq K} ({\rZ}_{T})^{\wedge k}\circ (\pi_{K\setminus T})^{\wedge k}(\det \mathscr Z) \otimes (\pi_{K\setminus T})^{\wedge k}(\det \mathscr Z).
\end{align*}
For each term of the sum, \eqref{eq:pitcgamma} yields $(\pi_{K\setminus T})^{\wedge k}(\det \mathscr Z)=\pm e_{K\setminus T}$ and the result follows.
\end{proof}

\subsection{Acyclic spanning subgraphs}\label{sec:kirchhoff-tree}

Let us now choose an integer $k\geq 0$ and a subgraph $F\in \cF_{k}(\rG)$, that is, a spanning forest with $k+1$ connected components. The module $B^{1}(\rG,F)$ has rank~$k$ and we fix an integral basis $\mathscr B^{F}$ of it. We will give an expression of $\det \mathscr B^{F}\in \ext^{k}C^{1}(\rG)$, and of its tensor square. 

\begin{proposition}[First coboundary-tree identity]\label{lem:dimo}  
In $\ext^{k} C^{1}(\rG)$, the following identity holds:
\begin{equation}\label{eq:dimo}
\det \mathscr B^{F}=\sum_{T\in\cT(\rG), T\supseteq F} \det(\mathscr B^{F}/\mathscr B^{F}_{T})\, e^{\star}_{T\setminus F} \; . 
\end{equation}
\end{proposition}

\begin{proof} Let us write $\mathscr B^{F}=(\kappa_{1},\ldots,\kappa_{k})$ and decompose $\kappa_{1}\wedge \ldots \wedge \kappa_{k}\in \ext^{k}C^{1}(\rG)$ on the basis $\{e^{\star}_{S}: |S|=k\}$:
\[\kappa_{1}\wedge \ldots \wedge \kappa_{k}=\sum_{S\subset \rE^{+} : |S|=k} b_{S} e^{\star}_{S}.\]
Since $B^{1}(\rG,F)$ is a submodule of $C^{1}(F^{c})$, the sum involves only terms for which $S$ is disjoint from~$F$. Consider now a subgraph $S$ with $k$ edges, disjoint from $F$, and such that $F\cup S$ is not a spanning tree of $\rG$. Then the subgraph $F\cup S$ has $|\rV|-1$ edges but it is not a spanning tree, so that it has at least two connected components. Let $\alpha$ be the coboundary of the indicator of one of them. Since every connected component of $F$ is contained in a connected component of $F\cup S$, the $1$-cochain~$\alpha$ vanishes on $F$. In other words, $\alpha$ belongs to $B^{1}(\rG,F)$, of which it is a non-zero element. 
Let us decompose it on the integral basis $\kappa_{1},\ldots,\kappa_{k}$, as $\alpha=n_{1}\kappa_{1}+\ldots +n_{k}\kappa_{k}$.
Reordering $\kappa_{1},\ldots,\kappa_{k}$ if needed, we may assume that $n_{1}\neq 0$. Then, since $\alpha$ vanishes on $S$, we have 
\[0=(\pi_{S})^{\wedge k}(\alpha\wedge \kappa_{2}\wedge \ldots \wedge \kappa_{k})=n_{1}(\pi_{S})^{\wedge k}( \kappa_{1}\wedge \ldots \wedge \kappa_{k})=n_{1} b_{S}e^{\star}_{S},\]
so that $b_{S}=0$.

Consider now a spanning tree $T$ of $\rG$ such that $T\supseteq F$. Using \eqref{eq:bpiid}, we find that 
\[\kappa_{1}\wedge \ldots \wedge \kappa_{k}=(\rB_{T}\circ \pi_{T\setminus F})^{\wedge k}( \kappa_{1}\wedge \ldots \wedge \kappa_{k})=\sum_{S\subset \rE^{+} : |S|=k} b_{S}\, (\rB_{T}\circ \pi_{T\setminus F})^{\wedge k}(e^{\star}_{S})\]
and the only non-zero term of the last sum is that corresponding to $S=T\setminus F$, so that
\[\det \mathscr B^{F}=\kappa_{1}\wedge \ldots \wedge \kappa_{k}=b_{T}\, (\rB_{T})^{\wedge k}e^{\star}_{T\setminus F}.\]
From the definition of $\mathscr B^{F}_{T}$ (see \eqref{eq:defBFT}), we see that $(\rB_{T})^{\wedge k}e^{\star}_{T\setminus F}=\det\mathscr B^{F}_{T}$. This identifies the coefficient $b_{T}$ as $\det(\cB^{F}/\cB^{F}_{T})$ and concludes the proof.
\end{proof}

For every spanning tree $T$ containing $F$, applying $(\pi_{T\setminus F})^{\wedge k}$ to \eqref{eq:dimo}, we find
\begin{equation}\label{eq:pitkappa}
(\pi_{T\setminus F})^{\wedge k}(\det \mathscr B^{F})=\det(\cB^{F}/\cB^{F}_{T}) \, e^{\star}_{T}.
\end{equation}
Then, since $\ext^{k}B^{1}(\rG,F)$ is a module of rank~$1$ generated by $\det \mathscr B^{F}$, 
the last identity combined with~\eqref{eq:dimo} implies that on the submodule $\ext^{k}B^{1}(\rG,F)$ of $\ext^{k}C^{1}(\rG)$,
\begin{equation}\label{eq:equalityextB}
\sum_{T\in \cT(\rG), T\supseteq F}(\pi_{T\setminus F})^{\wedge k}={\rm id}.
\end{equation}

\begin{proposition}[Second coboundary-tree identity]\label{prop:sum2-kir}
In $\big(\ext^{k}C^{1}(\rG)\big)^{\otimes 2}$, one has
\[(\det \mathscr B^{F})^{\otimes 2}=\sum_{\substack{T\in \cT(\rG), T\supseteq F}}(\rB_{T})^{\wedge k} (e^{\star}_{T\setminus F}) \otimes e^{\star}_{T\setminus F}.\]
\end{proposition}

\begin{proof} 
Let us apply \eqref{eq:equalityextB} to the second factor of the left-hand side, and then the $k$-th exterior power of \eqref{eq:bpiidlocal} to the first factor of each term of the sum, to find
\begin{align*}
(\det \mathscr B^{F})^{\otimes 2}&=\sum_{T\in \cT(\rG), T\supseteq F}
\det \mathscr B^{F} \otimes (\pi_{T\setminus F})^{\wedge k}(\det \mathscr B^{F})\\
&=\sum_{T\in \cT(\rG), T\supseteq F} ({\rB}_{T})^{\wedge k} \circ (\pi_{T\setminus F})^{\wedge k}(\det \mathscr B^{F}) \otimes (\pi_{T\setminus F})^{\wedge k}(\det \mathscr B^{F}).
\end{align*}
In each term of the sum, \eqref{eq:pitkappa} implies that $(\pi_{T\setminus F})^{\wedge k}(\det \mathscr B^{F})=\pm e^{\star}_{T\setminus F}$ and the result follows.
\end{proof}


\section{Metric structures on cochains}\label{sec:metricstructures}

The probability measures on subgraphs of $\rG$ that we will consider in Section \ref{sec:dss} depend on more than the topological structure of the graph $\rG$: they depend on a certain amount of geometric structure, that is embodied by the choice of a collection $\ul{\x}=(x_{e})_{e\in \rE}$ of positive real weights indexed by the edges of the graph $\rG$, such that $x_{e^{-1}}=x_{e}$. These weights can be interpreted as inverse lengths or, viewing $\rG$ as an electrical network, as conductances. These weights will allow us to endow the spaces of cochains with a Euclidean or Hermitian structure. 

Let us choose a base field $\K$ that is equal to $\R$ or $\C$. We will consider chains and cochains with coefficients in $\K$, and to denote the corresponding spaces we will use the letter $\Omega$ instead of $C$. Thus, for $i\in \{0,1\}$, we set
\[\Omega^{i}(\,\cdot\,)=C^{i}(\,\cdot\,) \otimes \K \ \text{ and }\ \Omega_{i}(\,\cdot\,)=C_{i}(\,\cdot\,) \otimes \K.\]
We will use cochains with coefficients in $\K$ more often than chains, and instead of {\em cochains} we will often call them {\em forms} or, to be more precise, $0$-forms and $1$-forms. 

\subsection{Inner products}\label{sec:euclidean}

Let us endow $\Omega^{0}(\rG)$ with the inner product 
\begin{equation}\label{eq:inner-product0}
\lanx f,g\ranx=\sum_{v\in \rV} \overline{f(v)} g(v)
\end{equation} 
and $\Omega^{1}(\rG)$ with the inner product 
\begin{equation}\label{eq:inner-product1}
\lanx\alpha,\beta\ranx=\sum_{e\in\dE} \x_e \, \overline{\alpha(e)} \beta(e).
\end{equation} 
This definition is independent of the choice of orientation of the edges. 

Note that the inner product \eqref{eq:inner-product1} depends on $\ul \x$, whereas \eqref{eq:inner-product0} does not. In the following, we will not stress this dependence more explicitly, but it plays a key role in several proofs, where identities between polynomials in $\ul\x$ are considered, see in particular Sections \ref{sec:routine}, \ref{sec:polynomials}, and \ref{sec:matroid-setting}.
It would be possible, but not very useful for us, to endow the spaces of chains $\Omega_{0}(\rG)$ and $\Omega_{1}(\rG)$ with inner products. See Section~\ref{sec:conventions} for a brief discussion of this point.

Every $1$-chain defines a linear form on $\Omega^{1}(\rG)$, that can be represented, thanks to the inner product on this space, by an element of $\Omega^{1}(\rG)$ itself. We denote by 
\begin{equation}\label{eq:defJx1}
{\Jx}: \Omega_1(\rG) \to \Omega^{1}(\rG), \ e \mapsto x_{e}^{-1} e^{\star}
\end{equation}
this antilinear isomorphism. A similar but simpler antilinear isomorphism 
\begin{equation}\label{eq:defJ0}
{\sf J}^{0}: \Omega_0(\rG) \to \Omega^{0}(\rG), \ v \mapsto \1_{v}
\end{equation}
exists, that does not depend on the conductances. For all $v\in \rV$, $e\in \rE$, $f\in \Omega^{0}(\rG)$ and $\alpha\in \Omega^{1}(\rG)$, the following relations hold:
\begin{equation}\label{eq:J01}
\lanx {\sf J}^{0} v, f\ranx=(f,v)=f(v) \ \ \text{and} \ \ \lanx{\Jx} e,\alpha \ranx=(\alpha,e)=\alpha(e).
\end{equation}

In the following, we will use adjoints of linear maps between inner product spaces. We will denote by $u^{*}$ the adjoint of $u$, keeping in mind that it may depend on the choice of conductances.

For instance, we define 
\[d=\delta \otimes {\rm id}_{\K}:\Omega^{0}(\rG) \to \Omega^{1}(\rG) \ \text{ and } \ d^{*}:\Omega^{1}(\rG)\to \Omega^{0}(\rG),\]
the usual discrete differential and its adjoint. The operators $\partial$ and $d^{*}$ are related by the equation 
\begin{equation}\label{eq:partiald*}
{\sf J}^{0} \circ (\partial \otimes {\rm id}_{\K}) = d^{*}\circ {\Jx}.
\end{equation}
Indeed, for all $e\in \rE$, $f\in \Omega^{0}(\rG)$, we have $\lanx d^{*}{\Jx}e,f\ranx=\lanx {\Jx}e, df \ranx=(df,e)=(f,\partial e)=\lanx f,{\sf J}^{0} \partial e\ranx$.

\subsection{Splittings and projections}\label{sec:splitproject}

We have the orthogonal decomposition
\begin{equation}
\Omega^1(\rG)=\im d \oplus \ker d^* = (B^{1}(\rG)\otimes \K) \oplus {\Jx}(Z_{1}(\rG)\otimes \K),
\end{equation}
the equality ${\Jx}(Z_{1}(\rG)\otimes \K)=\ker d^{*}$ coming from \eqref{eq:partiald*}. This equality shows how the inner product allowed us to bring cycles and coboundaries together in a single space. 

Let $T\in \cT(\rG)$ be a spanning tree. We have ${\Jx}(\Omega_{1}(T))=\Omega^{1}(T)$, so that starting from~\eqref{eq:decC1Z}, tensoring by $\K$ and applying ${\Jx}$, we find the splitting
\begin{equation}\label{eq:splitOmega1d*T}
\Omega^1(\rG)=\Omega^1(T)\oplus \ker d^*,
\end{equation}
which is not orthogonal, unless $\rG$ is a tree. We denote by $\rP^{\ker d^*}_{T}$ the associated projection on $\ker d^{*}$. This projection is of course related to the map $\rZ_{T}$: we have 
\begin{equation}\label{eq:PJZJ}
\rP^{\ker d^{*}}_{T} ={\Jx} (\rZ_{T}\otimes {\rm id}_{\K})({\Jx})^{-1}.
\end{equation}

In a similar way, tensoring \eqref{eq:decC1B} by $\K$ yields the splitting
\begin{equation}\label{eq:splitOmega1dTc}
\Omega^1(\rG)=\im d \oplus \Omega^1(T^{c}),
\end{equation}
which again is not orthogonal, unless~$\rG$ is a tree. Note however that $\Omega^{1}(T^{c})$ is the orthogonal of~$\Omega^{1}(T)$, regardless of the choice of $\ul{\x}$. We will denote by $\rP^{\im d}_{T^{c}}$ the projection on $\im d$ corresponding to the splitting \eqref{eq:splitOmega1dTc}. It is related to $\rB_{T}$ through the simple equality
\begin{equation}\label{eq:PB}
\rP^{\im d}_{T^{c}}=\rB_{T}\otimes {\rm id}_{\K}.
\end{equation}

\subsection{Exterior powers}

The inner products that we defined on spaces of forms induce inner products on exterior powers of these spaces (see e.g.~\cite[Section 5.1]{KL3} for details): tensors of distinct degrees are orthogonal, and the inner product of two pure tensors of the same degree $k$ is given by 
\begin{equation} 
\lanx u_1\wedge \ldots \wedge u_k, v_1\wedge \ldots \wedge v_k\ranx = \det\left[\big(\lanx u_i,v_j \ranx\big)_{1\le i,j\le k}\right]\,.
\end{equation}
We will often use this equality and express determinants as scalar products in an exterior algebra.

For any collection $A\subseteq\rE$ of edges, let us define 
\[\ul{\x}^A=\prod_{e\in (A\cup A^{-1})\cap \rE^+} \x_e.\]

For every $k\geq 0$, the family $\{e^{\star}_{S} : S\subseteq \dE, |S|=k\}$ is an orthogonal basis of $\ext^{k}\Omega^{1}(\rG)$, with 
\begin{equation} \label{eq:monomials}
\|e^{\star}_{S}\|^{2}=\ul\x^{S}
\end{equation}
for every $S\subseteq \dE$,

\subsection{Inner products on chains}\label{sec:conventions}
We chose not to endow $\Omega_{0}(\rG)$ nor $\Omega_{1}(\rG)$ with an inner product.
If we had to, we would take the one which turns the maps ${\Jx}$ and ${\sf J}^0$ into isometries. On~$\Omega_0(\rG)$, this would be the canonical inner product, since there is no dependency on $\ul\x$ there. On $\Omega_1(\rG)$, this would be the inner product for which the canonical basis is orthogonal and for which $\|e\|^{2}=x_{e}^{-1}$ for every edge $e$. 
Under this inner product, for every $k\geq 0$, the family $\{e_{S} : S\subseteq \dE, |S|=k\}$ is orthogonal in $\ext^{k}\Omega_{1}(\rG)$, and for each $S\subseteq \dE$,
\begin{equation}
\|e_{S}\|^{2}=(\ul\x^{S})^{-1}.
\end{equation}
Since we think of $\ul\x$ as conductances (or inverse lengths), their inverses are resistances (or lengths), and this justifies the following invented rule of thumb: chains resist, cochains conduct.

\subsection{Kirchhoff polynomial and the classical matrix-tree theorem}\label{sec:pythagore-kir}

Let us now take a Euclidean (or Hermitian) look at some of the identities that we proved in Section \ref{sec:trees}. Let us start with the first coboundary-tree identity \eqref{eq:dimo} in the case where $F=\varnothing$. In this case, the module $B^{1}(\rG,F)$ is equal to $B^{1}(\rG)$ and has rank $k=|\rV|-1$, the number of edges of a spanning tree. Given any integral basis $\mathscr B$ of $B^{1}(\rG)$, the identity reads
\begin{equation}\label{eq:decomposedetB}
\det \mathscr B=\sum_{T\in \cT(\rG)} \det(\mathscr B/\mathscr B_{T}) \, e^{\star}_{T}.
\end{equation}
Read in the exterior algebra of $\Omega^{1}(\rG)$, the right-hand side is a sum of pairwise orthogonal terms, so that taking the square of the norm of both sides, we find
\begin{equation}\label{eq:kir0}
\|\det \mathscr B \, \|^{2}=\sum_{T\in \cT(\rG)} \ul\x^{T}.
\end{equation}
The right-hand side of \eqref{eq:kir0} is famously called the \emph{Kirchhoff polynomial} of $\rG$ 

To recover the classical matrix-tree theorem (see for instance \cite{KL2,KL2C}), let us apply \eqref{eq:kir0} to a particular integral basis of $B^{1}(\rG)$. For this, let us pick a vertex $v_{0}\in \rV$. Then the family 
\begin{equation}\label{eq:defBv0}
\mathscr B_{v_{0}}=\big\{\delta \1_{v} : v\in \rV\setminus \{v_{0}\}\big\}
\end{equation}
is an integral basis of $B^{1}(\rG)$ and
\[\|\det \mathscr B_{v_{0}}  \|^{2}=\det\big((\lanx d \1_{v},d \1_{w} \ranx)_{v,w\in \rV\setminus \{v_{0}\}} \big)=\det\big((\lanx \1_{v},d^{*}d \1_{w} \ranx)_{v,w\in \rV\setminus \{v_{0}\}} \big).
\]
Let us define the combinatorial Laplacian on $\Omega^{0}(\rG)$ as the operator
\[\Delta=d^{*}d.\]
Then the determinant that we just computed is that of the principal submatrix of $\Delta$ in the canonical basis of $\Omega^{0}(\rG)$, where the row and column corresponding to $v_{0}$ have been erased. Therefore, with our choice of basis, and a self-explanatory notation, we find that \eqref{eq:kir0} can be rewritten as
\begin{equation}\label{eq:mtt-classical}
\det \Delta^{\widehat{v_{0}}}_{\widehat{v_{0}}}=\sum_{T\in\cT(\rG)} \ul\x^{T}\,,
\end{equation}
which is one of the classical forms of the matrix-tree theorem.

\subsection{Symanzik polynomial}\label{sec:pythagore-sym}

Let us repeat what we just did, now with the first cycle-tree identity~\eqref{eq:omid}, when $K=\rG$. In this case, $Z_{1}(\rG)$ has rank $b_{1}$, the first Betti number of $\rG$, and for any integral basis $\mathscr Z$ of this module, the identity reads
\begin{equation}\label{eq:decomposedetZ}
\det \mathscr Z=\sum_{T\in \cT(\rG)} \det(\mathscr Z/\mathscr Z_{T}) \, e_{T^{c}}.
\end{equation}
Reading this equality in the tensor algebra of $\Omega_{1}(\rG)$, applying $(\Jx)^{\wedge b_{1}}$ to both sides and taking the square norm, we find\footnote{This is one of the rare places where it would have been useful to endow the spaces of chains with inner products.}
\begin{equation}\label{eq:sym1}
\|(\Jx)^{\wedge b_{1}} \det \mathscr Z \, \|^{2}=\sum_{T\in \cT(\rG)} \big(\ul\x^{\rE \setminus T}\big)^{-1}=\big(\ul\x^{\rE}\big)^{-1} \sum_{T\in \cT(\rG)} \; \ul\x^T\ .
\end{equation}

In contrast with the situation for coboundaries, there is no obvious integral basis of $Z_{1}(\rG)$. The best we can do is to choose one, say $\mathscr Z=(z_{1},\ldots,z_{b_{1}})$, and to describe the left-hand side of \eqref{eq:sym1} as
\[\|(\Jx)^{\wedge b_{1}} \det \mathscr Z \, \|^{2}=\det\big(\lanx \Jx z_{i},\Jx z_{j}\ranx_{i,j\in \{1,\ldots,b_{1}\}}\big).\]
Then \eqref{eq:sym1} can alternatively be phrased as follows: let $X=\mathrm{diag}(\x_e : e\in\dE)$ be the diagonal matrix of weights $\x_e$ and $M$ the $|\dE|\times b_{1}$ matrix formed by writing the cycles $z_{1},\ldots,z_{b_{1}}$ in the canonical basis of $C_{1}(\rG)$. Then
\begin{equation}\label{eq:sym1-omid}
\det ({}^{t}\!M X^{-1} M) =  \sum_{T\in\cT(\rG)} \; \big(\ul\x^{\rE\setminus T}\big)^{-1}\, , 
\end{equation}
an equality which appears in \cite[Lemma 3.1]{Amini-exchange}.
 
In the theory of Feynman integrals of Euclidean quantum field theory and associated graph polynomials \cite{Bogner-Weinzierl}, the right-hand side of~\eqref{eq:sym1-omid} is called the \emph{first Symanzik polynomial} (applied here to $\ul\x^{-1}$, see Section~\ref{sec:symanzik-kirchhoff-pols}).

\subsection{Real tori and finite abelian groups}\label{sec:sandpile}

Let us conclude this section by a short digression which is not needed for the rest of the paper. For simplicity, let us assume in this section that $\K=\R$.

Recall Equations \eqref{eq:kir0} and \eqref{eq:sym1}, and let us use the same notations. Using \eqref{eq:kir0} and the fact that ${\Jx}(Z_{1}(\rG)\otimes \R)=\ker d^{*}$ is orthogonal to $B^{1}(\rG)\otimes \R=\im d$, we find
\begin{align*}
\big\|\big((\Jx)^{\wedge b_{1}}\det \mathscr Z\big) \wedge \det \mathscr B\, \big\|^{2}&=
\big\|(\Jx)^{\wedge b_{1}}\det \mathscr Z\big\|^{2}\ \big\|\det \mathscr B\, \big\|^{2}
=\big(\ul\x^{\rE}\big)^{-1} \, \Big(\sum_{T\in\cT(\rG)} \; \ul{\x}^{T}\Big)^{2}\,.
\end{align*}

On the other hand, applying $(\Jx)^{\wedge b_{1}}$ to the left-hand side of \eqref{eq:omid}, and using \eqref{eq:dimo}, we find
\[\big((\Jx)^{\wedge b_{1}}\det \mathscr Z\big) \wedge \det \mathscr B
=\sum_{T\in\cT(\rG)} \pm \big(\ul\x^{T^{c}}\big)^{-1} \, e^{\star}_{T^{c}}\wedge e^{\star}_{T}=\big(\ul\x^{\rE}\big)^{-1} \Big(\sum_{T\in\cT(\rG)} \pm \, \ul\x^{T} \Big) e^{\star}_{\rE}\,,\]
with the notation introduced in \eqref{eq:defeS}. 

Comparing with the previous equality, we deduce that all signs in the last sum are the same.\footnote{
The fact that all signs are the same in the above formula is also a reflection of the fact that ${\Jx}^{\wedge b_{1}}\det \cZ$ and $\det \cB$ are Hodge dual of each other, up to a constant.}
Thus, we  proved the following proposition.

\begin{proposition} Let $\mathscr Z$ be a basis of $Z_{1}(\rG)$ and $\mathscr B$ a basis of $B^{1}(\rG)$. Then in the line $\ext^{|\dE|}\Omega^{1}(\rG)$, we have
\[\big((\Jx)^{\wedge b_{1}}\det \mathscr Z\big) \wedge \det \mathscr B=\pm \big(\ul\x^{\rE}\big)^{-1} \Big(\sum_{T\in\cT(\rG)} \ul\x^{T} \Big)\, e^{\star}_{\rE}\,.\] 
In particular, 
\[\big\|\big((\Jx)^{\wedge b_{1}}\det \mathscr Z\big) \wedge \det \mathscr B\, \big\|=\big(\ul\x^{\rE}\big)^{-\frac{1}{2}} \, \sum_{T\in\cT(\rG)} \; \ul{\x}^{T}.\]
\end{proposition}

In the vector space $\Omega^{1}(\rG)$, the images by $\Jx$ of the elements of $\mathscr Z$, together with the elements of~$\mathscr B$,
form a basis, and generate a lattice. This lattice, as a discrete abelian subgroup of $\Omega^{1}(\rG)$, does not depend on our choice of basis, indeed it is equal to ${\Jx}(Z_{1}(\rG))\oplus B^{1}(\rG)$. The quotient 
\[\Omega^{1}(\rG)/\big({\Jx}(Z_{1}(\rG))\oplus B^{1}(\rG)\big)\]
is a real torus, of which the second assertion of the last proposition computes the volume. In the case where all the weights $x_e$ are taken to be equal to $1$, this volume is equal to the number of spanning trees of $\rG$.

Still in the case where $\ul\x$ is identically equal to $\ul{1}$, this volume is also equal to the cardinal of the finite group
\[C^{1}(\rG)/\big({\sf J}^{1}_{\ul1}(Z_{1}(\rG))\oplus B^{1}(\rG)\big).\]
Taking $\ul\x$ identically equal to $\ul{1}$ blurs the distinction between chains and cochains and encourages us to identify them. If we do so, we can write the last group as
\[C_{1}(\rG)/\big(Z_{1}(\rG)\oplus \delta C_{0}(\rG)\big).\]
The boundary map $\partial$ descends to an injective map on this quotient, and induces an isomorphism with the group  
\[\partial C_{1}(\rG)\big/\partial\delta C_{0}(\rG),\] 
sometimes called the {\em Jacobian} group of the graph, itself isomorphic to the \emph{sandpile group}
\[C_{0}(\rG\setminus \{v_{0}\})\big/\partial\delta C_{0}(\rG\setminus \{v_{0}\}),\]
where $v_{0}$ is an arbitrarily chosen vertex and $\rG\setminus \{v_{0}\}$ is the graph obtained from $\rG$ by removing the vertex $v_{0}$ and all incident edges, see \cite[Corollary 13.15]{Corry-Perkinson} and \cite{Nagnibeda, Biggs-sand}. Kotani and Sunada also define the Jacobian torus \cite{Kotani-Sunada}, see also \cite{Baker}.


\section{Determinantal toolbox}\label{sec:det-toolbox}

We will now let probability enter the picture, under the form of determinantal probability measures. In this section, we give a self-contained presentation of a few basic facts about determinantal probability measures that we will need, in a form that is adapted to our purposes, and collect a few less basic or less classical properties that will be useful.

\subsection{Compressions and submatrices}

Let $(E,\langle\cdot,\cdot\rangle)$ be a finite-dimensional inner product space over $\K=\R$ or $\C$. Let $(e_i)_{i\in S}$ be an orthonormal basis of $E$, indexed by some finite set $S$. 

For every subset $I$ of $S$, we set $E_{I}=\Vect(e_{i}:i\in I)$. If $a$ is a linear map with values in $E$, we denote by $a^{I}$ the $E_{I}$-valued linear map obtained by composing~$a$ with the orthogonal projection of $E$ onto $E_{I}$. Similarly, if $J$ is a subset of $S$, and $a$ is a linear map defined on $E$, we denote by~$a_{J}$ the linear map on $E_{J}$ defined by pre-composing $a$ with the injection of~$E_{J}$ into $E$. If $a$ is an endomorphism of $E$, we also define $a^{I}_{J}$ as the linear map from $E_{J}$ to $E_{I}$ obtained by injecting $E_{J}$ into $E$, applying $a$, and projecting orthogonally onto $E_{I}$.

We adopt a similar notation for matrices: if $M$ is a matrix, with rows (resp. columns, resp. both) indexed by $S$, we define $M^{I}$ (resp. $M_{J}$, resp. $M^{I}_{J}$) as the submatrix of $M$ obtained by keeping only the rows with indices in $I$ (resp. the columns with indices in $J$, resp. both).

\subsection{Determinantal measures on finite sets}\label{sec:dpp}
In our finite dimensional inner product space $(E,\langle\cdot,\cdot\rangle)$ endowed with an orthonormal basis $(e_i)_{i\in S}$, let us consider a linear subspace~$H$. Let $\proj{H}$ denote the orthogonal projection on $H$.

A random subset $\X$ of $S$ is said to be \emph{determinantal} associated to $H$ (or determinantal with kernel~$\proj{H}$) in the basis $(e_i)_{i\in S}$ if for all subset $J\subseteq S$, we have
\begin{equation}\label{eq:incidence}
\P(J\subseteq \X)=\det \big[(\proj{H})_{J}^{J}\big].
\end{equation}

It is true, but not obvious, that such a random subset exists. Let us recall one construction of it using Pythagoras' theorem in the exterior algebra of $E$, endowed with the inner product inherited from that of $E$. We also give a translation of this construction in terms of determinants and the Cauchy--Binet formula.  

Set $d=\dim E$ and $n=\dim H$. Let $(h_{1},\ldots,h_{n})$ be a basis of $H$. In $\ext^{n} E$, set $\eta=h_{1}\wedge \ldots \wedge h_{n}$. Also, for every subset $I$ of $S$, let us denote by $e_{I}$ the exterior product of the basis elements $\{e_{i}:i\in I\}$. Since $\eta$ is a non-zero element of $\ext^{n} E$, and $\{e_{I} : I\subseteq S, |I|=n\}$  is an orthonormal basis of this space, we have 
\begin{equation}\label{eq:pythagor}
\|\eta\|^{2}=\sum_{I\subseteq S, |I|=n}  |\langle \eta,e_{I}\rangle|^{2}.
\end{equation}
Thus, there exists a random subset $\X$ of $S$ such that $|\X|=n$ almost surely and
\begin{equation}\label{eq:defloiX}
\P(\X=I)= |\langle \eta,e_{I}\rangle|^{2} \ \big/ \ \|\eta\|^{2} 
\end{equation}
for every subset $I$ of $S$ with $|I|=n$.

Let us rewrite these equalities in matricial terms. Let $A$ be the $d\times n$ matrix whose columns are the vectors $(h_{1},\ldots,h_{n})$ written in the basis $(e_{i})_{i\in S}$ of $E$. For a reason that will appear later in Section~\ref{sec:routine}, we denote by~$A^{\dagger}$ the adjoint matrix of $A$. Then the matrix $A^{\dagger}A$ is invertible, and more precisely, $\det (A^{\dagger}A)=\|\eta\|^{2}$. Moreover, for every subset $I$ of $S$ with $|I|=n$, we have $|\det(A^{I})|^{2}=  |\langle \eta,e_{I}\rangle|^{2}$. Therefore, \eqref{eq:pythagor} and \eqref{eq:defloiX} are equivalent to
\[\det(A^{\dagger}A)=\sum_{I\subseteq S, |I|=n} |\det(A^{I})|^{2} \ \text{ and } \ \P(\X=I)= |\det(A^{I})|^{2}\  /\  \det(A^{\dagger}A),\]
the first equality being an instance of the Cauchy--Binet formula.
 
\begin{proposition}\label{eq:propdetelem}
A random subset $\X$ of $S$ with distribution defined by \eqref{eq:defloiX} is determinantal associated to $H$ in the orthonormal basis $(e_{i})_{i\in S}$ of $(E,\langle\cdot,\cdot\rangle)$. 
\end{proposition}

This is a classical fact, but for the convenience of the reader, and in preparation for what follows, we give a short proof of it. 

\begin{proof} The matrix $A^{\dagger}A$ is invertible and $A(A^{\dagger}A)^{-1}A^{\dagger}$ is the matrix of the orthogonal projection~$\proj{H}$ in the basis $(e_{i})_{i\in S}$. Let us choose a set of indeterminates $\ul\x=(x_{i}:i\in S)$ and introduce the $d\times d$ diagonal matrix $X={\diag(x_{i} : i\in S)}$. Then on the one hand, the Cauchy--Binet formula yields
\begin{align*}
\det(A^{\dagger}(I_{d}+X)A)&=\sum_{I\subseteq S, |I|=n} |\langle \eta,e_{I}\rangle|^{2} \prod_{i\in I}(1+x_{i})\\
&=\sum_{I\subseteq S, |I|=n}\|\eta\|^{2} \,  \P(\X=I) \prod_{i\in I}(1+x_{i})=\det(A^{\dagger}A)\sum_{J\subseteq S} \P(J\subseteq \X) \prod_{j\in J} x_{j}
\end{align*}
and on the other hand, elementary identities for determinants give
\begin{align*}
\det(A^{\dagger}(I_{d}+X)A)&=\det(A^{\dagger}A)\det(I_{n}+(A^{\dagger}A)^{-1}A^{\dagger}XA)\\
&=\det(A^{\dagger}A)\det(I_{d}+{\proj{H}}X)=\det(A^{\dagger}A)\sum_{J\subseteq S} \det \big[(\proj{H})_{J}^{J}\big]\, \prod_{j\in J} x_{j}.
\end{align*}
Comparing the last two expressions, both of which are polynomial in $\ul\x$, we find that $\X$ is indeed determinantal associated to $H$ in the basis $(e_{i})_{i\in S}$.
\end{proof}
We will retain some of the notation introduced in this proof until the end of this section, namely the notation $\dagger$ for the adjoint of a matrix, the indeterminates $\ul{\x}=(x_{i}:i\in S)$, and the diagonal matrix $X={\diag(x_{i} : i\in S)}$.

The formula \eqref{eq:defloiX} is well adapted to situations where it is easy to find a basis of $H$, for example if $H$ is described as the span of a family of vectors of $E$, or if it is given as the range of a linear map into $E$. However, situations occur where $H$ is primarily specified through its orthogonal $H^{\perp}$. The following lemma will allow us to handle this case.

\begin{lemma} \label{lem:orthominors} Let $H$ be a linear subspace of $E$ of dimension $n$. Let $\eta\in \ext^{n}E$ be the exterior product of the elements of a basis of $H$. Let $\zeta\in \ext^{d-n}E$ be the exterior product of the elements of a basis of~$H^{\perp}$. Then for every subset $I$ of $S$ with $|I|=n$, we have
\[|\langle \eta,e_{I}\rangle|^{2}/ \|\eta\|^{2}=|\langle \zeta,e_{S\setminus I}\rangle|^{2}/ \|\zeta\|^{2}.\]
\end{lemma}

\begin{proof} In the language of the exterior algebra, this equality is a consequence of the fact that the Hodge operator on $\ext E$ (see \cite[Section 5.6]{KL3}) is an isometry which sends $\eta/\|\eta\|$ to $\pm\zeta/\|\zeta\|$ and~$e_{I}$ to $\pm e_{S\setminus I}$. 

Let us also give a matricial argument. 
Let $(z_{1},\ldots,z_{d-n})$ be a basis of $H^{\perp}$ of which $\zeta$ is the exterior product. Let $B$ be the $d\times (d-n)$ matrix whose columns are the vectors $(z_{1},\ldots,z_{d-n})$  in the basis $(e_{i})_{i\in S}$. Let $M=(A|B)$ be the $d\times d$ matrix obtained by juxtaposing $A$ and $B$. Let~$D$ be the block diagonal matrix $\diag\big((A^{\dagger}A)^{-\frac{1}{2}},(B^{\dagger}B)^{-\frac{1}{2}})$. Finally, set $R=DM$. Then the matrix $R$ is unitary. 

We now use the fact that complementary minors of a unitary matrix have the same modulus, a consequence of the classical Jacobi formula relating minors of a matrix and its inverse. Therefore, for every subset $I$ of $S$ with $|S|=n$, 
\begin{align*}
|\langle \eta,e_{I}\rangle|^{2}/ \|\eta\|^{2}=|\det(A^{I})|^{2}/\det(A^{\dagger}A)&=|\det(R^{I}_{\{1,\ldots,n\}})|^{2}\\
&\hspace{-2cm}=|\det(R^{S\setminus I}_{[d]\setminus [n]})|^{2}=|\det(B^{I})|^{2}/\det(B^{\dagger}B)=|\langle \zeta,e_{S\setminus I}\rangle|^{2}/ \|\zeta\|^{2}
\end{align*}
and the equality is proved.
\end{proof}

\subsection{A criterion of determinantality}\label{sec:routine}

In order to prove that a random subset $\X$ of $S$ is determinantal, we can check that it satisfies \eqref{eq:incidence}, which describes its incidence measure, or \eqref{eq:defloiX}, which describes its distribution, or the variant of \eqref{eq:defloiX} given by Lemma \ref{lem:orthominors}. Between the lines of the proof of Proposition \ref{eq:propdetelem}, one can read a third characterisation of the distribution of~$\X$, by its generating function. It is a classical property of determinantal processes, that we will not use, and quote here without proof, as an equality of polynomials in the indeterminates $\ul\x$:
\begin{equation}\label{eq:fonctiongeneratriceDPP}
\E\big[\ul\x^{\X}\big] =\det\big(I_{n}+(X-1)\proj{H}\big).
\end{equation}
We will use yet a different characterisation, also in the spirit of generating functions, and that we now describe. 

Let us think of $\ul{\x}$ as a collection of positive weights, and use these weights to twist the inner product on $E$. More precisely, on $(E,\langle\cdot,\cdot\rangle)$, let us define a new inner product $\langle\cdot,\cdot\rangle_{\ul\x}$ by setting
\begin{equation}\label{eq:psx}
\forall u,v\in E, \ \langle u, v\rangle_{\ul\x} =\sum_{i\in S} \x_i \,\overline{u_i} v_i , 
\end{equation}
where $(u_i)_{i\in S}$ and $(v_i)_{i\in S}$ are the coefficients of $u$ and $v$ in the basis $(e_i)_{i\in S}$. Note that $(e_i)_{i\in S}$ is still an orthogonal basis with respect to the twisted inner product.

For endomorphisms defined on $E$ or with values in $E$, we will denote by $\dagger$ the adjunction with respect to the original inner product $\langle\cdot,\cdot\rangle$ and by $*$ the adjunction with respect to the twisted inner product $\langle\cdot,\cdot\rangle_{\ul\x}$. It should therefore be borne in mind that an expression such as $a^{*}$ depends on $\ul\x$.

The criterion of determinantality that we will use is the following. It comes in two versions, that we state and prove together. 

\begin{proposition}\label{prop:routine} 
1. Let $a$ be an injective linear map from some inner product space into~$E$, with rank $n$. As a function of $\ul\x$, the determinant $\det(a^{*}a)$ is a homogeneous polynomial of degree $n$, with non-negative coefficients: 
\begin{equation}\label{eq:assumption-mtt}
\det (a^* a)=\sum_{I\subseteq S, |I|=n}  w(I)\; \ul{\x}^I
\end{equation} 
and for any $\ul\x\in (0,\infty)^{S}$, a random subset $\X$ of $S$ such that
\begin{equation}\label{eq:probaDPP}
\forall I\subseteq S, \ \ \P(\X=I)=w(I)\; \ul{\x}^I \, \big/ \,  \det(a^{*}a)
\end{equation}
is determinantal, associated with the subspace $\im a$ and the orthonormal basis $(e_{i}/\sqrt{x_{i}})_{i\in S}$ of the inner product space $(E,\langle \cdot,\cdot\rangle_{\ul\x})$.

2. Let $b$ be a surjective linear map from $E$ onto some inner product space, with kernel of dimension~$n$.
Then $\ul\x^{S} \det(bb^{*})$ is a homogeneous polynomial of degree $n$ in $\ul\x$ with non-negative coefficients:
\begin{equation}\label{eq:assumption-mtt2}
\ul\x^{S}\det (b b^{*})=\sum_{I\subseteq S, |I|=n}  w(I)\; \ul{\x}^{I}
\end{equation} 
and for any $\ul\x\in (0,\infty)^{S}$, a random subset ${\sf Y}$ of $S$ such that
\begin{equation}\label{eq:probaDPP2}
\forall I\subseteq S, \ \ \P({\sf Y}=I)=w(I)\; \ul{\x}^{I} \, \big/ \, ( \ul\x^{S} \det(bb^{*}))
\end{equation}
is determinantal, associated with the subspace $\ker b$ and the orthonormal basis $(e_{i}/\sqrt{x_{i}})_{i\in S}$ of the inner product space $(E,\langle \cdot,\cdot\rangle_{\ul\x})$.
\end{proposition}

\begin{proof}
1. Let $F$ be an inner product space, and let $a:F\to E$ be an injective linear map. Let us equip $F$ with an orthonormal basis and write matrices with respect to this basis of $F$, and the basis $(e_{i})_{i\in S}$ of $E$. Let $A$ be the matrix of $a$. Then the matrix of $a^{\dagger}$ is $A^{\dagger}$ (hence our unusual notation for the matricial adjunction introduced in Section \ref{sec:dpp}), and the matrix of $a^{*}$ is $A^{\dagger}X$. 

Now, letting $\eta$ be the exterior product of the columns of $A$, we have the equalities
\[\det(a^{*}a)=\det(A^{\dagger} X A)=\sum_{I\subseteq S, |I|=n}|\langle \eta,e_{I}\rangle|^{2} \ \ul\x^{I}=\sum_{I\subseteq S, |I|=n} \big|\big\langle \eta,e_{I}/\sqrt{\ul\x^{I}}\big\rangle_{\ul\x}\big|^{2}=\|\eta\|_{\ul\x}^{2},\]
from which we can read several things. 

The first is that $\det(a^{*}a)$, as a function of $\ul\x$, is indeed a polynomial, and more precisely a homogeneous polynomial of degree $n$.

Secondly, let us choose $\ul\x$ and introduce a random subset $\sf Z$ of $S$, which is determinantal associated to $\im a$ in the basis $(e_{i}/\sqrt{x_{i}})_{i\in S}$ of  $(E,\langle \cdot,\cdot\rangle_{\ul\x})$. Then for every subset $I$ of $S$, the monomial of multidegree $I$ of $\det(a^{*}a)$ is equal to 
\[w(I)\; \ul\x^{I}=\big|\big\langle \eta,e_{I}/\sqrt{\ul\x^{I}}\big\rangle_{\ul\x}\big|^{2}= \|\eta\|_{\ul\x}^{2}\  \P({\sf Z}=I)=\det(a^{*}a)\,\P({\sf Z}=I).\]
Therefore, $\sf Z$ and $\X$ have the same distribution, and the proof of the first statement is finished.
\medskip

2. Let again $F$ be an inner product space, but let now $b:E\to F$ be a surjective linear map. Let us as before equip $F$ with an orthonormal basis and write matrices with respect to this basis of $F$, and the basis $(e_{i})_{i\in S}$ of $E$. Let $B$ be the matrix of $b$. Then the matrix of $b^{\dagger}$ is $B^{\dagger}$, and the matrix of $b^{*}$ is $X^{-1}B^{\dagger}$. We set $d=\dim E$.

Let $\rho\in \ext^{d-n}E$ denote the exterior product of the columns of $B^{\dagger}$. Let $\eta\in \ext^{n}E$ be the exterior product of the elements of an orthonormal basis of the kernel of $B$.  According to Lemma \ref{lem:orthominors}, for every subset $I$ of $S$, we have $|\langle \rho,e_{I}\rangle|^{2}=\|\rho\|^{2} |\langle \eta,e_{S\setminus I}\rangle |^{2}$. 

Therefore, we have the equalities
\begin{align*}
\ul\x^{S} \det(bb^{*})&=\ul\x^{S} \det(BX^{-1}B^{\dagger})=\sum_{|J|=d-n}|\langle \rho,e_{J}\rangle|^{2} \ \ul\x^{S\setminus J}\\
&\hspace{1cm}= \|\rho\|^{2} \sum_{|I|=n}|\langle \eta,e_{I}\rangle|^{2} \ \ul\x^{I}
= \|\rho\|^{2}\sum_{|I|=n} \big|\big\langle \eta,e_{I}/\sqrt{\ul\x^{I}}\big\rangle_{\ul\x}\big|^{2}=\|\rho\|^{2}\|\eta\|_{\ul\x}^{2},
\end{align*}
from which we see that $\ul\x^{S} \det(bb^{*})$ is a homogeneous polynomial of degree $n$ in $\ul\x$.

Let us fix $\ul\x$ and introduce a random subset $\sf Z$ of $S$, which is determinantal associated to $\ker b$ in the basis $(e_{i}/\sqrt{x_{i}})_{i\in S}$ of $(E,\langle \cdot,\cdot\rangle_{\ul\x})$. Then for every subset $I$ of $S$, the monomial of multidegree $I$ of $\ul\x^{S} \det(bb^{*})$ is equal to 
\[w(I)\; \ul\x^{I}=\|\rho\|^{2} \big|\big\langle \eta,e_{I}/\sqrt{\ul\x^{I}}\big\rangle_{\ul\x}\big|^{2}=\|\rho\|^{2} \|\eta\|_{\ul\x}^{2}\  \P({\sf Z}=I)=\det(bb^{*})\,\P({\sf Z}=I).\]
Therefore, $\sf Z$ and $\sf Y$ have the same distribution.
\end{proof}

\subsection{The mean projection theorem}\label{sec:mean-proj}

In this section, we record a variant of a very general property of determinantal processes that we proved in \cite{KL3}, and that seems not to have been previously known. 

As before, we work on a finite-dimensional inner product space $(E,\langle\cdot,\cdot\rangle)$  equipped with an orthonormal basis $(e_{i})_{i\in S}$ and in which a linear subspace $H$ is fixed. We denote by $\X$ the corresponding determinantal random subset of $S$. We will also use its complement $\X^{c}=S\setminus \X$.

Then $E_{\X}=\Vect(e_{i}:i\in \X)$ is a random linear subspace of $E$ and it follows for example from~\eqref{eq:defloiX} that $\P(H^{\perp}\cap E_{\X}\neq\{0\})=0$, so that $E=H^{\perp}\oplus E_{\X}$ almost surely. For the sake of completeness, and in preparation for Section \ref{sec:matroid-setting}, let us state and prove the converse.

\begin{proposition} \label{prop:supportdpp} Let $I$ be a subset of $S$. Then $\P(\X=I)>0$ if and only if $E=H^{\perp}\oplus E_{I}$, if and only if $E=H\oplus E_{I^{c}}$.
\end{proposition}

\begin{proof} The equivalence of the last two conditions follows from the fact that $(E_{I})^{\perp}=E_{I^{c}}$. With the observation made just before the statement of the proposition, there only remains to consider a subset $I$ of $S$ such that $E=H^{\perp} \oplus E_{I}$ and prove that  $\P(\X=I)>0$. For this, we consider two bases of $E$, the first built by joining an orthonormal basis of $H$ and an orthonormal basis of $H^{\perp}$, the second by joining the basis $\{e_{i} : i\in I\}$ of $E_{I}$ and an orthonormal basis of $H^{\perp}$. The change of basis is invertible, blockwise triangular, and $\P(\X=I)$ is the determinant of one of the two diagonal blocks.
\end{proof}

We denote by 
\[\op{P}^{H^{\perp}}_{\X}:E=H^{\perp}\oplus E_{\X} \longrightarrow H^{\perp}\]
the projection on $H^{\perp}$ parallel to $E_{\X}$. Similarly, $E=H\oplus E_{\X^{c}}$ almost surely and we denote by 
\[\op{P}^{H}_{\X^{c}}:E=H\oplus E_{\X^{c}} \longrightarrow H\]
the projection on $H$ parallel to $E_{\X^{c}}$. 
We keep the notation $\proj{H}$ and $\proj{H^{\perp}}$ for the orthogonal projections on $H$ and $H^{\perp}$ respectively.

We use the notation $\ext \op{a}$ for the endomorphism of the exterior algebra $\ext E$ of $E$ induced by an operator $\op{a}\in\End(E)$ (see \cite[Section 5.4]{KL3}). In the basis $\{e_{I} : I\subseteq S\}$ of $\ext E$, the coefficient $(I,J)$ of the matrix of $\ext \op{a}$ is the $(I,J)$-minor  of the matrix of $a$ in the basis $(e_{i})_{i\in S}$.

\begin{theorem}\label{thm:projection} The following equalities hold:
\[\E\big[\ext \op{P}^{H}_{\X^{c}}\big]=\ext \proj{H} \ \text{ and } \ \E\big[\ext \op{P}^{H^\perp}_{\X}\big]=\ext \proj{H^\perp}.\]
\end{theorem}

\begin{proof} The first equality is exactly Theorem 5.9 of \cite{KL3}. To prove the second equality, we apply the same theorem \cite[Theorem 5.9]{KL3} to the determinantal subset of $S$ associated with $\proj{H^{\perp}}$ and use the fact, proved in \cite[Proposition 4.2]{KL3}, and which is also a consequence of Lemma \ref{lem:orthominors}, that~$\X$ is the complement of this determinantal subset.
\end{proof}

\subsection{A conditional probability measure}

The following result will only be needed in Section~\ref{sec:matroid-conditional}, but we record it here for future reference.\footnote{After posting the first version of this article, we noticed the presence of a similar formula in \cite[Eq. (6.5)]{Lyons-DPP}; the statement does not seem to be formulated in the same way, and the proof, by induction in Lyons' paper, is also different, so that we think the following lemma is relevant for our purposes.}

The situation is exactly the same as in the previous section, with the space $E$, the basis $(e_{i})_{i\in S}$, the subspace $H$ and the random subset $\X$ of $S$. We want to describe the distribution of~$\X$ conditional on being included in a subset $K$ of $S$. For this to make sense, we need to have $\P(\X\subseteq K)>0$, which  is equivalent to $E_{S\setminus K}\cap H=\{0\}$, or $E_{K}+H^{\perp}=E$.

\begin{proposition}\label{lem:dpp-conditional-general}
Let $K\subseteq S$ be such that $\P(\X\subseteq K)>0$. Then the random subset $\X$ conditioned on being included in $K$ is a determinantal random subset of $K$ associated to the subspace $\proj{E_K}(H)$ of the inner product space $E_{K}$ with the orthonormal basis $(e_{i})_{i\in K}$.
\end{proposition}

\begin{proof} Set $n=\dim H$. Let $(h_{1},\ldots,h_{n})$ be a basis of $H$. For each $i\in S$, set $f_{i}=\proj{E_{K}}(h_{i})$. Since $\P(\X\subseteq K)>0$, the restriction to~$H$ of the projection $\proj{E_{K}}$ is injective, so that $(f_{1},\ldots,f_{n})$ is a basis of $\proj{E_{K}}(H)$. Let us define 
\[\eta=h_{1}\wedge \ldots \wedge h_{n} \ \text{ and }\ \phi=f_{1}\wedge \ldots \wedge f_{n}.\]
Note that $\phi=(\proj{E_{K}})^{\wedge n}\eta$.

Let us decompose $\eta$ on the orthonormal basis $\{e_{I} : I\subseteq S, |I|=n\}$ of $\ext^{n}E$ and apply $(\proj{E_{K}})^{\wedge n}$ to both sides of the equality. We find
\[\phi=\sum_{I\subseteq K, |I|=n} \langle \eta,e_{I}\rangle e_{I}.\]
Taking squared norms, and recalling \eqref{eq:defloiX}, we find
\[\|\phi\|^{2}=\sum_{I\subseteq K, |I|=n}|\langle \eta,e_{I}\rangle|^{2}=\|\eta\|^{2}\, \P(\X\subseteq K).\]
Therefore, for every subset $I\subseteq K$ with $|I|=n$, we have
\[\P(\X=I | \X\subseteq K)=\frac{|\langle \eta,e_{I}\rangle|^{2}}{\|\eta\|^{2}}\frac{\|\eta\|^{2}}{\|\phi\|^{2}} =\frac{|\langle \eta,e_{I}\rangle|^{2}}{\|\phi\|^{2}}=\frac{|\langle \phi,e_{I}\rangle|^{2}}{\|\phi\|^{2}},\]
and the last quantity is the probability that the determinantal random subset of $K$ associated to $\proj{E_{K}}(H)$ is equal to $I$.
\end{proof}

\subsection{A Schur complement formula}
Let us conclude this section with a statement of the Schur complement formula, under a form that is not quite the standard one, but that will be useful for us. We use the following notation: for every subspace $G$ of an inner product space $E$, we denote by~$1_{G}:G\to E$ the inclusion of $G$ in $E$ and by $1^{G}:E\to G$ its adjoint, that is, the orthogonal projection on $G$.

\begin{lemma} \label{lem:Schurcomp}
Let $a:E_{0}\to E_{1}$ be an injective linear map between finite-dimensional inner product spaces. Let $G$ be a subspace of $E_{0}$. Then 
\[\det(a^{*}a)=\det\big(1^{G^{\perp}}a^{*}a1_{G^{\perp}}\big) \det\big(1^{G}a^{*}\proj{\ker (1^{G^{\perp}}a^{*})}a1_{G}\big)\,.\] 
\end{lemma}

The determinant $\det(a^{*}a)$ is the square of the volume in $E_{1}$ of the image by $a$ of a unit parallelotope in $E_{0}$. This result expresses this volume as a product of two lower-dimensional volumes corresponding to a decomposition of $E_{0}$ into the sum of two orthogonal spaces. 

\begin{proof} Let us start by writing $a=a1_{G^{\perp}}+a1_{G}$ and $a^{*}=1^{G^{\perp}}a^{*}+1^{G}a^{*}$.
Using the Schur complement formula, we find
\begin{align*}
\det(a^{*}a)&=\det\big(1^{G^{\perp}}a^{*}a1_{G^{\perp}}\big) \det\big(1^{G}a^{*}a1_{G}-1^{G}a^{*}a1_{G^{\perp}}  (1^{G^{\perp}}a^{*}a1_{G^{\perp}})^{-1} 1^{G^{\perp}}a^{*}a1_{G}\big)\\
&=\det\big(1^{G^{\perp}}a^{*}a1_{G^{\perp}}\big) \det\big(1^{G}a^{*}\big[ \id_{E_{1}} - a1_{G^{\perp}}  (1^{G^{\perp}}a^{*}a1_{G^{\perp}})^{-1} 1^{G^{\perp}}a^{*}\big]a1_{G}\big).
\end{align*}
Between the square brackets, we have the identity of $E_{1}$ minus the orthogonal projection on the range of $a1_{G^{\perp}}$, that is, the orthogonal projection on $\ker (1^{G^{\perp}}a^{*})$.
\end{proof}


\section{The uniform spanning tree}\label{sec:spanningtrees}
As an important example and in preparation for our next steps, let us apply some of the tools that we collected to the case of random spanning trees. Let us introduce a notation for the Kirchhoff polynomial that we already mentioned in Section \ref{sec:pythagore-kir}: we set
\begin{equation}\label{eq:defTx}
\rT(\ul\x)=\sum_{T\in \cT(\rG)} \ul\x^{T}.
\end{equation}

The following fundamental fact was discovered by Burton and Pemantle, see \cite{Burton-Pemantle}, and we give two proofs of it.

\begin{proposition} \label{prop:lawUST} 
The determinantal random subset $\X$ of $\dE$ associated with the subspace $\im d$ of~$\Omega^{1}(\rG)$  in the orthonormal basis $(e^{\star}/\sqrt{x_{e}})_{e\in \dE}$ is the `uniform' spanning tree of $\rG$. More precisely, the random subgraph $\X$ is almost surely a spanning tree, and for every spanning tree $T\in \cT(\rG)$, we have
\[\P(\X=T)=\ul\x^{T}\big/\, \rT(\ul\x).\]
\end{proposition}

\begin{proof}[First proof] 
Let us choose a vertex $v_{0}$ and denote by $\Omega^{0}(\rG,\{v_{0}\})$ the subspace of $\Omega^{0}(\rG)$ formed by functions that vanish at $v_{0}$. Let $a$ be the restriction to $\Omega^{0}(\rG,\{v_{0}\})$ of the operator $d$. Then the range of $a$ is $\im d$, and \eqref{eq:mtt-classical} gives us an expression of $\det(a^{*}a)$ which, in view of Proposition \ref{prop:routine}, immediately implies the result.
\end{proof}

\begin{proof}[Second proof] According to \eqref{eq:defloiX}, in order to understand the distribution of $\X$, we need to find the decomposition of the exterior product of the elements of a basis of $\im d$ on the orthogonal basis $\{e^{\star}_{S} : S\subseteq \dE, |S|=|\rV|-1\}$ of $\ext^{|\rV|-1}\Omega^{1}(\rG)$. This is exactly what is given to us by \eqref{eq:decomposedetB} and we read on this equality, firstly, that $\X$ is almost surely a spanning tree of $\rG$, and secondly, for all spanning tree $T$, that $\P(\X=T)\propto \|e^{\star}_{T}\|^{2}=\ul\x^{T}$, which is
the expected result.
\end{proof}

Incidentally, we observed in this proof that $\rT(\ul\x)$ is the determinant of the compression on $\Omega^{0}(\rG,\{v_{0}\})$ of $d^{*}d$. Let us record a similar but more intrinsic equality.

\begin{proposition}[Matrix-tree theorem] \label{prop:MTT2}The following equalities hold:
\[\rT(\ul\x)=|\rV|^{-1} \det\Big((d^{*}d)^{(\ker d)^{\perp}}_{(\ker d)^{\perp}}\Big)=|\rV|^{-1} \det\big((d^{*}d)^{\im d^{*}}_{\im d^{*}}\big).\]
\end{proposition}

\begin{proof} The kernel of $d^{*}d$, which is equal to the kernel of $d$, is the line of constant functions on $\rG$. Therefore, $0$ is a simple eigenvalue of $d^{*}d$, and the determinant of the compression of $d^{*}d$ on the orthogonal of its kernel is equal to the coefficient of the monomial of degree $1$ of its characteristic polynomial. This is the sum of all the minors of size $|\rV|-1$ of its matrix in a basis of $\Omega^{0}(\rG)$, and we just observed that in the canonical basis, these minors are all equal, namely to $\rT(\ul\x)$.
\end{proof}

Let us now state the mean projection theorem (Theorem \ref{thm:projection}) in this special case of the uniform spanning tree. Recall the projections $\op{P}^{\ker d^{*}}_{T}$ and $\op{P}^{\im d}_{T^{c}}$ that appear respectively in \eqref{eq:PJZJ} and \eqref{eq:PB}. 

\begin{proposition}\label{prop:projUST} For every integer $k\geq 0$, the following equalities of endomorphisms of $\ext^{k} \Omega^{1}(\rG)$ hold:
\begin{equation}\label{eq:projection-kir}
\sum_{T\in \cT(\rG)} \ul{\x}^T (\rP^{\im d}_{T^{c}})^{\wedge k} = \rT(\ul{\x}) \; (\proj{\im d})^{\wedge k},
\end{equation}
\begin{equation}\label{eq:projection-sym}
\sum_{T\in \cT(\rG)} \ul{\x}^T (\rP^{\ker d^*}_T)^{\wedge k} = \rT(\ul{\x}) \; (\proj{\ker d^*})^{\wedge k} \,.
\end{equation}
\end{proposition}

\begin{proof}
We apply Theorem \ref{thm:projection} and observe that on the event where $\X=T$, the projections  $\op{P}^{H^{\perp}}_{\X}$ and $\op{P}^{H}_{\X^{c}}$
are respectively equal to $\op{P}^{\ker d^{*}}_{T}$ and $\op{P}^{\im d}_{T^{c}}$.
\end{proof}


\section{Geometric multilinear identities}\label{sec:gmi}

In this short section, we combine the results of Sections \ref{sec:metricstructures}, \ref{sec:det-toolbox}, and \ref{sec:spanningtrees} to establish two identities which will play an important role in the proof of our main results in Section \ref{sec:dss}. 

Let us fix an integer $k\geq 0$. For every connected spanning subgraph $K\in \cC_{k}(\rG)$ of $\rG$ with~$k$ linearly independent cycles, the exterior product $\det \mathscr Z^{K}$ of the elements of an integral basis of $Z_{1}(K)$ is defined only up to a sign, but its tensorial square is well defined, and we gave an expression of it in Proposition \ref{prop:sum2-sym}. The $k$-th exterior power of the antilinear isomorphism $\Jx$ defined by \eqref{eq:defJx1} sends $\det \mathscr Z^{K}$ in the inner product space $\ext^{k}\Omega^{1}(\rG)$, of which the family $(e^{\star}_{S}/\sqrt{\ul\x^{S}} : S\subseteq \dE, |S|=k)$ is an orthonormal basis.

The following proposition is a geometric variant of Proposition \ref{prop:sum2-sym}. The identity  is stated in the tensor square {\em over $\R$} of $\ext^{k}\Omega^{1}(\rG)$, even if the base field is $\K=\C$, in order to deal later with the sesquilinearity of the inner product. 

\begin{proposition}\label{prop:geommultconnexe} In $\ext^{k}\Omega^{1}(\rG)\otimes_{\R}\ext^{k}\Omega^{1}(\rG)$, the following equality holds:
\[\sum_{K\in \cC_{k}(\rG)} \ul\x^{K} \big((\Jx)^{\wedge k} \det \mathscr Z^{K}\big)^{\otimes 2}=\rT(\ul\x)\  \sum_{S\subseteq \dE, |S|=k} (\proj{\ker d^{*}})^{\wedge k}\big(e^{\star}_{S}/\sqrt{\ul\x^{S}}\big)\otimes e^{\star}_{S}/\sqrt{\ul\x^{S}}.\]
\end{proposition}

\begin{proof} Taking a weighted sum of the result of Proposition \ref{prop:sum2-sym} over all possible subgraphs $K$, we find
\[\sum_{K\in \cC_{k}(\rG)} \ul\x^{K} \det \mathscr Z^{K} \otimes \det \mathscr Z^{K}=\sum_{T\subseteq K}  \ul\x^{K} ({\rZ}_{T})^{\wedge k} (e_{K\setminus T}) \otimes e_{K\setminus T},\]
where the sum is over all pairs formed by a subgraph $K\in \cC_{k}(\rG)$ and a spanning tree $T$ of $\rG$ contained in $K$. Let us re-index this sum as a sum over pairs formed by a spanning tree $T$ of $\rG$ and a set $S$ of~$k$ edges disjoint from $T$, this set $S$ playing the role of $K\setminus T$. We find
\[\sum_{K\in \cC_{k}(\rG)} \ul\x^{K} \det \mathscr Z^{K} \otimes \det \mathscr Z^{K}=\sum_{T,S}\ul\x^{S} \ul\x^{T} ({\rZ}_{T})^{\wedge k}(e_{S}) \otimes e_{S}\]
and since $({\rZ}_{T})^{\wedge k}(e_{S})=0$ whenever the subset $S$ is not disjoint from the spanning tree $T$, we can read the last sum as a double sum over all spanning trees $T$ of $\rG$ and all $k$-subsets $S$ of $\dE$. Thus, 
\[\sum_{K\in \cC_{k}(\rG)} \ul\x^{K} \det \mathscr Z^{K} \otimes \det \mathscr Z^{K}=\sum_{S}\ul\x^{S} \Big(\sum_{T} \ul\x^{T} ({\rZ}_{T})^{\wedge k}\Big)(e_{S}) \otimes e_{S},\]
an equality originally in $\ext^{k}C_{1}(\rG)\otimes_{\Z}\ext^{k}C_{1}(\rG)$, and that still holds in $\ext^{k}\Omega_{1}(\rG)\otimes_{\R}\ext^{k}\Omega_{1}(\rG)$ after replacing $\rZ_{T}$ by $\rZ_{T}\otimes \id_{\K}$.

Let us now use \eqref{eq:PJZJ} to express $\rZ_{T}\otimes \id_{\K}$ as $\rP^{\ker d^*}_T$ conjugated by the antilinear isomorphism $\Jx$. Let us also use the form of the mean projection theorem given by Proposition \ref{prop:projUST}. We find
\[\sum_{K\in \cC_{k}(\rG)} \ul\x^{K} \det \mathscr Z^{K} \otimes \det \mathscr Z^{K}=\rT(\ul\x)\ \sum_{S}\ul\x^{S} \big((\Jx)^{-1}\proj{\ker d^{*}}\Jx\big)^{\wedge k}(e_{S}) \otimes e_{S}.\]
Applying the $\R$-linear map $(\Jx\big)^{\wedge k}\otimes (\Jx\big)^{\wedge k}$ on both sides of this equality, and using twice the fact that $(\Jx\big)^{\wedge k}(e_{S})=e^{\star}_{S}/\ul\x^{S}$, we find the desired result.
\end{proof}

In order to explain the significance of this result, and to prepare the way in which we will use it, consider an inner product space $(V,\langle\cdot,\cdot\rangle)$ with an orthonormal basis $(v_{1},\ldots,v_{d})$ and a symmetric or Hermitian endomorphism $f$. In the space $V\otimes_{\R}V$, consider the tensor $\sum_{i=1}^{d} f(v_{i})\otimes v_{i}$. Choose now a vector $w\in V$. Then the form $(u,v)\mapsto \langle w,u\rangle \langle v,w\rangle$, which is $\R$-bilinear on $V\times V$ (but not $\C$-bilinear in the Hermitian case), induces a linear form on $V\otimes_{\R}V$, that we denote by $\langle w,\cdot \rangle \otimes \langle \cdot ,w\rangle$. Then we have the equality
\begin{equation}\label{eq:expliqueps}
\big|\langle f(w),w\rangle\big|^{2}=\big(\langle w,\cdot \rangle \otimes \langle \cdot ,w\rangle\big) \Big(\sum_{i=1}^{d} f(v_{i})\otimes v_{i}\Big),
\end{equation}
and the second factor of the right-hand side is precisely the quantity which appears in Proposition~\ref{prop:geommultconnexe}, as well as in Proposition \ref{prop:geommultacyclique} below, to which we now turn.

Let us treat the case of acyclic spanning subgraphs. The integer $k\geq 0$ is still fixed and for every spanning forest $F\in \cF_{k}(\rG)$, the tensor square of the determinant of the elements of an integral basis of $B^{1}(\rG,F)$ is defined without ambiguity.

\begin{proposition}\label{prop:geommultacyclique} In $\ext^{k}\Omega^{1}(\rG)\otimes_{\R}\ext^{k}\Omega^{1}(\rG)$, the following equality holds:
\[\sum_{F\in \cF_{k}(\rG)} \ul\x^{F} (\det \mathscr B^{F})^{\otimes 2}=\rT(\ul\x)\  \sum_{S\subseteq \dE, |S|=k} (\proj{\im d})^{\wedge k}\big(e^{\star}_{S}/\sqrt{\ul\x^{S}}\big)\otimes e^{\star}_{S}/\sqrt{\ul\x^{S}}.\]
\end{proposition}

\begin{proof} The proof is almost the same as that of Proposition \ref{prop:geommultconnexe}, but simpler. We start by applying Proposition \ref{prop:sum2-kir} and summing over all subgraphs $F\in \cF_{k}(\rG)$. We re-index the sum as a sum over pairs formed by a spanning tree $T$ of $\rG$ and a set $S$ of $k$ edges of $T^{+}$. We then observe that the summand makes sense, and is zero, if $S$ is a $k$-susbet of $\dE$ not contained in $T$. Thus, we find
\[\sum_{F\in \cF_{k}(\rG)} \ul\x^{F} \det \mathscr B^{F}\otimes \det \mathscr B^{F}=\sum_{S} (\ul\x^{S})^{-1} \Big(\sum_{T} \ul\x^{T} (\rB_{T})^{\wedge k}\Big)(e^{\star}_{S})\otimes e^{\star}_{S}.\]
We apply \eqref{eq:PB} and then the second statement of Proposition \ref{prop:projUST} to find the result.
\end{proof}


\section{Determinantal spanning subgraphs}\label{sec:dss}

\begin{figure}[!ht]
\centering
\includegraphics[width=6cm]{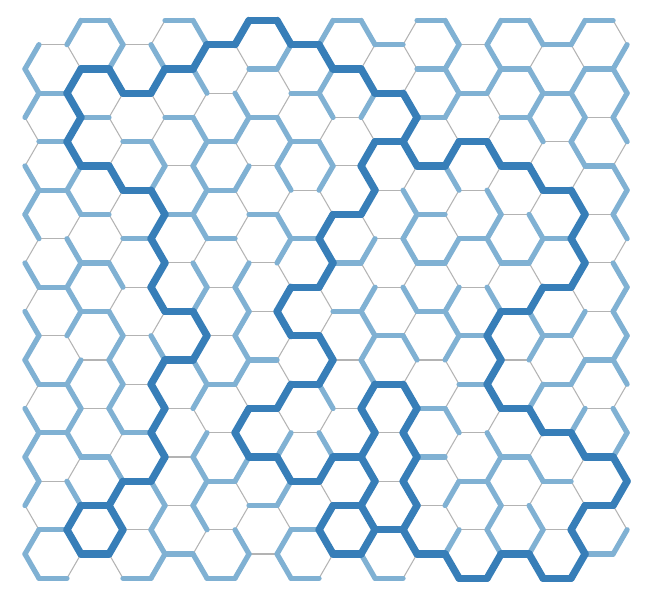} \hspace{1cm} \includegraphics[width=6cm]{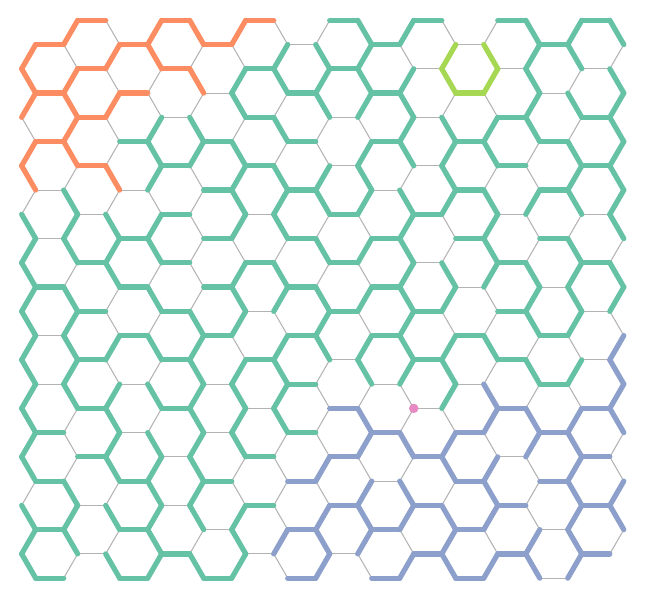}
\caption{\small A random element of $\cC_4(\rG)$ (its $2$-core is represented by thickened edges) and a random element of $\cF_{4}(\rG)$ on a hexagonal grid determined by~$4$ random $1$-forms (note that one of the trees has zero edges). The sampling algorithm we used here, and for the other figures, is the classical one of~\cite[Algorithm 18]{HKPV} for sampling determinantal probability measures.}\label{fig:dpp}
\end{figure}

We now introduce two families of determinantal random subgraphs, which are supported respectively by $\cC_{k}(\rG)$ and $\cF_{k}(\rG)$ for some integer $k\geq 0$, see Figure~\ref{fig:dpp} below for an illustration. Later, in Section~\ref{sec:betti-fixed}, we will extend these definitions to a family of determinantal measures on more general families of subgraphs with constrained Betti numbers. 

In general, and in contrast with the case of spanning trees, the probability measures that we will define are \emph{not} uniform. 
In fact, in view of the complexity results stated below (see Section \ref{sec:complexity}), the uniform measure on $\cC_k(\rG)$ and~$\cF_k(\rG)$ cannot be determinantal in general (otherwise, summing over $k$, there would be a determinantal formula for enumerating $\cC(\rG)$ and $\cF(\rG)$, contradicting the $\# P$-hardness of these counts). Studying the uniform measure would be a more difficult task; on that matter, see~\cite{Grimmett-Winkler} for conjectures about the uniform measure on connected subgraphs, or spanning forests, without constraint on the Betti numbers.

\subsection{Determinantal connected spanning subgraphs}\label{sec:connected}

Our first result describes the determinantal random subset of $\rE$ associated, in the orthonormal basis $(e^{\star}/\sqrt{x_{e}})_{e\in \dE}$ of $\Omega^{1}(\rG)$, to a linear subspace $H$ that {\em contains} $\im d$, and with a dimension exceeding that of $\im d$ by some integer $k\geq 0$.

Let us therefore choose an integer $k\geq 0$ and a linear subspace $H$ of $\Omega^{1}(\rG)$, such that 
\[\im d \subseteq H \ \text{ and } \ \dim H = \rank(d)+k=|\rV|-1+k.\]
Let us choose $k$ vectors $\theta_{1},\ldots,\theta_{k}$ of $\Omega^{1}(\rG)$ such that
\[H=\im d \oplus \Vect(\theta_{1},\ldots,\theta_{k})\]
and set
\[\vartheta=\theta_{1}\wedge \ldots \wedge \theta_{k} \in \ext^{k}\Omega^{1}(\rG).\]

For every subgraph $K\in \cC_{k}(\rG)$ of $\rG$, the $\Z$-module $Z_{1}(K)$ has rank $k$. The exterior product of the elements of an integral basis $\mathscr Z^{K}$ of $Z_{1}(K)$ depends on the choice of this basis only up to a sign: we denote it by $\det \mathscr Z^{K}$ and see it as an element of $\ext^{k}\Omega_{1}(\rG)$.

Using this element, we can associate to $K$ the non-negative real weight $|(\vartheta,\det \mathscr Z^{K})|^{2} \ul\x^{K}$ which, once an integral basis $(z_{1}, \ldots, z_{k})$ of $Z_1(K)$ is chosen, can be written more concretely using
\begin{equation}\label{eq:weightSym}
|(\vartheta,\det \mathscr Z^{K})|^{2} =\big\vert\det\big(\theta_i(z_{j})\big)_{1\le i,j\le k}\big\vert^2.
\end{equation}

Let us define the generating polynomial 
\begin{equation}\label{eq:defCpol}
{\rC}_{k}(\vartheta,\ul\x)=\sum_{K\in \cC_{k}(\rG)} |(\vartheta,\det \mathscr Z^{K})|^{2} \, \ul\x^{K}
\end{equation}
of weighted connected spanning subgraphs with $k$ independent cycles.  
Let us emphasize that the definition \eqref{eq:defCpol} makes sense for every $\vartheta\in \ext^{k}\Omega^{1}(\rG)$. 

\begin{theorem}\label{thm:DPP1}
Let $\X$ be a determinantal random subset of $\rE$ associated to $H$ in the orthonormal basis $(e^{\star}/\sqrt{x_{e}})_{e\in \dE}$ of  $\Omega^{1}(\rG)$. Then $\X$ belongs to $\cC_{k}(\rG)$ almost surely and for every $K\in \cC_{k}(\rG)$, 
\begin{equation}\label{eq:distribK}
\P(\X=K)= |(\vartheta,\det\mathscr Z^{K})|^{2}\  \ul\x^{K} \ \big/ \ \rC_{k}(\vartheta,\ul\x).
\end{equation}
\end{theorem}

Just as we gave two proofs of the fundamental Proposition \ref{prop:lawUST}, we will give two proofs of this theorem. The first relies on Proposition \ref{prop:routine}, for the application of which we need to introduce a linear map with range equal to $H$. To this end, let us consider the map 
\[\omega_{\vartheta} : \C^k \to \Omega^1(\rG), \ \ (t_{1},\ldots, t_{k})\mapsto t_{1} \theta_1+\ldots+ t_{k} \theta_k.\]
This map does not only depend on the tensor $\vartheta$, but on the whole family $(\theta_{1},\ldots,\theta_{k})$, so that our notation is slightly abusive. 
The space $\C^{k}$ being endowed with the usual Hermitian inner product, let us consider the orthogonal direct sum $(\ker d)^{\perp} \oplus \C^{k}$ and the linear maps
\[d\oplus \omega_{\vartheta}:(\ker d)^{\perp}\oplus \C^{k} \to \Omega^{1}(\rG) \ \text{ and } \ \Delta_{\vartheta}=(d\oplus \omega_{\vartheta})^*(d\oplus \omega_{\vartheta}).\]

The following proposition is the key to our first proof of Theorem \ref{thm:DPP1}.

\begin{proposition} \label{prop:egalitesconnexes} One has the equalities
\begin{equation}\label{eq:connexetriangle}
\rC_{k}(\vartheta,\ul\x)=\rT(\ul\x) \ \big\| (\proj{\ker d^{*}})^{\wedge k}(\vartheta)\big\|^{2}=|\rV|^{-1} \det \Delta_{\vartheta}.
\end{equation}
\end{proposition}

\begin{proof}[Proof of Proposition \ref{prop:egalitesconnexes}]
Let us start by proving the second equality of \eqref{eq:connexetriangle}. For this, let us apply the Schur complement formula under the form given by Lemma \ref{lem:Schurcomp}, with $E_{0}=(\ker d)^{\perp}\oplus \C^{k}$, $E_{1}=\Omega^{1}(\rG)$, $G=\C^{k}$ and $a=d\oplus \omega_{\vartheta}$. We find
\[\det \Delta_{\vartheta} =\det(d^*d)^{\im d^{*}}_{\im d^{*}}\  \det\big(\omega_{\vartheta}^{*}\proj{\ker d^{*}} \omega_{\vartheta}\big).
\]
The first factor is equal to $|\rV| \rT(\ul\x)$, by Proposition \ref{prop:MTT2}. Let us compute the second. For this, let us observe that the adjoint of $\omega_{\vartheta}$ is given, for all $\alpha\in \Omega^{1}(\rG)$, by $\omega_{\vartheta}^{*}(\alpha)= (\lanx\theta_{1},\alpha\ranx,\ldots,\lanx\theta_{k},\alpha\ranx)$, so that the matrix in the canonical basis of $\C^{k}$ of $\omega_{\vartheta}^{*}\proj{\ker d^{*}} \omega_{\vartheta}$ is
$ \big(\lanx \theta_{i},\proj{\ker d^{*}}\theta_{j} \ranx\big)_{1\leq i,j \leq k}$
and 
\[\det\big(\omega_{\vartheta}^{*}\proj{\ker d^{*}} \omega_{\vartheta}\big)= \big\| (\proj{\ker d^{*}})^{\wedge k} \vartheta \big\|^{2}.\]

Let us now prove the first equality of \eqref{eq:connexetriangle}. Starting from the definition of $\rC_{k}(\vartheta,\ul\x)$ and using Proposition \ref{prop:geommultconnexe}, we find that $\rC_{k}(\vartheta,\ul\x)$ is equal to 
\[\sum_{K\in \cC_{k}(\rG)}\ul\x^{K} \lanx \vartheta,(\Jx)^{\wedge k} \det \mathscr Z^{K} \ranx \lanx (\Jx)^{\wedge k} \det \mathscr Z^{K},\vartheta\ranx =\big(\lanx \vartheta,\cdot\ranx \otimes \lanx \cdot ,\vartheta\ranx\big) \sum_{K\in \cC_{k}(\rG)}  \ul\x^{K}\big((\Jx)^{\wedge k} \det \mathscr Z^{K}\big)^{\otimes 2}.\]
As explained in the paragraph preceding \eqref{eq:expliqueps}, this is equal to
\[\rT(\ul\x) \ \blanx (\proj{\ker d^{*}})^{\wedge k} \vartheta,\vartheta \branx=\rT(\ul\x) \ \big\|(\proj{\ker d^{*}})^{\wedge k} \vartheta\big\|^{2}\]
and the proof is complete.
\end{proof}

Let us now see how this proposition implies the theorem.

\begin{proof}[First proof of Theorem \ref{thm:DPP1}] In view of the equality of the first and third terms of \eqref{eq:connexetriangle}, an application of the first part of Proposition \ref{prop:routine} to the operator 
 $d\oplus \omega_{\vartheta}$, whose range is the subspace $H$ of $\Omega^{1}(\rG)$, implies that the random subset $\X$ is determinantal with distribution given by \eqref{eq:distribK}.
\end{proof}

Let us give a second proof of Theorem \ref{thm:DPP1}. It seems longer than the first one, but it is more elementary, in the sense that it relies less heavily on the content of the previous sections. More specifically, it depends neither on Proposition \ref{prop:egalitesconnexes} nor on the results of Section \ref{sec:gmi}.

\begin{proof}[Second proof of Theorem \ref{thm:DPP1}] This proof relies directly on the description given by \eqref{eq:defloiX} of the distribution of the determinantal random subgraph~$\X$.

The first thing to check is that $\X$ almost surely belongs to $\cC_{k}(\rG)$. The number of edges of $\X$ is almost surely $\dim H=|\rV|+k-1$, so that its Betti numbers satisfy $b_{0}-b_{1}=1-k$ (see \eqref{eq:Euler}). A subgraph $S$ satisfying $b_{0}-b_{1}=1-k$ and that does not belong to $\cC_{k}(\rG)$ must satisfy $b_{1}>k$, which means that it must have at least $k+1$ linearly independent cycles. By an argument of dimension, there exists a non-trivial $\K$-linear combination of these cycles  that annihilates each of the $1$-forms $\theta_{1},\ldots,\theta_{k}$ and also, because this linear combination is still a cycle, every $1$-form of $\im d$. Therefore, there exists a non-zero element of $\Omega^{1}(S)$ that is orthogonal to $H$, and $\P(\X=S)=0$.

Let us now fix a subgraph $K\in \cC_{k}(\rG)$. We will compute $\P(\X=K)$ using \eqref{eq:defloiX}, after observing that the right-hand side of this equality is unaffected by the multiplication of the tensor $\eta$ by a scalar, and therefore does not depend on the choice of the basis of $H$ used to produce $\eta$. We have therefore the liberty, for every subgraph $K$, to apply this formula with any particular basis of $H$ that we find suitable.
 
Let us, then, construct a basis of $H$ that is well-suited to the computation of $\P(\X=K)$. Let us start by picking a vertex $v_{0}\in \rV$ and considering the basis $\mathscr B_{v_{0}}=\{d\1_{v} : v\neq v_{0}\}$ of $\im d$. In order to complete $\mathscr B_{v_{0}}$ into a basis of $H$, let us choose a spanning tree $T$ of $\rG$ contained in $K$. According to \eqref{eq:splitOmega1dTc}, any $1$-form on $\rG$ can be made to vanish on $T$ by adding to it an appropriate coboundary. To be more specific, if $\theta$ is a $1$-form, then $(\rB_{T}\otimes \id_{\K})\theta$ is a coboundary and $\theta'=\theta-(\rB_{T}\otimes \id_{\K})\theta$ vanishes on $T$. Moreover, for every $1$-chain $c\in C_{1}(\rG)$, Proposition \ref{prop:decBZ} implies 
\[(\theta',c)=(\theta,\rZ_{T}(c)).\]
The $1$-forms $\theta_{1},\ldots,\theta_{k}$ still being those used to define $H$, and with the obvious notation, let us consider the basis $\mathscr H_{T}=\{d\1_{v} : v\neq v_{0}\}\cup \{\theta'_{1},\ldots,\theta'_{k}\}$ of $H$. Let us observe that since $\theta'_{i}$ differs from $\theta_{i}$ by a coboundary, the exterior product of the elements of $\mathscr H_{T}$ is, up to a sign, equal to that of the elements of $\mathscr H=\{d\1_{v} : v\neq v_{0}\}\cup \{\theta_{1},\ldots,\theta_{k}\}$. In particular, $\|\det \mathscr H_{T}\|^{2}=\|\det \mathscr H\|^{2}$ does not depend on $T$.

The last ingredient needed to apply \eqref{eq:defloiX} is to recognize the tensor denoted there by $e_{I}$ as being, in the present setting, the exterior product of the elements corresponding to $K$ in the fixed orthonormal basis of $\Omega^{1}(\rG)$. Once this is done, we find
\[\P(\X=K)=\big|\lanx \det \mathscr H_{T},e^{\star}_{K}/\sqrt{\ul\x^{K}}\ranx\big|^{2} \ / \ \|\det \mathscr H_{T}\|^{2}.\]
Thanks to our choice of $\mathscr H_{T}$, the scalar product in the numerator is the determinant of a block triangular matrix. Labelling the edges of $K\setminus T$ as $\{e_{1},\ldots,e_{k}\}$, we have 
\[\lanx \det \mathscr H_{T},e^{\star}_{K}/\sqrt{\ul\x^{K}}\ranx=\pm(\ul\x^{K})^{-\frac{1}{2}}\lanx \det \mathscr B_{v_{0}},e^{\star}_{T}\ranx \det(\lanx \theta'_{i},e^{\star}_{j}\ranx)_{1\leq i,j \leq k}.\]
According to \eqref{eq:decomposedetB}, the second factor on the right-hand side is equal to $\|e^{\star}_{T}\|^{2}=\ul\x^{T}$. To compute the third factor, let us observe that for each $i,j\in \{1,\ldots,k\}$, 
\[\lanx \theta'_{i},e^{\star}_{j}\ranx=x_{j} (\theta'_{i},e_{j})=x_{j}\big(\theta_{i},\rZ_{T}(e_{j})\big),\]
so that the determinant is equal to $\pm \ul\x^{K\setminus T} (\vartheta,z_{K})$.

Therefore,
\[\P(\X=K)= |(\vartheta,z_{K})|^{2}\  \ul\x^{K} \ \big/ \ \|\det \mathscr H_{T}\|^{2}\] 
and since $\|\det \mathscr H_{T}\|^{2}$ does not depend on $T$, the result is proved.
\end{proof}

\subsection{Determinantal acyclic spanning subgraphs}\label{sec:forests}

We now turn our attention to the study of the determinantal random subset of $\rE$ associated, in the orthonormal basis $(e^{\star}/\sqrt{x_{e}})_{e\in \dE}$ of $\Omega^{1}(\rG)$, to a linear subspace $H$ that is {\em contained} in $\im d$, and with a dimension equal to that of $\im d$ minus some integer $k\geq 0$. 

Let us choose an integer $k\geq 0$ and a linear subspace $H$ of $\Omega^{1}(\rG)$, such that 
\[H \subseteq \im d \ \text{ and } \ \dim H = \rank(d)-k=|\rV|-1-k.\]
An appropriate way to describe $H$ is as the intersection of the kernels of $k$ linear forms on $\im d$. Let us therefore choose $k$ chains $c_{1},\ldots,c_{k}$ of $\Omega_{1}(\rG)$ such that
\[H=\im d \cap \Vect(c_{1},\ldots,c_{k})^{\circ}=\{\alpha\in \im d : (\alpha,c_{1})=\ldots=(\alpha,c_{k})=0\}.\]
We set
\[\calc=c_{1}\wedge \ldots \wedge c_{k} \in \ext^{k}\Omega_{1}(\rG).\]

For every subgraph $F\in \cF_{k}(\rG)$ of $\rG$, the $\Z$-module $B^{1}(\rG,F)$ has rank $k$ and the exterior product of the elements of an integral basis $\mathscr B^{F}$ of $B^{1}(\rG,F)$ 
is an element of the exterior algebra of $\Omega^{1}(\rG)$, defined up to a sign, that we denote by $\det \mathscr B^{F}$.

We associate to $F$ the non-negative real weight $|(\det \mathscr B^{F},\calc)|^{2} \ul\x^{F}$ which, once chosen an integral basis $(\beta_{1},\ldots,\beta_{k})$ of $B^{1}(\rG,F)$, can be written as
\begin{equation}\label{eq:weight-kir}
|(\det \mathscr B^{F},\calc)|^{2}=\big\vert\det\big(\beta_i(c_j)\big)_{1\le i,j\le k}\big\vert^2.
\end{equation}

Let us define the generating polynomial 
\begin{equation}\label{eq:defApol}
{\rA}_{k}(\calc,\ul\x)=\sum_{F\in \cF_{k}(\rG)}|(\det \mathscr B^{F},\calc)|^{2} \, \ul\x^{F}
\end{equation}
of weighted acyclic spanning subgraphs with $k+1$ connected components.  
Let us emphasize that the definition \eqref{eq:defCpol} makes sense for every $\vartheta\in \ext^{k}\Omega^{1}(\rG)$. 

\begin{theorem}\label{thm:DPP2}
Let $\X$ be a determinantal random subset of $\rE$ associated to $H$ in the orthonormal basis $(e^{\star}/\sqrt{x_{e}})_{e\in \dE}$ of  $\Omega^{1}(\rG)$. Then $\X$ belongs to $\cF_{k}(\rG)$ almost surely and for every $F\in \cF_{k}(\rG)$, 
\begin{equation}\label{eq:distribF}
\P(\X=F)= |(\det \mathscr B^{F},\calc)|^{2}\  \ul\x^{F} \ \big/ \ \rA_{k}(\calc,\ul\x).
\end{equation}
\end{theorem}

See the right of Figure~\ref{fig:dpp} for a sample of this measure in an example.

As in the connected case, we will give two proofs of this theorem. The first relies on Proposition~\ref{prop:routine} and begins with the construction of a linear map of which $H$ is the kernel. In this proof, we will make use of the inner product on $\Omega^{1}(\rG)$ associated with conductances identically equal to $1$. We call this the untwisted inner product, and use the adjective untwisted to qualify the related notions and objects, such as orthogonality and adjunction.  

Let us define the linear subspace 
\[Q={\sf J}^{1}_{\ul 1}(Z_{1}(\rG)\otimes \K)=\ker d^{\dagger}\]
of $\Omega^{1}(\rG)$. Let us endow $Q$ with the untwisted inner product and consider the untwisted orthogonal projection
\[p:\Omega^{1}(\rG)\to Q=\ker d^{\dagger}.\]
Let us also endow $\C^{k}$ with the usual inner product and define 
\[\pi_{\calc}:\Omega^{1}(\rG) \to \C^{k}, \ \alpha \mapsto \big((\alpha,c_{1}),\ldots,(\alpha,c_{k})\big).\]
Let us finally form the orthogonal direct sum $Q\oplus \C^{k}$ and consider the operator
\[b=p+\pi_{\calc}:\Omega^{1}(\rG) \to Q\oplus \C^{k}.\]
Then the kernel of $b$ is the subspace $H$.

We will use the equality $p^{*}\circ {\sf J}^{1}_{\ul 1}=\Jx$, which follows from the fact that for every $1$-form $\alpha\in Q$ and every $1$-chain $c$,
\[\lanx p^{*}{\sf J}^{1}_{\ul 1}(c),\alpha\ranx_{\ul\x}=\lanx {\sf J}^{1}_{\ul 1}(c),\alpha\ranx_{\ul 1}=(\alpha,c)=\lanx \Jx \alpha,c\ranx_{\ul\x},\]
where, for the sake of clarity, we kept explicit track of the edge weights $\ul\x$.

\begin{proposition} \label{prop:egalitesacycliques} One has the equalities
\begin{equation}\label{eq:acycliquetriangle}
\rA_{k}(\calc,\ul\x)=\rT(\ul\x) \, \big\| (\proj{\im d}\Jx)^{\wedge k}(\calc)\big\|^{2}= \rT(\ul{1})\, \ul\x^{\rE} \det (bb^{*}).
\end{equation}
\end{proposition}

\begin{proof}[Proof of Proposition \ref{prop:egalitesacycliques}]
Let us start by proving the second equality of \eqref{eq:acycliquetriangle}. For this, let us apply the Schur complement formula under the form given by Lemma \ref{lem:Schurcomp}, with $E_{0}=Q\oplus \C^{k}$, $E_{1}=(\Omega^{1}(\rG),\lanx \cdot , \cdot \ranx)$, $G=\C^{k}$ and $a=b^{*}$. We have $\ker (1^{Q}b)=\im d$, and we find
\[\det (bb^{*}) =\det(1^{Q}bb^{*}1_{Q})\det(\pi_{\calc} \proj{\im d} \pi_{\calc}^{*}).\]
The first determinant on the right-hand side is the square of the volume of the image by $b^{*}$ of a basis of volume $1$ of $Q$. We do not have a basis of volume $1$ of $Q$, but we can consider an integral basis $(z_{1},\ldots,z_{b_{1}})$ of $Z_{1}(\rG)$ and take its image by ${\sf J}^{1}_{\ul 1}$. This produces a basis of $Q$, of which the square of the volume is computed by \eqref{eq:sym1-omid} with all edge weights equal to $1$, and is equal to $\rT(\ul 1)$. Let us take the image by $b^{*}$ of this basis of $Q$. Since $p^{*}\circ {\sf J}^{1}_{\ul 1}=\Jx$, it is the image by $\Jx$ of $(z_{1},\ldots,z_{b_{1}})$,  the volume of which is also computed by \eqref{eq:sym1-omid}, now with the edge weights $\ul\x$, and is equal to $\rT(\ul\x)/\ul\x^{\rE}$.

Therefore, we have
\[\det(1^{Q}bb^{*}1_{Q})= \frac{\rT(\ul x)}{\ul x^{\rE} \rT(\ul 1)}.\]
To compute the second determinant, let us first observe that the adjoint of $\pi_{\calc}$ is given by 
\[\pi_{\calc}^{*}(t_{1},\ldots,t_{k})=t_{1}\Jx c_{1} + \ldots + t_{k}\Jx c_{k}.\]
Therefore, 
\[\det(\pi_{\calc} \proj{\im d} \pi_{\calc}^{*})=\det\big(\blanx\Jx c_{i},\proj{\im d} \Jx c_{j}\branx_{1\leq i,j\leq k}\big)=\big\| (\proj{\im d}\Jx)^{\wedge k}(\calc)\big\|^{2}\]
and the second equality of \eqref{eq:acycliquetriangle} is proved.

Let us now prove the first equality of \eqref{eq:acycliquetriangle}. Let us introduce the notation $\varsigma=(\Jx)^{\wedge k} \calc$. Starting from the definition of $\rA_{k}(\calc,\ul\x)$ and using Proposition \ref{prop:geommultacyclique}, we find
\[\rA_{k}(\calc,\ul\x)=\sum_{F\in \cF_{k}(\rG)}\ul\x^{F} \lanx\varsigma,\det \mathscr B^{F}\ranx \lanx \det \mathscr B^{F},\varsigma\ranx =\big(\lanx \varsigma,\cdot\ranx \otimes \lanx \cdot ,\varsigma\ranx\big) \Big(\sum_{F\in \cF_{k}(\rG)}  \ul\x^{F}\big(\det \mathscr B^{F}\big)^{\otimes 2}\Big).\]
As explained in the paragraph preceding \eqref{eq:expliqueps}, this is equal to
\[\rT(\ul\x) \ \blanx (\proj{\im d})^{\wedge k}\varsigma ,\varsigma \branx=\rT(\ul\x) \ \big\|(\proj{\im d}\Jx)^{\wedge k} \calc\big\|^{2}\]
and the proof is complete.
\end{proof}

As in the connected case, this proposition implies immediately our theorem.

\begin{proof}[First proof of Theorem \ref{thm:DPP2}] In view of the equality of the first and third terms of \eqref{eq:acycliquetriangle}, an application of Proposition \ref{prop:routine} to the operator $b$, whose kernel is the subspace $H$ of $\Omega^{1}(\rG)$, implies that the random subset $\X$ is determinantal with distribution given by \eqref{eq:distribF}.
\end{proof}

Let us give a second proof of Theorem \ref{thm:DPP2}, relying directly on \eqref{eq:defloiX}, and on Lemma \ref{lem:orthominors}.

\begin{proof}[Second proof of Theorem \ref{thm:DPP2}] Let us start by checking that $\X$ almost surely belongs to $\cF_{k}(\rG)$. It has almost surely $\dim H=|\rV|-1-k$ edges, so that its Betti numbers (see \eqref{eq:Euler}) satisfy $b_{0}-b_{1}=k+1$. A subgraph $S$ satisfying $b_{0}-b_{1}=k+1$ and that does not belong to $\cF_{k}(\rG)$ must have $b_{1}\geq 1$, which means that it must have a non-trivial cycle. This cycle annihilates every $1$-form of $\im d$, and in particular every element of $H$. Therefore, there exists a non-zero element of $\Omega^{1}(S)$ that is orthogonal to $H$, and $\P(\X=S)=0$.

Let us now fix a subgraph $F\in \cF_{k}(\rG)$ and compute $\P(\X=F)$ using Lemma \ref{lem:orthominors}. For this, we need to construct a basis of the orthogonal of $H$ in $\Omega^{1}(\rG)$. Let us emphasize that in this proof,~$\Omega^{1}(\rG)$ is always endowed with the (twisted) inner product defined by \eqref{eq:inner-product1}.

Let us start by taking an integral basis~$\mathscr Z$ of $Z_{1}(\rG)$, and considering its image by $\Jx$ in $\Omega^{1}(\rG)$, where it becomes a basis of $\ker d^{*}$. To complete this basis into a basis of the orthogonal of $H$, let us choose a spanning tree $T$ containing $F$. According to \eqref{eq:splitOmega1d*T}, any $1$-form on $\rG$ can be made to vanish on $T^{c}$ by the addition of a suitable element of $\ker d^{*}$. More precisely, if $\alpha$ is a $1$-form, then ${\Jx} (\rZ_{T}\otimes {\rm id}_{\K})({\Jx})^{-1} \alpha$ belongs to $\ker d^{*}$ and $\alpha'=\alpha - {\Jx} (\rZ_{T}\otimes {\rm id}_{\K})({\Jx})^{-1} \alpha$ vanishes on $T^{c}$. Moreover, for every $1$-chain $c\in \Omega_{1}(\rG)$, Proposition \ref{prop:decBZ} implies
\[(\alpha',c)=(\alpha,(\Jx)^{-1}(\rB_{T}\otimes \id_{\K})\Jx(c))=\lanx \alpha,(\rB_{T}\otimes \id_{\K})\Jx(c)\ranx.\]

For each $i\in \{1,\ldots,k\}$, let us define $\sigma_{i}=\proj{\im d}\Jx (c_{i})$. Then 
\[\mathscr H^{\perp}=\Jx (\mathscr Z )\cup \{\sigma_{1},\ldots,\sigma_{k}\} \ \text{ and } \ \mathscr H^{\perp}_{T}=\Jx (\mathscr Z )\cup \{\sigma'_{1},\ldots,\sigma'_{k}\}\]
are two bases of $H^{\perp}$ which differ by a triangular change of basis with unit diagonal, and therefore satisfy $\|\det \mathscr H_{T}^{\perp}\|^{2}=\|\det \mathscr H^{\perp}\|^{2}$.

According to \eqref{eq:defloiX} and Lemma \ref{lem:orthominors}, we have
\[\P(\X=F)=\big|\lanx \det \mathscr H^{\perp}_{T},e^{\star}_{F^{c}}/\sqrt{\ul\x^{F^{c}}}\ranx\big|^{2} \ / \ \|\det \mathscr H^{\perp}_{T}\|^{2}.\]
Labelling the edges of $T\setminus F$ as $\{e_{1},\ldots,e_{k}\}$, we have
\[\lanx \det \mathscr H^{\perp}_{T},e^{\star}_{F^{c}}/\sqrt{\ul\x^{F^{c}}}\ranx=\pm (\ul\x^{F^{c}})^{-\frac{1}{2}} \lanx \det \Jx (\mathscr Z),e^{\star}_{T^{c}}\ranx \det\big(\lanx \sigma_{i}',e^{\star}_{j}\ranx \big)_{1\leq i,j\leq k}.
\]
By \eqref{eq:decomposedetZ}, and recalling the definition \eqref{eq:defJx1} of $\Jx$, the first scalar product on the right-hand side is equal to $\|e^{\star}_{T^{c}}\|^{2}/\ul\x^{T^{c}}=1$. To compute the last factor of the right-hand side, let us compute
\begin{align*}
\lanx \sigma_{i}',e^{\star}_{j}\ranx=x_{j} (\sigma'_{i},e_{j}) =x_{j}\lanx \proj{\im d} \Jx c_{i}, (\rB_{T}\otimes \id_{\K})\Jx(e_{j})\ranx 
&=x_{j}\lanx  \Jx c_{i}, (\rB_{T}\otimes \id_{\K})\Jx(e_{j})\ranx \\
&=\lanx  \Jx c_{i}, (\rB_{T}\otimes \id_{\K}) e^{\star}_{j}\ranx 
=(\rB_{T} e^{\star}_{j},c_{i})
\end{align*}
so that the last factor is equal to $\pm \ul\x^{T\setminus F} (\det \mathscr B^{F},\calc)$. Finally, we find
\[\P(\X=F)=\ul\x^{F} |(\det \mathscr B^{F},\calc)|^{2} \ \big/ \ \ul x^{\rE} \|\det \mathscr H^{\perp}_{T}\|^{2}\]
and since $\|\det \mathscr H_{T}^{\perp}\|^{2}$ does not depend on $T$, the result is proved.
\end{proof}

\subsection{Duality and the point of view of matroids}\label{sec:duality}

The two situations we described above (Sections \ref{sec:connected} and \ref{sec:forests}) are in fact dual from one another: in the first case, we consider subspaces containing a reference subspace~$H=\im d$, in the second case, we consider subspaces contained in this reference subspace. In terms of configuration of edges, it corresponds to adding edges to the random tree or removing edges. This duality (direct sum versus intersection, set union versus set intersection) is more apparent within the language of matroids.

In Section \ref{sec:matroid-setting}, we will revisit these situations from the point of view of the theory of matroids. Precise definitions will be given there, and can also be found in \cite{Oxley}, but we conclude this section with a brief preview of this approach. 

It is a general fact that the support of a determinantal point process on a finite set is the set of bases of a matroid on this set, see \cite{Lyons-DPP}. For example, the set of spanning trees of $\rG$ is the set of bases of the {\em circular matroid} on the set of edges of $\rG$.

The support of the distribution of the determinantal subgraph of $\rG$ described by Theorem~\ref{thm:DPP1} is the set of bases of a matroid contained in 
$\cC_k(\rG)$, and which is equal to $\cC_k(\rG)$ for a generic choice of~$\vartheta$. In particular, $\cC_k(\rG)$ is the set of bases of a matroid on the set of edges of $\rG$. This matroid is obtained by forming the {\em matroid union} of the circular matroid and the uniform matroid of rank $k$. 

Similarly, the support of the distribution of the determinantal subgraph of $\rG$ described by Theorem~\ref{thm:DPP2} is the set of bases of a matroid contained in $\cF_k(\rG)$, which is equal to $\cF_k(\rG)$ for a generic choice of $\calc$. In this case, $\cF_k(\rG)$ is the {\em dual matroid} of the union of the dual of the circular matroid and the uniform matroid of rank $k$. It is also called the \emph{$r$-truncation} of the circular matroid, with $r=\vert\rV\vert-1-k$. 


\section{Duality in the two-dimensional case}\label{sec:two-dimensional-case}

\begin{figure}[!ht]
\centering
\includegraphics[width=6cm]{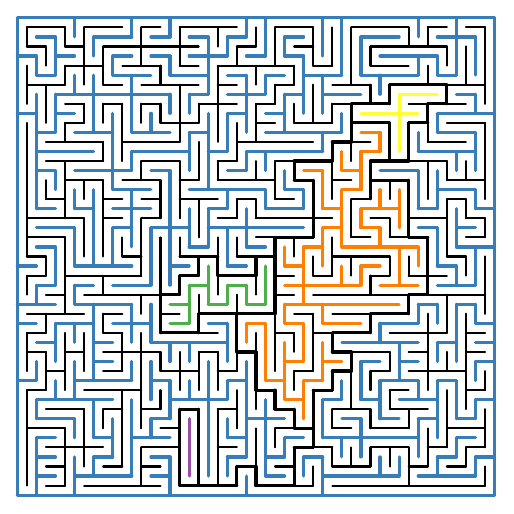}
\caption{\small Duality between the probability distributions on spanning forests and spanning connected subgraphs, illustrated here in the case $k=4$ for a square grid graph.}
\label{fig:duality}
\end{figure}

The matroidal duality between connected spanning subgraphs and spanning forests (see  Section~\ref{sec:duality-mat}) holds further in a topological sense in the case where $\rG$ is the $1$-dimensional skeleton of a $2$-dimensional complex. Duality is then given by an explicit map which acts as an involution on configurations. Figure \ref{fig:duality} illustrates this fact in the case $k=4$. We discuss this important case in this short section.

\subsection{The spherical case}\label{sec:planar-duality}

Let us start by assuming that $\rG$ is a graph embedded in an oriented sphere. Let $\rG^\dagger$ be the dual graph. On a set-theoretic level, the orientation of the sphere induces a bijection between the oriented edges of $\rG$ and those of $\rG^{\dagger}$, and we denote simply by $e^{\dagger}$ the oriented edge associated to $e$. This bijection determines two isomorphisms
\begin{align*}
\sharp:C_{1}(\rG,\Z)&\longrightarrow C^{1}(\rG^{\dagger},\Z)\hspace{1cm} \text{ and } \hspace{-1cm} &\flat : C^{1}(\rG,\Z)\longrightarrow C_{1}(\rG^{\dagger},\Z)\\
e &\longmapsto e^{\sharp}=(e^{\dagger})^{\star} &e^{\star} \longmapsto (e^{\star})^{\flat}=e^{\dagger}
\end{align*}
which are related by $(c^{\sharp},\alpha^{\flat})=(\alpha,c)$, for all $c\in C_{1}(\rG,\Z)$ and $\alpha\in C^{1}(\rG,\Z)$.

The first isomorphism sends $Z_{1}(\rG,\Z)$ to $B^{1}(\rG^{\dagger},\Z)$, and the second sends $B^{1}(\rG,\Z)$ to $Z_{1}(\rG^{\dagger},\Z)$. Extending the scalars from $\Z$ to $\K$ yields isomorphisms from $\Omega_{1}(\rG)$ to $ \Omega^{1}(\rG^{\dagger})$, and from $\Omega^{1}(\rG)$ to~$\Omega_{1}(\rG^{\dagger})$.

Let us consider an integer $k\geq 0$ and choose $\calc\in \Omega_{1}(\rG)^{\wedge k}$. Let us consider $\vartheta=\calc^{\sharp} \in \Omega^{1}(\rG^{\dagger})$ the tensor associated to $\calc$ by the $k$-th exterior power of the first isomorphism.

Let $F\in \cF_{k}(\rG)$ be a spanning forest of $\rG$ with $k+1$ connected components. Then $K=(F^{c})^{\dagger}$ belongs to~$\rC_{k}(\rG^{\dagger})$ and, with the notation of the previous sections, $z_{K}=\pm (b_{F})^{\flat}$. Then,
\[(\vartheta,z_{K})=\pm(\calc^{\sharp},b_{F}^{\flat})=\pm(\calc,b_{F}).\]

Let $\ul\x$ be a set of positive weights associated with the edges of our graphs. The following proposition is then a consequence of the definitions \eqref{eq:defApol}
 and \eqref{eq:defCpol} of the generating polynomials.
 
 \begin{proposition}\label{prop:planar-duality}For all $\calc\in \Omega_1(\rG)^{\wedge k}$, we have the equality of polynomials
\begin{equation}\label{eq:dualiteAC}
\rA^{(k)}_{\rG}(\calc,\ul\x)=\ul\x^{\dE} \rC^{(k)}_{\rG^{\dagger}}(\calc^{\sharp},\ul\x^{-1})\,.
\end{equation}
\end{proposition}

It follows from \eqref{eq:dualiteAC} that the determinantal measures described in Theorem \ref{thm:DPP1} on~$\cF_{k}(\rG)$ and in Theorem \ref{thm:DPP2} on $\cC_{k}(\rG^{\dagger})$ are in correspondence, via the map $S\mapsto (S^{c})^{\dagger}$, up to the replacement of the subspace~$H_\vartheta$ of $\Omega_{1}(\rG)$ by the subspace $H_\calc=H_{\vartheta}^{\sharp}$ of $\Omega^{1}(\rG^{\dagger})$, and of the positive weights by their inverses.

Equation \eqref{eq:dualiteAC} can also be seen as a consequence of Propositions~\ref{prop:egalitesconnexes} and~\ref{prop:egalitesacycliques}, using a relation of conjugation between the operators $\Box_{\calc}=bb^*$ on $\rG^{\dagger}$ and $\Delta_{\vartheta}$ on $\rG$. 

\subsection{The case of two-dimensional simplicial complexes}\label{sec:cellulation}

Assume now that $\rG= (X_0,X_1)$ is the $1$-skeleton of a simplicial complex $X=(X_{0},X_{1},X_{2})$ of dimension~$2$. Let $d':\Omega^1(X)\to \Omega^2(X)$, be the coboundary map from $1$-forms to $2$-forms. Since $d'\circ d=0$, we have $\im d \subset \ker d'$. Let $H^{1}(X,\Z)=\ker_\Z d' /\im_\Z d$. Then $\dim H^{1}(X,\Z)\otimes \K=b_1(X)$, the first Betti number of $X$. 

Applying Theorem \ref{thm:DPP1} with $H=\ker d' \simeq \im d\oplus H^{1}(X,\K)$ yields a determinantal probability measure on~$\cC_{b_1(X)}(\rG)$ supported on the set $\cH_1(X)$ of all $K\in \cC_{b_1(X)}(\rG)$ such that the natural map $Z_1(K,\Z)\to H_1(X,\Z)$ is an isomorphism.
In particular, $\cH_1(X)$, which is the support of this measure, is the set of bases of a matroid. This case was considered by Lyons in \cite[Section 3]{Lyons-Betti}, under the name~$\mathbf{P}^1$.

Note that the `planar dual' of the random subgraph in that case is the uniform spanning tree, whose law is denoted by $\mathbf{P}_1$ by Lyons in this setup.

Let $\sigma_1, \ldots, \sigma_{b_1(X)}$ be an integral basis of $H^{1}(X,\Z)$ and set $\varsigma=\sigma_1\wedge\ldots\wedge\sigma_{b_1(X)}$. Then for any $K\in \cH_1(X)$, we have $(\varsigma, z_K)=\pm 1$ and our construction yields the uniform measure on $\cH_1(X)$.

For example, take $X$ to be a $2$-cellulation of a closed surface $\Sigma$ of genus $g\ge 1$. Then $b_1(X)=2g$ and $H_1(X, \Z)\simeq H_1(\Sigma, \Z)\simeq \Z^{2g}$. Elements of the support of our measure, $\cH_1(X)$, are then sometimes called $2g$-quasitrees of the map $X$ in the combinatorics literature.\footnote{In yet other parts of the literature, $2g$-quasitrees of a map $X$ of genus $g$ are also refered to as unicellular maps, that is maps whose $1$-skeleton is a subgraph of the $1$-skeleton of the map such that the complement of this subgraph in the map $X$ is homeomorphic to a disc. We can obtain determinantal measures on these lower genera quasi-trees by taking other subspaces $H$ such that $\im d \subseteq H \subseteq \ker d'$. An special example is given in Section \ref{sec:random-g-quasitrees-periods}.} Along with lower genera quasitrees, they appear in the definition of the Bollob\'as--Riordan polynomial of the cellulation~\cite{CKS}, which is known to fit in the general framework of Tutte polynomials of matroids~\cite{Moffatt-Smith}. 

\begin{figure}[!ht]
\centering
\includegraphics[width=.5\textwidth]{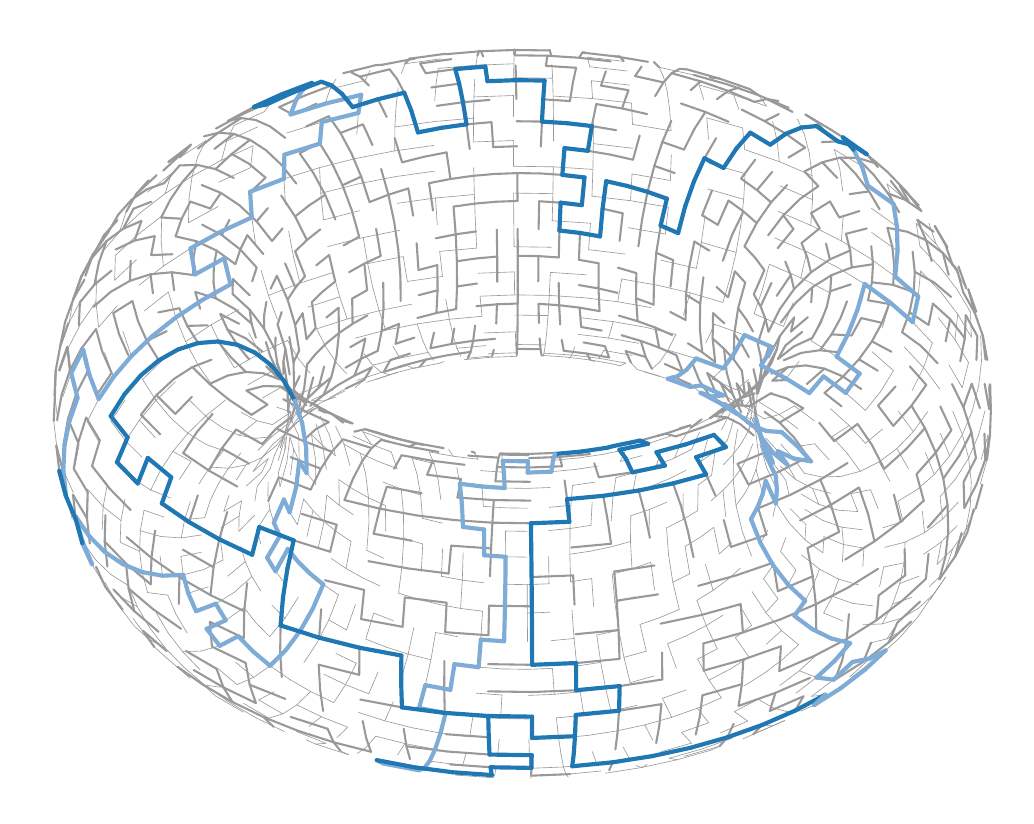}
\caption{\small A uniform element of $\cH_1(X)\subset \cC_{2}(X_{1})$ where $X=(X_0,X_1,X_2)$ is a quadrangulation of the flat torus of genus $1$.}
\label{fig:torus}
\end{figure}

\begin{figure}[!ht]
\centering
\includegraphics[width=.6\textwidth]{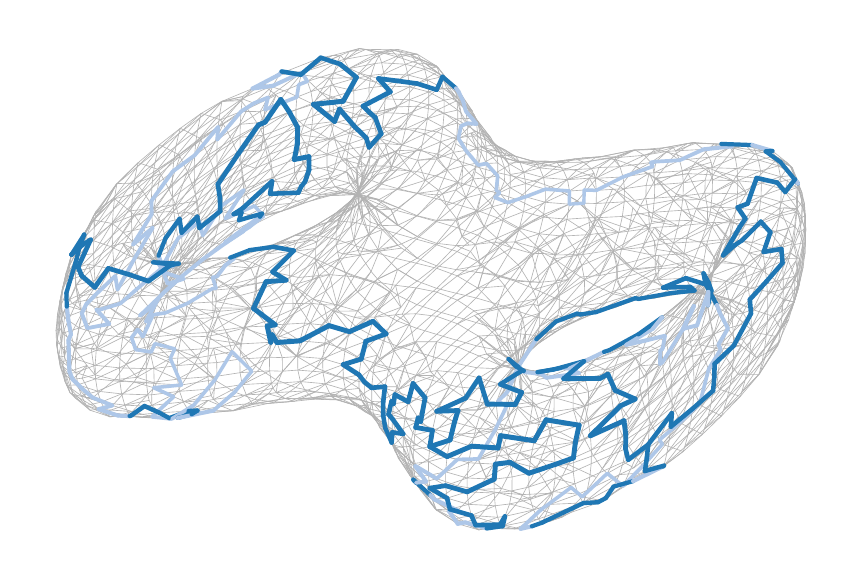}
\caption{\small The $2$-core of a uniform element of $\cH_1(X)\subset \cC_{4}(X_{1})$ where $X=(X_0,X_1,X_2)$ is a discretization of a surface of genus $2$.}
\label{fig:torus2}
\end{figure}

A natural question is to compute the distribution on $\Z^{2g}$ induced by the uniform measure on~$\cH_1(X)$.
See Figure~\ref{fig:torus} and Figure \ref{fig:torus2} for an illustration in genus $g=1$ and $g=2$, respectively. 

Our construction can be generalized to higher dimensional complexes, following \cite{Lyons-Betti}, see the brief discussion in Section \ref{sec:examples}. The above-mentioned duality in the planar case then is the Poincar\'e duality $X\mapsto X^\dagger$ between~$k$-cells and $(d-k)$-cells in a $d$-dimensional complex.

\subsection{Choice of convention}

Note that we could have given alternative definitions for the polynomials $\rA_\rG(\vartheta, \ul\x)$, replacing the terms $\ul\x^F$ by $\ul\x^{F^c}$, like in the definition of Symanzik polynomials (Section \ref{sec:polynomials-higher}). This would have had the advantage of simplifying certain formulas, notably those that make use of the planar duality, such as \eqref{eq:dualiteAC}. However, we have chosen to endow only cochains (that is,~$0$-forms and $1$-forms) with inner products (see Section \ref{sec:conventions}) and we also prefer to define determinantal processes on the set of edges with respect to subspaces of $\Omega^1(\rG)$ (not of $\Omega_1(\rG)$), so as to compare more easily with the classical cases of the uniform spanning tree. Defining the polynomials so as they would be generating functions for these determinantal probability measures (and not their dual determinantal probability measures) was a further argument in favor of this choice. This choice will also be apparent in the way we treat with the matroid generalisation in Section \ref{sec:matroid-setting}.


\section{Multivariate homogeneous real stable polynomials}\label{sec:polynomials}

In this short section, we take a closer look at the multivariate polynomials $\rC_{k}(\vartheta, \ul\x)$ and~$\rA_{k}(\calc, \ul\x)$, derive some of their properties from their relation with determinantal random subgraphs, and emphasize their link with the Symanzik polynomials of theoretical physics.

\subsection{Real stability}

A multivariate polynomial with real coefficients is called {{\em stable} if it does not vanish when all the variables have strictly positive imaginary part; see e.g. \cite[Definition~2.9]{Borcea-Branden-Liggett}. The stability of real multivariate polynomials is a generalization of the property for a univariate polynomial that all its roots be real.

Examples of real multivariate polynomials are given by the generating polynomials of probability measures on the power set of a finite set. The class of probability measures for which the corresponding generating polynomial is stable is called the class of \emph{strongly Rayleigh probability measures} (see \cite[Definition 2.10]{Borcea-Branden-Liggett}), and determinantal probability measures belong to this class by \cite[Proposition 3.5]{Borcea-Branden-Liggett}.

Since the multivariate polynomials  $\rC_{k}(\vartheta, \ul\x)$ and~$\rA_{k}(\calc, \ul\x)$ are, up to a normalizing constant, the generating functions of determinantal probability measures, these polynomials are real stable.

\begin{proposition}
The real multivariate polynomials $\rC_{k}(\vartheta, \ul\x)$ and~$\rA_{k}(\calc, \ul\x)$ defined in~\eqref{eq:defCpol} and~\eqref{eq:defApol}, respectively, are real stable.
\end{proposition}

\begin{proof}
By Proposition 3.5 of \cite{Borcea-Branden-Liggett}, this is a consequence of Theorems~\ref{thm:DPP1} and \ref{thm:DPP2}. Let us   briefly recall the main steps of the proof of this result. The first step consists in writing the generating function of the determinantal probability measure on a set of cardinality $n$ associated with a self-adjoint matrix $0\le K\le 1$, which according to \cite{Lyons-DPP} (see also \eqref{eq:fonctiongeneratriceDPP}) is equal to
\begin{equation}
\E\big[\ul{x}^{\X}\big]=\det\big(I_{n}+({\rm diag}(\ul{x})-I_{n})K \big).
\end{equation}
The second step is to put this expression in the `prototypical' form for stable polynomials, that is, for some appropriately chosen positive semi-definite matrices $A_1,\ldots,A_{n}$ and a self-adjoint matrix~$B$, the form
\begin{equation}
P:\underline{x}\mapsto \det\Big(B+ \sum_{i=1}^{n} x_i A_i \Big)\,.
\end{equation}
In turn, a polynomial of this form is proven, in \cite[Proposition 3.2(1)]{Borcea-Branden-Liggett}, to be either zero or stable, based on the fact (\cite[Proposition 3.1(1)]{Borcea-Branden-Liggett}) that stability is equivalent to the property that for all $\ul a\in \R_+^{n}$ and $\ul b\in \R^{n}$, the univariate polynomial 
\begin{equation}
z\mapsto P(\ul a z +\ul b)
\end{equation} 
has all its zeros on the real line.
\end{proof}

Homogeneous stable polynomials as above are a special case of Lorentzian polynomials, a family of polynomials with deep connections to matroid theory, see \cite{Branden-Huh}.

\subsection{Symanzik and Kirchhoff polynomials}\label{sec:symanzik-kirchhoff-pols}

Symanzik polynomials appear in Feynman integrals associated with finite graphs. We refer to the introduction of \cite{ABBGF} for a mathematical presentation of these integrals. Combining the (slightly modified) notations of these authors with ours, the amplitude associated with an unweighted graph $\rG$ endowed with external momenta $\vv{q}\in (\R^D)^\rV$ is the real number defined by
\begin{equation}\label{eq:feynman}
I_\rG(\vv{q})=\int_{(\R_+)^{\dE}} \exp\bigg(\! -i \frac{\Psi_{2,\vv{q}}(\ul{y})}{\Psi_{1}(\ul{y})}\bigg) \frac{d\ul{y}}{\Psi_{1}(\ul{y})^{\frac{1}{2}}}
\end{equation}
where
\[\Psi_{1}(\ul{y})=\sum_{T\in \cT(\rG)} {\ul{y}}^{\rE\setminus T}\quad\text{and}\quad\Psi_{2,\vv{q}}(\ul{y})=\sum_{F=\{T,T'\}\in \cF_2(\rG)} -\big\langle  \vv{q}(\rV(T)),\vv{q}(\rV(T'))\big\rangle \; {\ul{y}}^{\rE\setminus F}\,,\]
and $\langle \cdot,\cdot \rangle$ is a Minkowski bilinear form on $\R^D$. Here we used the notation 
\[\vv{q}(\rV(T))=\sum_{v\in \rV(T)} \vv{q}(v)\in \R^D\,.\]
As already alluded to in Section \ref{sec:pythagore-sym}, these polynomials are called the first and second Symanzik polynomials in the literature, see~\cite{Bogner-Weinzierl}.

In the remainder of this subsection and the next, when we write $\ul{x}$, we mean $\ul{y}^{-1}$. The reason for this distinction is that we think of $y_e$ as edge-lengths, or resistances, whereas $x_e$ are conductances. The discrete analogue of a continuous Laplacian is the operator $\Delta$ associated with conductances~$x_e$, when the edge lengths are $y_e$.

In physical terms, $D$ represents the dimension of space-time, so that the case $D=1$ that we will now consider seems to have little physical relevance.
When $D=1$ and considering the usual norm on $\K$ instead of a Minkowski bilinear form, taking $q\in (\ker d)^\perp=\im d^*$ and writing it~$q=d^* {\Jx} c$ for some~$c\in \Omega_1(\rG)$, we can re-express the second Symanzik polynomial in terms of the polynomial $\rA_1((c, \ul{x})$ defined in~\eqref{eq:defApol}, as
\begin{equation}\label{eq:second-sym-pol}
\Psi_{2,q}(\ul{y})=\sum_{F=\{T,T'\}\in \cF_{2}(\rG)} |q(\rV(T))|^{2}\,  \ul{y}^{\rE\setminus F}
=\ul{y}^{\rE}\rA_{1}(c, \ul{y}^{-1})\,.
\end{equation}
 
\subsection{A reparametrization}\label{sec:reparametrization}

Since the polynomials $\rA_1$ have a generalization to $\rA_k$, equation~\eqref{eq:second-sym-pol} suggests the following generalization of the above computation, based on a reparametrization of the weight of forests \eqref{eq:weight-kir} appearing in $\rA_k$. This reparametrization will be further used in Section~\ref{sec:symanzik-forests}.

Let $k\ge 2$ be an integer. Let $(c_1, \ldots, c_k)$ be elements of $\Omega_1(\rG)$ and set $q_i={\sf J}^0 \partial c_i=d^*{\Jx}c_i$ for all $i\in \{1\ldots, k\}$. For all $F\in \cF_{k}(\rG)$ we choose $T_1, \ldots, T_k$ an enumeration of the trees of $F$ except one, and consider the integral basis of~$B^{1}(\rG,F)$ consisting in the set of cuts (coboundaries) $\{\beta_i=\delta(\1_{\rV(T_i)}): 1\le i \le k\}$. Thus, for all $j\in \{1,\ldots,k\}$,
\[\big(\beta_{j}, c_i\big)=\blanx {\Jx}c_{i}, d(\1_{\rV(T_{j})})\ranx=\blanx d^{*}{\Jx}c_{i},\1_{\rV(T_{j})}\branx=\sum_{v\in \rV(T_j)} q_i(v)\]
a quantity which we denote by $q_{i}(\rV(T_{j}))$.
The weight of the forest $F$ given by~\eqref{eq:weight-kir} is then
\begin{equation}\label{eq:weight-kir-beta-q}
\big\vert \det \big((\beta_j,c_{i})\big)_{1\le i,j\le k} \big\vert^2=\big\vert \det \big( q_i(\rV(T_j)) \big)_{1\le i,j\le k}\big\vert^2\,.
\end{equation}

\subsection{Higher order Symanzik polynomials}\label{sec:polynomials-higher}
Setting ${\sf q}=(q_1,\ldots, q_k)$, we may thus consider the polynomial
\begin{equation}\label{eq:sym-pol-k}
\Psi_{k+1,{\sf q}}(\ul{y})=\sum_{F\in \cF_{k}(\rG)} \big\vert \det \big( q_i(\rV(T_j))\big)_{1\le i,j\le k}\big\vert^2 \; \ul{y}^{\rE\setminus F} \;\; \in \; \R[\ul\y]
\end{equation}
as a natural generalization of the second Symanzik polynomial \eqref{eq:second-sym-pol} to higher order $k+1\ge 3$. This polynomial is simply $\ul{y}^{\rE}\rA_{k}(\calc, \ul{y}^{-1})$ defined in \eqref{eq:defApol} above, where $\calc=c_1\wedge\ldots\wedge c_k$.

Symanzik polynomials, and their `duals', Kirchhoff polynomials, have also been generalized to higher order, and extended from graphs to matroids, by Piquerez \cite{Piquerez}, where a link to determinantal (and even hyperdeterminantal) probability measures is also briefly mentioned in the introduction. Along with the family of polynomials $\rA_{k}(\calc, \ul\x)$, another natural generalization of these polynomials is the family of polynomials $\rC_{k}(\vartheta, \ul\x)$ defined in~\eqref{eq:defCpol}.

\subsection{Ratios of Symanzik polynomials and Amini's strong stability theorem}\label{sec:ratio-omid}

The first and second Symanzik polynomials are known to have interesting analytic properties.
In particular, Omid Amini has shown in \cite[Theorem 1.1]{Amini-exchange} that the ratio of the two first Symanzik polynomials, seen as a rational function of the weights $\ul\x$, which appears in the computation of the Feynman integral~\eqref{eq:feynman}, has bounded variation at infinity. This has applications to tropical geometry~\cite{ABBGF}.

We may rewrite this ratio of polynomials, for all $q\in \im d^*={\sf J}^0(\im \partial)$\footnote{Note that $\im d^*$ is independent of $\ul \x$, although $d^*$ depends on those weights. One short proof of this fact is that $\im d^*$ is the orthogonal in $\Omega^0(\rG)$ of $\ker d$, but neither $d$, nor the inner product on $\Omega^0(\rG)$ depend on $\ul{x}$.}, using our notations and considering $c\in C_1(\rG,\K)$ such that~$q={\sf J}^0 \partial c$, setting $\ul{x}=\ul{y}^{-1}$, as
\begin{equation}\label{eq:ratio-omid}
\frac{\Psi_{2,q}(\ul{y})}{\Psi_{1}(\ul{y})}=\frac{\rA_{1}(c, \ul{y}^{-1})}{\rT(\ul{y}^{-1})}=\frac{\rA_{1}(c, \ul{x})}{\rT(\ul\x)} = \big\|\proj{\im d}({\sf J}^{1}_{\ul{x}} c)\big\|^{2}\,,
\end{equation} 
where, to prove the second equality, we used Proposition \ref{prop:egalitesacycliques} with $k=1$.

Let us define the discrete Green function
\begin{equation}
G_{\ul\x}=\left(\big(\Delta\big)_{\im d^*}^{\im d^*}\right)^{-1}\in \End(\Omega^0(\rG))
\end{equation}
to be the inverse of the compression of the Laplacian $\Delta=(d^*d)$ on the orthogonal of its kernel.
Since $\proj{\im d}=dG^{-1}_{{\ul\x}}d^*$, and since $d^*({\Jx} c)={\sf J}^0 (\partial \otimes {\rm id}_\K) c$ by \eqref{eq:partiald*}, we can simplify \eqref{eq:ratio-omid} further to the following expression.
\begin{proposition}Setting $\ul{x}=\ul{y}^{-1}$, we have
\begin{equation}\label{eq:height-pairing}
\frac{\Psi_{2,q}(\ul{y})}{\Psi_{1}(\ul{y})}=\lanx q, G_{\ul\x} q\ranx\,.
\end{equation}
\end{proposition}

Note that $\Delta f(v)=\sum_{e:\ol{e}=v} x_e(f(v)-f(w))$ for each $f\in \Omega^0(\rG)$ and vertex $v$. Moreover, since both $\im d^* = {\sf J}^0 (\im \partial \otimes \K)$ and the inner product on $\Omega^0(\rG)$ are independent of $\ul\x=\ul\y^{-1}$, the dependence in $\ul\y$ of the right-hand side of \eqref{eq:height-pairing} is only via its dependence inside $\Delta$.

The expression \eqref{eq:height-pairing} seems to be a discrete analogue of the `archimedean height' pairing (seen as a quadratic form here, evaluated on $q$ seen as a degree-zero divisor on a curve) considered in \cite{ABBGF}, in view of its expression in terms of the Green function of the Riemann surface whose `dual graph' (in the sense of algebraic geometry, not of graph theory) is $\rG$ (\cite[Lemma~6.3]{ABBGF}). This archimedean height pairing is shown in that paper to be equal to the ratio of Symanzik polynomials in a certain limit, which suggests from the above it has to do with the continuous Green function converging in that limit (where the lengths of the edges of the associated metric graphs converge to infinity) to the discrete Green function.

From \eqref{eq:height-pairing} above, we obtain the following reformulation of a theorem of Amini.

\begin{theorem}[O. Amini, \cite{Amini-exchange}] \label{thm:stablility-infinity}
For all $q\in \Omega^0(\rG)$ such that $\sum_{v\in \rV} q(v)=0$, and all collections of positive weights $\ul\y_{0}$, the rational function $F_{q}:(\R_{+}^{*})^{\dE} \to \R_+,\  \ul\y\mapsto  \lanx q, G_{\ul\y^{-1}}q \ranx_{\ul\y^{-1}}$ satisfies the relation $F_{q}(\ul\y+\ul\y_{0})-F_{q}(\ul\y)=O_{\ul\y}(1)$ as~$\ul\y\to \infty$.
\end{theorem}
\begin{proof}
We combine \eqref{eq:height-pairing} above and a special case of~\cite[Theorem 1.1]{Amini-exchange}.
\end{proof}

In view of the expression \eqref{eq:height-pairing} for the ratio of polynomials in terms of the Green function, one may wonder if there is an alternative proof of (this special case of) Amini's stability theorem based on the study of variations of the Green function when changing edge weights in the limit of zero edge weights, and if his stability result extends to other ratios of multivariate polynomials, such as the ones appearing in Propositions~\ref{prop:egalitesconnexes} and \ref{prop:egalitesacycliques}.

Let us be a bit more specific and ask a concrete question. For that matter, let us start with the following proposition. 

\begin{proposition}\label{prop:amini-k}
Let $k\ge 1$ be an integer. Let $\ul{x}=\ul{y}^{-1}$. Let ${\sf q}=(q_1,\ldots, q_k)$ be a collection of $k$ elements of $\im d^{*}$. We have
\begin{equation}\label{eq:sym-green-k}
\frac{\Psi_{k+1,{\sf q}}(\ul{y})}{\Psi_{1}(\ul{y})}=\lanx q_1\wedge \ldots \wedge q_k, (G_{\ul\x})^{\wedge k}(q_1\wedge \ldots \wedge q_k)\ranx_{\ul{x}}=\det\left[\big(\lanx q_i, G_{\ul{x}} q_j\ranx_{\ul{x}}\big)_{1\le i,j\le k}\right]\,.
\end{equation}
\end{proposition}

It would be interesting to know if this expression has a meaning in geometry, following~\cite{Zhang,Amini-exchange}, as a discrete analogue of the exterior powers of the height pairing on degree zero divisors. In particular, could it be used, combined with the convergence results of \cite{Amini-Nicolussi} to give another proof of the convergence statement in \cite{ABBGF}?

\begin{proof}
For each $i\in\{1,\ldots, k\}$, we let $c_i\in C_1(\rG,\K)$ be such that $q_i={\sf J}^0 (\partial\otimes{\rm id}_\K) c_i$. By \eqref{eq:partiald*} we thus have $q_i=d^* \Jx c_i$. 

From the remarks in Subsection \ref{sec:symanzik-kirchhoff-pols} above, recalling that ${\sf q}=(q_1,\ldots, q_k)$, and setting $\calc=c_1\wedge\ldots\wedge c_k$, we have
\begin{equation}\label{eq:psi-k-norm}
\frac{\Psi_{k+1,{\sf q}}(\ul{y})}{\Psi_{1}(\ul{y})}=\frac{\rA_{k}(\calc, \ul\x)}{\rT(\ul\x)}= \big\| (\proj{\im d}{\Jx})^{\wedge k}(\calc)\big\|_{\ul\x}^{2}\,,
\end{equation} 
where, to prove the second equality, we used Proposition \ref{prop:egalitesacycliques} with $\ul\x=\ul\y^{-1}$. 

We may now rewrite the right-hand side of \eqref{eq:psi-k-norm} as 
\begin{align*}
\big\| (\proj{\im d}{\Jx})^{\wedge k}(\calc)\big\|_{\ul\x}^{2} & = \blanx (\Jx)^{\wedge k}(\calc), (d G_{\ul{x}} d^*)^{\wedge k} (\Jx)^{\wedge k}(\calc)\branx \\
&= \blanx q_1\wedge\ldots\wedge q_k , (G_{\ul{x}})^{\wedge k}(q_1\wedge\ldots\wedge q_k)\branx \\
& = \det\left[\big(\lanx q_i, G_{\ul{x}} q_j\ranx\big)_{1\le i,j\le k}\right]\,.
\end{align*}
This concludes the proof.
\end{proof}

Following O. Amini, we may thus ask the question of the analytic behaviour at infinity of the function $F_{{\sf q}}:(\R_{+}^{*})^{\dE} \to \R_+, \ \ul\y\mapsto  \det\big(\lanx q_i, G_{\ul\y^{-1}}q_j \ranx_{\ul{y}^{-1}}\big)_{1\le i,j\le k}$. Does a result similar to Theorem~\ref{thm:stablility-infinity} hold?

\subsection{A corollary}

Let us record the following consequence of Proposition \ref{prop:amini-k}, rewritten, as already given in \eqref{eq:sym-green-intro}, using our notations from Section~\ref{sec:forests} rather than the ones of Symanzik polynomials.

\begin{proposition}\label{prop:A-T-G}
For all integers $k\ge 2$ and $c_1, \ldots, c_k\in C_1(\rG,\K)$, setting $\calc=c_1\wedge\ldots\wedge c_k$, we have
\begin{equation}\label{eq:sym-green}
\frac{\rA_{k}(\calc, \ul\x)}{\rT(\ul\x)}= \det\left[\big(\lanx q_i, G_{\ul{x}} q_j\ranx\big)_{1\le i,j\le k}\right]\,,
\end{equation}
where $q_i=d^* \Jx c_i$ for all $i\in\{1,\ldots, k\}$. 
\end{proposition}

In particular, combining \eqref{eq:sym-green} for $k=2$ and for other values of $k\ge 3$, we obtain the following corollary.

\begin{corollary}
Let $k\ge 2$ be an integer. For all $c_1, \ldots, c_k\in C_1(\rG,\K)$, we have
\begin{equation}
\frac{\rA_{k}(c_1\wedge\ldots\wedge c_k, \ul\x)}{\rT(\ul\x)}= \det\left[\left( \frac{\rA_{2}(c_i\wedge c_j, \ul\x)}{\rT(\ul\x)} \right)_{1\le i,j\le k}\right]\,.
\end{equation}
\end{corollary}

Since \eqref{eq:sym-green} can be seen as a quadratic form in $q_1\wedge\ldots\wedge q_k$, we can, by polarisation, obtain the equality of the corresponding bilinear forms. That is, given $p_1,\ldots, p_k\in \im d^*$, we have

\begin{equation}\label{eq:bilinear-form}
\sum_{F\in \cF_{k}(\rG)}\det\big(\overline{q_{i}}(V_{j})\big)_{1\leq i,j\leq k} \det\big(p_{i}(V_{j})\big)_{1\leq i,j\leq k} \ \ul x^{F}= T_{\rG}(\ul{x}) \;\det\left[\big(\lanx q_i, G_{\ul{x}} p_j\ranx\big)_{1\le i,j\le k}\right]\,.
\end{equation}


\section{Spanning forests }\label{sec:symanzik-forests}

In this section, we take a closer look at the random spanning forests considered in Theorem~\ref{thm:DPP2} after observing that the weight of a forest, originally defined by \eqref{eq:weight-kir}, can be described in terms of zero-mean functions on vertices, according to \eqref{eq:weight-kir-beta-q}.

\subsection{Intersection of hyperplanes and random spanning forests}

Recall that $\im d^*=(\ker d)^\perp$ and that $\ker d$ is the space of constant functions on $\rV$ (since $\rG$ is assumed to be connected). Hence $\im d^*$ is the space of functions of zero average.

For ${\sf q}=(q_1,\ldots, q_k)\in(\im d^*)^k$, define the subspace $H_{{\sf q}}$ of $\im d$ by 
\begin{equation}\label{eq:Hq}
H_{{\sf q}}= \{df : f\in \Omega^0(\rG), \forall i\in \{1,\ldots, k\}, \langle q_i, f\rangle=0\}\,.
\end{equation}

\begin{lemma}\label{lem:Hc-Hq}
Let $k\ge 2$ be an integer. Let $(c_1, \ldots, c_k)$ be elements of $\Omega_1(\rG)$ and set $q_i={\sf J}^0 \partial c_i=d^*{\Jx}c_i$ for all $i\in \{1\ldots, k\}$. 
We have
\[\im d \cap \Vect(c_1,\ldots, c_k)^\circ = H_{{\sf q}}\,.\]
\end{lemma}

\begin{proof}For all $f\in \Omega^0(\rG)$ and $i\in \{1,\ldots, k\}$, we have 
$\big(df, c_i\big)=\blanx {\Jx}c_{i}, df\ranx=\blanx d^{*}{\Jx}c_{i},f \branx=\blanx q_i, f \branx$. The result follows.
\end{proof}

\begin{lemma}
The codimension of $H_{{\sf q}}$ in $\im d$ is $\dim \Vect(q_1,\ldots, q_k)$. In particular, this codimension is $k$ if and only if $(q_1,\ldots, q_k)$ are linearly independent. 
\end{lemma}
\begin{proof}
Since $d$ is $1$-to-$1$ on $\im d^*$ and since $H_{{\sf q}}=\{f\in \im d^*: \forall i\in \{1,\ldots, k\},  \langle f,q_i\rangle=0\}$ is an intersection of $k$ hyperplanes, the codimension of $H_{{\sf q}}$ in $\im d$ is $\dim \Vect(q_1,\ldots, q_k)$. 
\end{proof}

In view of the two above lemmas and in light of the reparametrization given in Section \ref{sec:reparametrization}, Theorem \ref{thm:DPP2} is equivalent to the following statement. 

\begin{theorem}\label{thm:sym-forests}
Let ${\sf q}=(q_1,\ldots, q_k)$ be linearly independent elements of $\im d^*$.
The random spanning forest which assigns to any $F\in \cF_k(\rG)$ a probability proportional to $\big\vert \det \big( q_i(\rV(T_j)) \big)_{1\le i,j\le k}\big\vert^2 \ul{x}^F$ is determinantal associated with the subspace $H_{{\sf q}}$ and the orthonormal basis $(e^{\star}/\sqrt{x_{e}})_{e\in \dE}$ of $\Omega^{1}(\rG)$. 
\end{theorem}

In view of Section \ref{sec:polynomials-higher}, we call this random forest, the \emph{Symanzik spanning forest} associated with the `external momenta'~${\sf q}$.

\subsection{Rooted spanning forests}\label{sec:rooted-spanning-forests}

Let us explain how specializing the subspace $H\subseteq\Omega^1(\rG)$ in Theorem~\ref{thm:DPP2} yields a model of random rooted spanning forest (the uniform spanning tree with Dirichlet boundary conditions at a fixed number of vertices), with which we can recover the classical model of massive spanning forests which has attracted considerable attention in the literature (see e.g. \cite{BdTR, Avena-Gaudilliere,Kenyon-massive}).

For that matter, we specialize the `external momenta' ${\sf q}$ in Theorem \ref{thm:sym-forests} to zero mean functions supported on two vertices only.

\subsubsection{Vertices as parameters}

Let us consider a collection of $k+1$ disjoint vertices ${\sf v}=(v_0,v_1,\ldots, v_k)$ and for any $F\in \cF_{k}(\rG)$, and for all $i\in\{0,\ldots, k\}$, let us call $T_i$ the connected component of $F$ containing~$v_i$, and let $V_i$ be its vertex-set. 

For $i\in\{1,\ldots, k\}$, define the function $q_i=\1_{v_i}-\1_{v_0}\in \im d^*$. Recall the notation $q_i(V_j)=\sum_{v\in V_j} q_i(v)$.

\begin{lemma}\label{lem:separating-forest}
The quantity $\det(q_i(V_j))_{1\le i,j\le k}$ is nonzero, equal to $\pm 1$, if and only if there is a bijection~$\sigma$ of $\{1,\ldots, k\}$ such that for all $i\in \{1,\ldots, k\}$, we have $v_{i}\in V_{\sigma(i)}$.
\end{lemma}

Note that spanning forests having a non-vanishing weight in the above lemma can be seen as the set of $k+1$-component \emph{rooted} spanning forests (that is, spanning forests with a single distinguished root-vertex per connected component), with roots $(v_i$, $i\in \{0,\ldots, k\})$. Let us denote by $\cF_k(\rG,\underline{v})$ this set.

\begin{proof}
Let us first show that the condition is necessary. For that matter, observe that if the determinant is nonzero, then each row and each column must be nonzero. This implies that for all~$i\in\{1,\ldots, k\}$, there is $j\in\{1,\ldots, k\}$ such that $v_i\in V_j$, and reciprocally that for all $j\in\{1,\ldots, k\}$, there is $i\in\{1,\ldots, k\}$ such that $v_i\in V_j$. Hence there is a bijection $\sigma$ such that for all $i\in \{1,\ldots, k\}$, $v_i\in V_{\sigma_i}$. This shows the necessity of the assumption.

Conversely, suppose we have such a bijection $\sigma$. Observe that we then have $q_i(V_{\sigma(i)})=1$ and $q_i(V_j)=0$ if $j\ne \sigma(i)$, which shows that the matrix of which we are taking the determinant is the permutation matrix associated with $\sigma$. In particular, this determinant is nonzero and equal to $\pm 1$. This concludes the proof.
\end{proof}

\subsubsection{Rooted spanning forests}\label{sec:rooted-sf}

For each $i\in \{1,\ldots, k\}$, let $c_i$ be a path from $v_0$ to $v_i$. We can for example pick a spanning tree of $\rG$ and look at the subtree connecting the $k+1$ vertices $v_0, \ldots, v_k$ it induces, and consider for all $i\in\{1,\ldots, k\}$, the shortest path $c_i$ in this subtree connecting $v_0$ to $v_i$. 

Let us define the subspace $H_{{\sf v}}$ of $\im d$ defined by  
\begin{equation}\label{eq:Hv}
H_{{\sf v}}=\{df:f\in \Omega^0(\rG), f(v_0)=\cdots=f(v_k)\}\,.
\end{equation}
Note from \eqref{eq:Hq} that $H_{{\sf v}}=H_{{\sf q}}$.

\begin{lemma}
We have 
\[\im d\cap \Vect(c_1,\ldots, c_k)^\circ=H_{{\sf v}}\,.\]
\end{lemma}
\begin{proof}This is a special case of Lemma \ref{lem:Hc-Hq} above. However, we give a more pedestrian and `graphical' proof.
Let $f\in\Omega^0(\rG)$. Note that for any path $c$ connecting vertex $u$ to vertex $v$, the result of a telescoping sum gives $df(c)=f(v)-f(u)$. Hence, $df(c_i)=0$ for all $1\le i\le k$, if and only if $f(v_i)=f(v_0)$ for all $1\le i \le k$. This concludes the proof.
\end{proof}

We hence deduce the following consequence of Theorem~\ref{thm:sym-forests}.

\begin{theorem}
For all $k+1$ disjoint vertices $v_0,\ldots, v_k$, the probability measure on $\cF_k(\rG,\ul{v})$ which to any $F$ assigns a probability proportional to $\ul{x}^F$ is the determinantal probability measure on $\rE$ associated with the subspace $H_{{\sf v}}\subseteq \Omega^1(\rG)$ and the orthonormal basis $(e^{\star}/\sqrt{x_{e}})_{e\in \dE}$ of $\Omega^{1}(\rG)$.
\end{theorem}

Therefore, writing $\calc=c_1\wedge\ldots\wedge c_k$ and recalling the definition of $\rA_{k}(\calc, \ul{x})$ in \eqref{eq:defApol}, we find
\begin{equation}
\rA_{k}(\calc, \ul{x})=\sum_{F\in \cF_k(\rG,\ul{v})}  \ul{x}^F\,.
\end{equation}

This identity gives the evaluation at this special choice of ${\sf q}$ of the higher order Symanzik polynomials defined in \ref{sec:polynomials-higher}: $\Psi_{k+1,{\sf q}}(\ul{\y})=\sum_{F\in \cF_k(\rG,\ul{v})}  \ul{y}^{\rE\setminus F}$.

\subsubsection{Massive spanning forests}

Let $m:\rV\to \R_{+}^*$ be a positive nonvanishing function on $\rV$ (often called a `mass' in the literature). The model of massive spanning forests is the restriction of the random spanning tree measure on an augmented graph, obtained by adding to $\rG$ a sink vertex and connecting all vertices of $\rG$ to it, and assigning these edges a weight $m(v)$ per vertex $v$. 

It gives to any rooted spanning forest $F=(T_0,\ldots, T_k) \in \bigsqcup_{\ell=0}^\infty \cF_\ell(\rG)$ the weight 
\[\prod_{i=0}^k m(v_i) \ul{x}^F\,.\]

In what follows, we keep the notation $\Delta$ to denote the matrix of the discrete Laplacian $\Delta$ in the canonical basis of~$\Omega^0(\rG)$ (associated with any fixed ordering of $\rV$). By the classical matrix-tree theorem, the partition function of this model is given by 
\[\det\left(\diag(m(v):v\in \rV) + \Delta\right) = \sum_{k=1}^\infty \sum_{(v_0,\ldots, v_k)}\sum_{F\in \cF(\rG, \ul{v})} \prod_{i=0}^k m(v_i) \ul{x}^F\]
so that, from the previous section, we find the identity
\begin{equation}\label{eq:delta-Apol}
\det\left(\diag (m(v):v\in \rV) + \Delta \right) =\sum_{k\ge 0}\sum_{(v_0,\ldots, v_k)} \rA_{k}(\calc_{\ul{v}}, \ul{x}) \prod_{i=0}^k m(v_i) \,,
\end{equation}
where $\calc_{\ul{v}}$ is the wedge-product of any choice of paths $c_i$ connecting $v_0$ to $v_i$ for $i\in \{1,\ldots, k\}$.

\subsubsection{Identity between two Green functions}

By identification of the monomials in both sides of~\eqref{eq:delta-Apol}, we thus find 
\begin{equation}\label{eq:Apol-delta}
 \rA_{k}(\calc, \ul{x})=\det \Delta_{\widehat{\{v_0,v_1,\ldots, v_k\}}}^{\widehat{\{v_0,v_1,\ldots, v_k\}}}\,.
\end{equation}
Interestingly, this quantity is thus invariant under reordering of the $k+1$ vertices, and independent of the choice of paths $c_i$. 

Dividing the right-hand side of \eqref{eq:Apol-delta} by $\det \Delta_{\widehat{v_0}}^{\widehat{v_0}}$ this is equal, by Cram\'er's formula for the entries of the inverse of a matrix, as well as Jacobi's complementary minor formula, to 
\begin{equation}
\det \left[(G^{(v_0)}_{v_i,v_j})_{1\le i,j\le k}\right]\,,
\end{equation}
where we introduced the notation for the Green function with Dirichlet boundary condition at $s$ defined by 
\begin{equation}
G^{(s)}_{u,v}=\frac{\det \Delta_{\widehat{s,u}}^{\widehat{s,v}}}{\det \Delta_{\widehat{s}}^{\widehat{s}}}\,.
\end{equation}

Comparing with Proposition \ref{prop:amini-k}, which gives an alternative expression for $\rA(\calc,\ul{x})$, we thus find 
\begin{equation}
\det \left[(\lanx q_i, G_{\ul{x}} q_j\ranx)_{1\le i,j\le k}\right]
=\det \left[(G^{(v_0)}_{v_i,v_j})_{1\le i,j\le k}\right]\,.
\end{equation}
In the case $k=1$, it is equivalent to the following lemma. 
\begin{lemma}
For all $s\in \rV$, and vertex $v\in \rV$, we have 
\begin{equation}
\lanx \1_v-\1_s,G_{\ul{x}} (\1_v-\1_s)\ranx = G^{(s)}_{v,v}\,.
\end{equation}
\end{lemma} 

\subsection{The Laplacian spanning forest}\label{sec:Laplacian-forest}

Several models of random spanning forests appear in the literature. One is the \emph{massive} (or rooted) spanning forest, discussed in Section \ref{sec:rooted-spanning-forests}. It corresponds to a random spanning tree on an augmented graph where a new `sink-vertex' is connected to all original vertices, via a bijection removing edges incident to the new vertex (the rooted model corresponds to conditioning the tree to contain certain edges connected to the sink-vertex; the massive spanning forest is the unconditional measure).\footnote{Another, more recent model is the \emph{arboreal gas}, an unrooted forest model derived from percolation conditioned to yield an acyclic subgraph. While the massive spanning forest model inherits many properties from the uniform spanning tree, the arboreal gas remains less understood (see~\cite{HH}).}

In this section, extending the setup from Section \ref{sec:forests}, we define the \emph{Laplacian spanning forest}. Unlike the family of determinantal measures with fixed component numbers discussed earlier where a choice is required\footnote{One might ask which of these measures is optimal—for instance, in terms of maximizing Shannon entropy.}, this spanning forest is canonically determined from the graph's structure, without further choices. The resulting distribution is a mixture of the fixed-component measures, and the number of components has the law of a sum of independent Bernoulli random variables.

In order to introduce this model, we need to introduce a few more operators on the graph.

\subsubsection{Various Laplacian operators and their spectra}

Let us define the pointwise multiplication operator 
\begin{equation}
\begin{split}
\chi & : \Omega^0(\rG) \rightarrow \Omega^0(\rG) \\
& \quad f \mapsto \Big\{v\mapsto \Big(\sum_{e: \ol{e}=v} x_e\Big) f (v)\Big\} \,.
\end{split}
\end{equation}
With this notation at hand, let us now define an operator on $\Omega^1(\rG)$ by
\begin{equation}\label{eq:k-op}
\k=\frac{1}{2} \left( d \circ \chi^{-1}\circ d^*\right) \in \End(\Omega^1(\rG))\,,
\end{equation}
whose kernel is $\ker d^*$ and range $\im d$. This is a self-adjoint operator, which is a sort of normalized Laplacian on $1$-forms. 
In the following, we let $n=\dim \im d=\vert \rV\vert-1$, and let 
\begin{equation}
0=\mu_{-(m-n-1)}=\ldots=\mu_0< \mu_1\le \ldots \le\mu_n\,,
\end{equation}
be the eigenvalues of $\k$ ordered from smallest to largest (we used the notation $m=\vert \rE^+\vert$ for the number of edges).

For the sake of the proofs below, let us in addition define the twisted Laplacian (sometimes called the signless Laplacian) to be the operator, defined, for all $f\in \Omega^0(\rG)$ and $v\in \rV$, by
\begin{equation}
\Delta_{(-1)}f(v)=\sum_{e\in \rE: \ol{e}=v} x_e \left(f(v)+f(\ol{e})\right)\,.
\end{equation}
In the language of \cite{Kenyon, Kassel-ESAIM, KL1, KL2, KL2C, KL3, KL6}, this is the covariant Laplacian associated to a connection on the trivial vector bundle of rank $1$ with connection equal to $-1$ on all edges. It is a self-adjoint nonnegative operator on $\Omega^0(\rG)$. 

We finally introduce the normalized symmetric Laplacian 
\begin{equation}
\cL=\chi^{-1/2} d^* d \chi^{-1/2}\in \End(\Omega^0(\rG))\,. 
\end{equation}
Its kernel is of rank $1$ (when $\rG$ is connected) and the rank of the operator $\cL$ is $n=\dim \im d$.
This operator is readily self-adjoint and nonnegative, and we let 
\begin{equation}
0=\lambda_0<\lambda_1\le \ldots \le \lambda_n
\end{equation}
be the eigenvalues of $\cL$ counted with multiplicity, and ordered from smallest to largest.

\begin{lemma}
We have $\chi^{1/2}(2\id-\cL)\chi^{1/2}=\Delta_{(-1)}$.
\end{lemma}
\begin{proof}
We have $2\chi-d^*d=\Delta_{(-1)}$ and the result follows.
\end{proof}

Since $\cL$ and $\Delta_{(-1)}$ are both self-adjoint nonnegative operators, it follows that $\cL$ is bounded above by the twice the identity. Hence the spectrum of $\cL$ lies in $[0,2]$.

\begin{lemma}\label{lem:lambda-mu}
The eigenvectors of $\k$ in $\im d$ are mapped to those of $\cL$ on $\chi^{-1/2}\im d^*$ by the map $\chi^{-1/2}d^*$. Moreover, for all $i\in \{1,\ldots, n\}$, we have $\lambda_i=2\mu_i$.
\end{lemma}
\begin{proof}
Let $\phi$ be an eigenvector of $\k$ associated to an eigenvalue $\mu \ne 0$. This means that
$\k \phi=\mu \phi$ so that if we set $f=\chi^{-1/2} d^* \phi$, we have 
\begin{align*}
\cL f  &= \chi^{-1/2} d^*d \chi^{-1/2}\chi^{-1/2} d^* \phi \\
& = \chi^{-1/2} d^* (2 \k) \phi \\
& = \chi^{-1/2} d^* (2 \mu \phi) = 2 \mu f\,.
\end{align*} 
Hence $2\mu$ is a nonzero eigenvalue in $\cL$ since $f\ne 0$ (because $\phi\notin \ker \k=\ker d^*$). Since $\chi^{-1/2} d^*$ is $1$-to-$1$ on $\im d$, by taking $\phi$ spanning a basis of $\im d$ consisting in eigenvectors of $\k$, we find $n$ corresponding eigenvectors for $\cL$. Since $n=\rank(\cL)$, we have found all nonzero eigenvalues and the claim is proved.
\end{proof}

Hence, we have proved the following.

\begin{corollary}\label{lem:k-op}
The self-adjoint operator $\k$ defined in \eqref{eq:k-op} has all its eigenvalues in $[0,1]$. Namely $0\leq \k \leq 1$.
\end{corollary}

\subsubsection{The Laplacian spanning forest}

In view of Corollary \ref{lem:k-op}, the operator $\k$ induces a determinantal random subgraph. 

Let $(\phi_i)_{1\le i\le n}$ be an orthonormal basis of $\im d$ consisting in eigenvectors of $\k$ corresponding to all its nonzero eigenvalues $\mu_i$. We have
\begin{equation}
\k= \sum_{i=1}^{n} \mu_i \proj{\K \phi_i}\,.
\end{equation}

By \cite{HKPV}, the determinantal probability measure associated to $\k$ is a mixture of the determinantal measures associated to the orthogonal projection on subspaces $H_I=\bigoplus_{i\in I} \K \phi_i$ of $\im d$ corresponding to subsets of $I\subseteq \{1,\ldots, n\}$.\footnote{If the eigenvalues have multiplicity greater than $1$, we can also use the continuous averaging explained in \cite[Section~3.3]{KL3}.}

By Theorem \ref{thm:DPP2} these projection determinantal measures are random spanning forests, hence so is the mixture. Thus, we have proved the following.

\begin{theorem}\label{thm:laplacian-sf}
The determinantal random spanning subgraph of $\rG$ induced by $\k$ and the orthonormal basis of~$\Omega^1(\rG)$ given by $(e^{\star}/\sqrt{x_{e}})_{e\in \dE}$ is almost surely a spanning forest. \end{theorem}

We call this random spanning forest, the \emph{Laplacian spanning forest}. Note that conditional on being connected, the random forest is the random spanning tree distribution of Burton and Pemantle: indeed the Laplacian spanning forest is a spanning tree if and only if it has $\vert \rV \vert-1$ edges, and conditional on this, the kernel of the determinantal measure is the orthogonal projection on $\im d$, hence the result follows from Proposition \ref{prop:lawUST}.

Spelling out that the Laplacian spanning forest distribution $\P_\k$ is the mixture of the Symanzik spanning forest distributions $\P_{H_I}$, we have, for all $F\in \cF(\rG)$:
\begin{equation}\label{eq:DPPk}
\P_{\k}(\F=F)=\sum_{\substack{I\subseteq\{1,\ldots, n\}\\\vert I\vert=\vert \rE^+(F)\vert}} \prod_{i\in I} \mu_j \prod_{i\notin I} (1-\mu_i) \P_{H_I}(\F=F)\,,
\end{equation}
In particular, the number of edges present in the forest has the distribution of a sum of independent Bernoulli random variables with parameters~$\mu_i$. The number of connected components of the spanning forest is $\vert\rV\vert-\vert \rE^+(F)\vert$, thus, for all $k\in\{0,\ldots, n\}$:
\begin{equation}\label{eq:number-edges}
\P_{\k}(\F\in\cF_k(\rG))=\sum_{\substack{I\subseteq\{1,\ldots, n\}\\ \vert I\vert=n-k}} \prod_{i\in I}\mu_i \prod_{i\notin I} (1-\mu_i) \,.
\end{equation}

Hence we have proved the following.
\begin{proposition}
The law of the number of connected components of the Laplacian spanning forest minus $1$ is that of a sum of independent Bernoulli random variables with parameters $1-\mu_i$.
\end{proposition}

In particular, we have
\begin{equation}\label{eq:number-cc}
\P_{\k}(\F\in \cT(\rG))=\prod_{i=1}^n \mu_i \quad \text{and} \quad \P_{\k}(\F\in \cF_n(\rG))=\P_{\k}(\rE^+(\F)=\varnothing)=\prod_{i=1}^n (1-\mu_i)\,.
\end{equation}

See Figure \ref{fig:flap} (left) for a sample of this distribution. Since the connected components of this forest seem rather small, it is tempting to look empirically at the process whose kernel is $\k^a$ where the exponent~$a$ is chosen between $0$ and $1$. See Figure \ref{fig:flap-a} (right).

\begin{figure}[!ht]
\centering
\includegraphics[width=6cm]{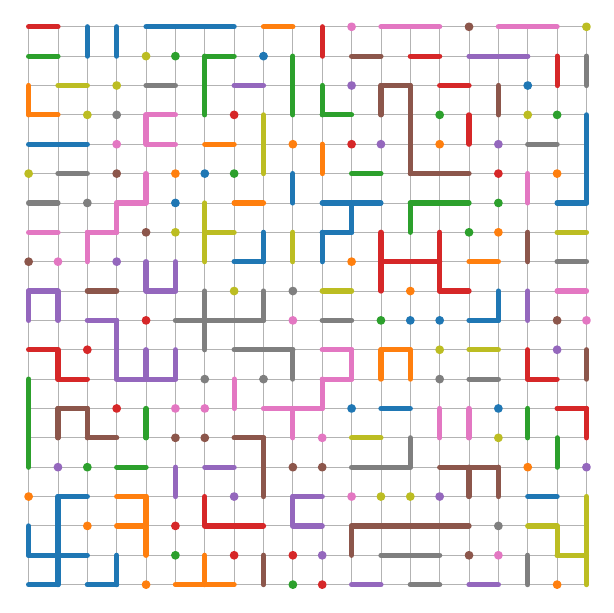} \hspace{1cm} \includegraphics[width=6cm]{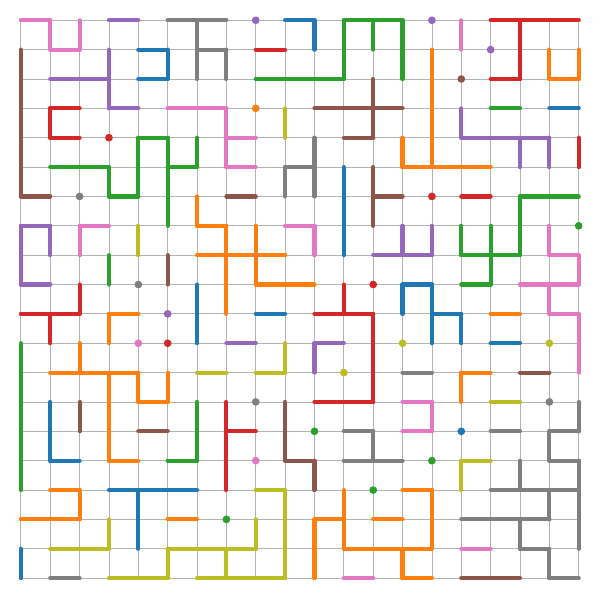}
\caption{\small Left: a sample of the Laplacian spanning forest on a finite square grid. Right: a sample of the $\k^a$-Laplacian spanning forest on a finite square grid for $a=0.3$.}
\label{fig:flap-a}\label{fig:flap}
\end{figure}

\subsubsection{Probability density}

Let us compute more precisely the probability density of the Laplacian spanning forest. We give two different formulas, Equations \eqref{eq:L-density} and \eqref{eq:density-mixture}.

\medskip

First, we know from the general theory of determinantal point processes (see for example \cite[Section 5.5]{KL3}) that 
there is a constant $\mathcal{Z}>0$ such that for any spanning forest $F\in \cF(\rG)$, we have
\begin{equation}\label{eq:L-density}
\P_{\k}(\F=F)=\frac{1}{\mathcal{Z}} \det\left({\oL}_{E^\star_F}^{E^\star_F}\right)
\end{equation}
where $\ell$ is a self-dual operator defined in terms of $\k$. 

In case $k<1$ (that is, when $\rG$ is not bipartite, by Lemma \ref{lem:bipartite}), this operator takes on a particular simple form
\begin{equation}
\oL=\k(\id-\k)^{-1}\,.
\end{equation}
Note that 
\begin{equation}
{\mathcal Z}=\det(\id + \ell)=\det(\id-k)^{-1}\,.
\end{equation}
In particular, from \eqref{eq:L-density}, we find $\P_{\k}(\rE^+(\F)=\varnothing)={\mathcal Z}^{-1}=\prod_i (1-\mu_i)$, as was noted in \eqref{eq:number-cc}.

\medskip

Second, we can compute exactly the density using our previous results. For every $i\in\{1,\ldots, n\}$, let $c_i$ be such that $\phi_i=\Jx c_i$ (this is just for notational purposes to introduce the polynomial $\rA_k(\calc,\ul{x})$). From Proposition \ref{prop:egalitesacycliques}, since $(\phi_i:1\le i\le n)$ forms an orthonormal basis of $\im d$, we have, for all $k\in \{1,\ldots, n\}$:
\begin{equation}
\rA_k(\calc,\ul{x})=\rT(\ul{x})\,,
\end{equation}  
Hence, for all $I\subseteq\{1,\ldots, n\}$, setting $H_I={\rm Vect}(\phi_{i} : i\in I)$, we have for all $F\in \cF_{n-|I|}(\rG)$:
\begin{equation}\label{eq:PHI}
\P_{H_I}(\F=F)=\frac{\vert\det(q_i(V_j(F))_{1\le i,j\le k} \vert^2}{\rA(\calc,\ul{x})} \; \ul{x}^F =\frac{\vert\det(q_i(V_j(F))_{1\le i,j\le k} \vert^2}{\rT(\ul{x})} \; \ul{x}^F\,.
\end{equation}

Thus, setting $q_i=d^*\phi_i$ for all $i\in\{1,\ldots, n\}$, we have proved the following statement.

\begin{proposition}For all $F\in \cF(\rG)$, we have
\begin{equation}\label{eq:density-mixture}
\P_{\k}(\F=F)=\frac{1}{\rT(\ul{x})}\sum_{\substack{I\subseteq\{1,\ldots, n\}\\\vert I\vert=\vert \rE^+(F)\vert}} \prod_{i\in I} \mu_i \prod_{i\notin I} (1-\mu_i)  \vert\det(q_i(V_j(F))_{1\le i,j\le k} \vert^2 \; \ul{x}^F \; \,.
\end{equation}
\end{proposition}

As follows from Lemma \ref{lem:lambda-mu}, the $q_i$ are in fact eigenvectors of $\cL$. For large random regular graphs, eigenvectors of $\cL$ are known to be `delocalized', hence heuristically, this constrasts with the very localized functions $(\1_{v_i}-\1_{v_0})$ considered for generating rooted forests in Section \ref{sec:rooted-sf} above.

\subsubsection{Gap probability and bipartiteness condition}

\begin{lemma}\label{lem:bipartite}
The operator $\k$ has eigenvalue $1$ if and only if $\rG$ is bipartite.
\end{lemma}

\begin{proof}
We note that $\k$ has eigenvalue $1$ if and only if the operator $2\chi-d^*d$ has eigenvalue $0$, but since this operator is equal to the twisted Laplacian $\Delta_{(-1)}$ with connection $-1$ on all edges, we know from~\cite{Forman} (see also~\cite{KL2, KL2C}) that it has a nontrivial kernel if and only if for all simple cycles, the product of weights is equal to $1$ (the signed graph is said to be balanced). Since for any cycle, this product is $(-1)^{\text{length of cycle}}$, $\k$ has eigenvalue $1$ if and only if all cycles are of even lengths, which means that $\rG$ is bipartite.
\end{proof}

\begin{corollary}\label{coro:vacuum-proba}
The probability that the Laplacian spanning forest has no edge is zero if and only if $\rG$ is bipartite.
\end{corollary}
\begin{proof}
From \eqref{eq:number-cc}, the probability that the Laplacian spanning forest is empty is equal to $\prod_{i=1}^n(1-\mu_i)$, where $\mu_i$ are the eigenvalues of $\k$. The product vanishes if and only if there is $i$ such that $\mu_i=1$ which occurs if and only if $\rG$ is bipartite by Lemma \ref{lem:bipartite}.
\end{proof}

This corollary highlights a difference with the uniform or massive spanning forests, as bipartiteness of the graph does not play this role in those models.

\section{Connected spanning subgraphs}\label{sec:connected-subgraphs}

In this short section which concludes the treatment of the graphical case, before moving on to the abstract matroidal generalization, we informally discuss a few examples concerning the connected case.

\subsection{Rooted spanning connected subgraphs}
In the following examples, we choose $H\supseteq \im d$ in such a way that it yields 
topologically relevant connected spanning subgraphs.

\subsubsection{Random $g$-quasitrees}\label{sec:random-g-quasitrees-periods}

The following example is inspired by the recent preprint~\cite{Lam-Lo-Yuen}, which appeared on the arXiv during the first revision of our paper. The first-named author thanks Wai Yeung Lam for a helpful conversation clarifying their work.

As in Section~\ref{sec:planar-duality}, let $\rG$ denote the $1$-skeleton of a cellular decomposition $X$ of a compact oriented surface $\Sigma$ of genus~$g\ge 1$. In this setting, a discrete analogue of the period matrix of $\Sigma$ was defined by Mercat (see~\cite{Mercat} for a survey). Unlike in the classical case where the period matrix is a $g\times g$ complex matrix, the discrete period matrix is a real $2g \times 2g$ matrix, which we denote by~$L$. The reason for this difference is that discrete holomorphic forms are best defined through their real and imaginary parts which are conjugate discrete harmonic forms on the graph and its Poincaré dual.

The main result of~\cite{Lam-Lo-Yuen}, Theorem~1.2, expresses all minors of $L$ as weighted (possibly signed) sums over $k$-quasitrees in the surface graph. In the terminology of these authors, a $k$-quasitree is a subgraph $K \in \cC_k(\rG)$ such that the induced map between homology groups $H_1(K,\Z) \to H_1(X,\Z)$ is injective.

Our observation is that \cite[Theorem~1.2]{Lam-Lo-Yuen} can be interpreted as an equality of bilinear forms on $\ext^k \Omega^1(\rG)$: one defined via the period matrix, and the other via a ratio of graph polynomials. Interestingly, the identities of our Proposition~\ref{prop:egalitesconnexes}, taken for all $k \in \{1, \ldots, 2g\}$, yield the corresponding equalities of quadratic forms, which by polarization encode the same information.

More precisely, this relationship becomes apparent for example by a reformulation of~\cite[Proposition~3.3]{Lam-Lo-Yuen}. In order to state it, let us introduce the matrix
\[
\Omega = \begin{pmatrix} 0 & I_g \\ -I_g & 0 \end{pmatrix}\,.
\]
Now, let $(c_1,\ldots, c_{2g})$ be a symplectic basis of the first homology group of $X$, that is a family of simple oriented loops in $C_1(\rG,\Z)$ which form an integral basis of $H_1(X,\Z)$ in which the matrix of the intersection form is $\Omega$. Let $(m_i:1\le i\le 2g)$ be the dual basis, that is the collection of $1$-forms such that $m_i(c_j) =\delta_{i,j}$ for all $i,j\in \{1,\ldots, 2g\}$. Then~\cite[Proposition~3.3]{Lam-Lo-Yuen} may be reformulated as an identity between a Gram matrix and the period matrix $L$:
\begin{equation}
\left(\langle m_i, \proj{\ker d^*} m_j \rangle \right)_{1 \le i, j \le 2g}=\Omega L\,.
\end{equation}
In our setup, computing principal minors of $\Omega L$ thus corresponds to applying Proposition~\ref{prop:egalitesconnexes} to subspaces of the form $H = \im d \oplus \Vect(m_{i_1}, \ldots, m_{i_k})$, for all $k \in \{1,\ldots,2g\}$ and all index sets $\{i_1 < \cdots < i_k\} \subseteq \{1,\ldots,2g\}$. Moreover, by Theorem~\ref{thm:DPP1}, the support of the determinantal probability measures corresponding to these subspaces are bases of a matroid. The special case of \cite[Corollary 1.3]{Lam-Lo-Yuen} corresponds to the case $k=1$ of \eqref{eq:partition-c}.
 
Although identities from both papers coincide in this setting, they arise from distinct perspectives. The work of~\cite{Lam-Lo-Yuen} is motivated by Lam’s earlier research on circle patterns and discrete complex structures, as well as the analogy between the period matrix and the response matrix of electrical networks, allowing the authors to resolve a question of Kenyon. Our approach, although also geometrically inspired, arose from the study of determinantal probability measures and applies in greater generality -- specializing here to the context of surface graphs.

This crossing of perspectives yields new insight. In particular, Section~7.1 of~\cite{Lam-Lo-Yuen} suggests a geometric refinement of harmonic 1-forms that we had not mentioned in the first version of our paper. Indeed it points to the importance of Lagrangian subspaces (which are of dimension~$g$) within the space of harmonic 1-forms (which has dimension~$2g$). These are the subspaces of harmonic $1$-forms obtained from the span of $g$ $1$-forms dual to a set of disjoint, homologically independent cycles. This highlights a distinguished family of determinant probability measures on random $g$-quasitrees associated with a surface discretization. These are parametrized by the compact manifold (of dimension $g(g+1)/2$) of all Lagrangian subspaces, which is called the Lagrangian Grassmannian. See Figure \ref{fig:lagrangian-g-quasitree} for an example.

\begin{figure}[!ht]
\centering
\includegraphics[width=.6\textwidth]{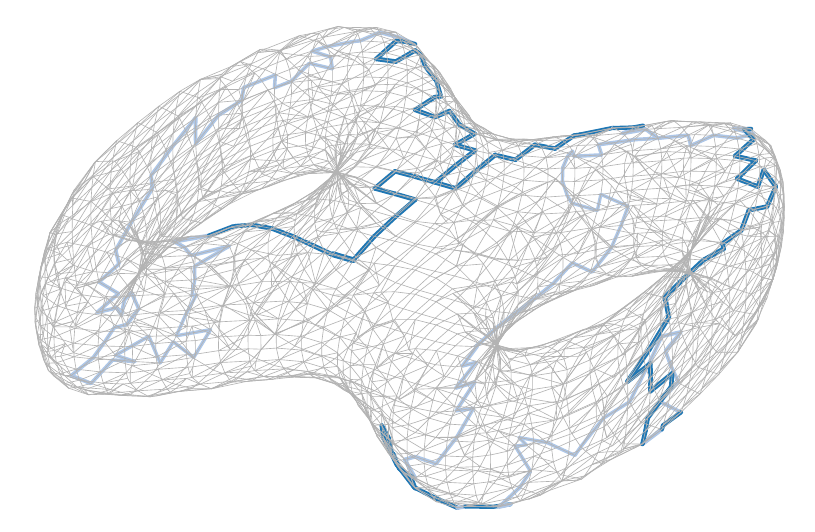}
\caption{\small The $2$-core of a determinantal $2$-quasitree associated with $1$-forms dual to the two meridians.}
\label{fig:lagrangian-g-quasitree}
\end{figure}

These measures are especially appealing since $g$-quasitrees are Poincaré dual to $g$-quasitrees on a surface of genus $g$. One may ask whether there exists a canonical choice of Lagrangian subspace yielding a determinantal measure that is self-dual under Poincaré duality. Indeed from the point of view of statistical mechanics, and in particular the Fortuin--Kasteleyn model, such a self-dual model may be critical. If it existed, such a measure could then be viewed as the natural random analogue, on a compact oriented surface of genus $g$, of the uniform spanning tree.

\subsubsection{Planar case}

In the planar case, a rooted spanning forest (like in Section \ref{sec:rooted-sf}) becomes, under planar duality, a connected spanning subgraph, whose cycles are conditioned to enclose certain faces; see Figure \ref{fig:rooted-cycle}.

\begin{figure}[!ht]
\centering
\includegraphics[width=6cm]{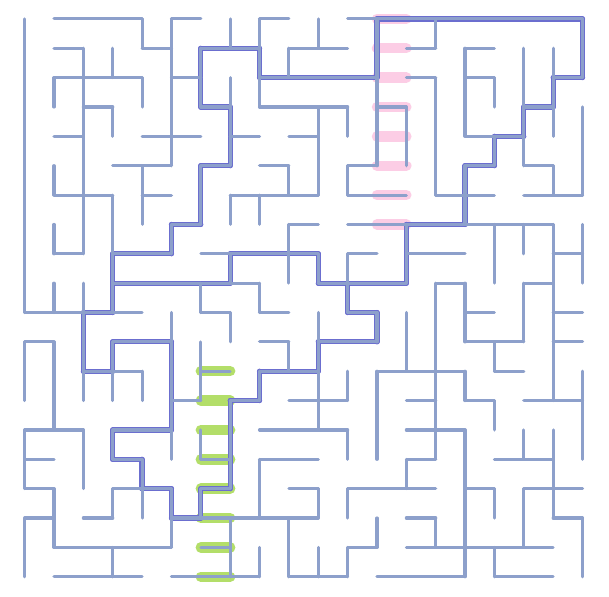}
\caption{\small A random element of $\cH_1(X)$ where $X$ is a quadrangulation of a twice-punctured disk. The colored edges represent the support of the $1$-forms $\theta_1$ (bottom zipper in green) and $\theta_2$ (top zipper in red) generating $H_1(X,\R)$}
\label{fig:rooted-cycle}
\end{figure}

\subsubsection{Three-dimensional cylinder case}

When the graph is no longer planar, there are other ways we can `root' a spanning connected subgraph. Here is an example in dimension $3$.

Let $M\subset \R^3$ be a three-dimensional filled cylindrical tube with a smaller radius cylindrical tube along its core axis removed (a `thick walled pipe section'). For concreteness, let us define it, for $r_2>r_1>0$ and $h>0$, as
\[M=\{(r,\theta,z)\in \R^3: r_1\le r\le r_2, \theta\in [0,2\pi], z\in [0,h]\}\,.\]
Let $X=(X_0,X_1,X_2,X_3)$ be a finite three-dimensional simplicial complex discretizing $M$, and let~$\rG=X_1$ be its $1$-skeleton.

The Hodge decomposition reads $\Omega^1(\rG)=\im d\oplus H_1 \oplus \im d'$, where $H_1$ is the space of harmonic $1$-forms. 

In this case, $H_1$ is $1$-dimensional. A generator of $H_1$ is a discretization of the form $d\theta$.
Let $H=\im d\oplus H_1$.

\begin{figure}[!ht]
\centering
\includegraphics[width=.5\textwidth]{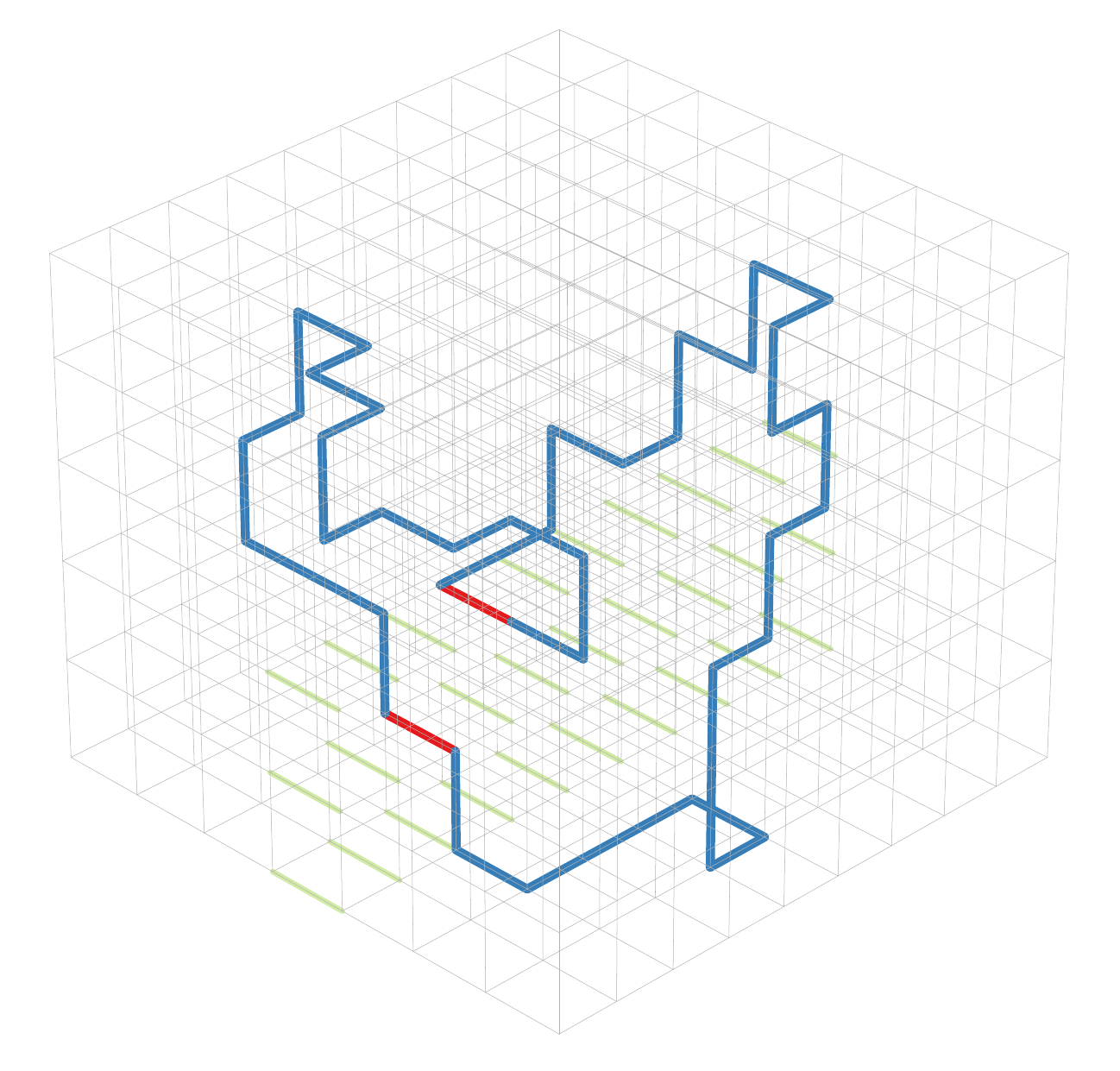}
\caption{\small Exact sample of the probability distribution from Theorem \ref{thm:cylinder} in a square grid approximation of a thick walled pipe $M=[0,10]^3\setminus [3,4]\times[0,10] \times [2,3]$. The cycle winds twice around the central tube. The $1$-form generating $H_1$ chosen for the simulation is supported on the half-wall of green edges. Red edges show places where the cycle takes one of those green edges.}
\label{fig:cylinder}
\end{figure}

Let $\cH_1(\rG)$ be the subset of $\cC_1(\rG)$ such that the unique simple cycle in the spanning subgraph is in a nontrivial homotopy class in the fundamental group $\pi_1(M)$: that is, the cycle winds around the inner tube. The winding number of a curve with respect to this inner tube is the integer representing its homology class.

In this setup, Theorem \ref{thm:DPP1} specializes to the following statement. See Figure \ref{fig:cylinder} for an illustration.

\begin{theorem}\label{thm:cylinder}
The random element of $\cC_1(\rG)$ obtained from the determinantal probability measure associated with $H$ is supported on $\cH_1(\rG)$. Under this distribution, the probability of obtaining a given spanning unicycle $F\in \cH_1(\rG)$ is proportional to $n(c)^2 \ul{x}^F$ where~$n(c)$ is the winding number of the simple cycle $c$ of $F$.
\end{theorem}

\subsubsection{Laplacian spanning connected subgraph}

We could not find a natural candidate for the analogue of the Laplacian spanning forest in the case of connected graphs (or in the general representable matroidal case, for that matter).
However, in the two-dimensional case, there is a way out.

Assume our graph $\rG$ is the $1$-skeleton of a two-dimensional simplicial complex $X=(X_0,X_1,X_2)$ with the coboundary maps summarized by the following diagram
\begin{equation}
\Omega^0(\X) 
\,
\overset{d}{\underset{d^*}{\substack{\longrightarrow\\ \longleftarrow}}} 
\,
\Omega^1(X) 
\,
\overset{d}{\underset{d^*}{\substack{\longrightarrow\\ \longleftarrow}}}
\,
\Omega^2(X)
\end{equation}
From the Hodge decomposition (see \cite[Section 1.5]{KL3}), in $\Omega^1(X)$, we have
\begin{equation}
\ker d^*=\im d^*\oplus H^1\,,
\end{equation} 
where $H^1$ are the harmonic $1$-forms, that is the kernel of 
\begin{equation}
d^*d+dd^*\,.
\end{equation} 

In order to define a random element of~$\cC(\rG)$, we need to exhibit a kernel $\k'$ defined on $\im d^*$. Taking advantage of the simplicial complex structure, let us note that the operator $\k$ (defined in~\eqref{eq:k-op}, with which we defined the Laplacian spanning forest) is a normalized `down-up' Laplacian on $1$-forms. Let us write $\k_0$ instead of $\k$ and define 
\begin{equation}
\k_2=\tfrac{1}{2}d^* \chi_2^{-1} d \in \End(\Omega^1(X))\,,
\end{equation} 
as a normalized up-down Laplacian on $1$-forms (here $\chi_2$ is the operator on $\Omega^2(X)$ which multiplies $2$-forms supported on a single face, by the sum of Euclidean weights of incident edges). Note that~$\k_0$ and~$\k_1$ are related by planar duality.

We have an analogue of the results of Section \ref{sec:Laplacian-forest}. 

\begin{proposition}
The complement $\rC$ of the determinantal random subset of $\rE^+$ associated with~$\k_2$ is a random element of~$\cC(\rG)$, which has at least $k=\dim H^1$ cycles. Its first Betti number is random and equal to $k$ plus a sum of independent Bernoulli variables of parameters $\mu_i^{(2)}$, where $(\mu_i^{(2)}:i)$ are the nonzero eigenvalues of $\k_2$.
\end{proposition}

\subsection{Deformed Hodge Laplacian and determinantal spanning subgraphs}

In the case of a two-dimensional simplicial complex, one can define a self-dual positive contraction $\k_1$ which appears as a deformation of the Hodge Laplacian $dd^*+d^*d$. 
Indeed, we can consider the sum of the normalized Laplacians in the up-down and down-up case, and add them up:
\begin{equation}
\k_1=\k_0+\k_2=\tfrac{1}{2} d \chi_{0}^{-1} d^* + \tfrac{1}{2} d^* \chi_{2}^{-1} d\,.
\end{equation}
Since both $\k_0$ and $\k_2$ are self-dual positive contractions on $\im d$ and $\im d^*$ respectively, which are orthogonal subspaces, we end up with an admissible determinantal process kernel $0\le \k_1\le 1$ on their direct sum. 

The corresponding determinantal probability distribution on $\rE^+$ is a sort of mixture of the spanning forest and spanning connected graph case. This model would deserve to be further investigated, but see Section \ref{sec:betti-fixed} for the type of subgraphs one could see arising. For this study, we may use our matroidal generalization (Theorem \ref{thm:matroid}), and we now turn to this general case, with the above examples as motivation.


\section{Measured matroids}\label{sec:matroid-setting}

In this section, we propose to revisit, in the light of the study of random subgraphs developed in this paper, the close link  
between matroids and determinantal probability measures on discrete sets, which was already made explicit and investigated by Lyons \cite[Sections 2 and 5]{Lyons-DPP}.

This section will follow and generalize to the context of linear matroids, at a much faster pace, the study that was done in Sections \ref{sec:graphs} to \ref{sec:dss}. We will show that the multilinear identities stated in Section~\ref{sec:trees} are specializations of identities that hold in general for linear matroids (see Propositions~\ref{lem:fundamental-global}, \ref{prop:sum2-matroid-global}, and \ref{prop:sumsum2}). The study of the uniform spanning tree made in Section \ref{sec:spanningtrees} and the mean projection theorem (Proposition \ref{prop:projUST}) will be paralleled by Propositions \ref{prop:wisDPP}, \ref{prop:proba-det} and \ref{coro:mean-proj-matroid}. The geometric multilinear identities proved in Section \ref{sec:gmi} will also find a more general form, namely Proposition \ref{prop:gmi-mat}. Finally, Theorem \ref{thm:DPP1} will be generalized as Theorem \ref{thm:DPPmat}.

The most significant difference between the graphical case and the general linear case that we consider in this section is the loss of an integral structure, which results in the fact that many of the matrices that we manipulate do not have integer coefficients anymore. In particular, the general expressions need to accommodate the seemingly trivial fact that the square of the determinant of an invertible matrix needs not be $1$ anymore. 

We start by briefly reviewing basic facts about matroids, then discuss determinantal measures on the set of bases of a matroid, and finally investigate the case of the union of a matroid with a uniform matroid, which corresponds to the operation that led us from spanning trees to connected spanning subgraphs, and which by a suitable operation of duality can also lead from spanning trees to spanning acyclic graphs.

\subsection{Linear matroids} For background on matroids (also known as combinatorial geometries\footnote{This terminology was proposed by Gian-Carlo Rota to replace the term matroid introduced by Hassler Whitney (1935) in his seminal study (independently carried out by Takeo Nakasawa); see \cite{Ardila}.}) we refer to the textbook \cite{Oxley} and the short introductory paper~\cite{Ardila}.

\subsubsection{Matroids and bases}
Let us recall that a matroid is, by definition, a pair $\cM=(S,\cI)$ formed by a finite set $S$ and a non-empty collection~$\cI$ of subsets of $S$, called {\em independent subsets}, such that
\begin{itemize}
\item if $J\subseteq I$ and $I\in \cI$, then $J\in \cI$,
\item for all $I, J\in \cI$ with $\vert J\vert < \vert I\vert$, there exists $i\in I\setminus J$ such that $J\cup \{i\}\in \cI$.
\end{itemize}
Elements of $\cI$ that are maximal with respect to inclusion are called {\em bases} of the matroid, and the set of all bases of $\cM$ is denoted by $\B(\cM)$. It follows from the second defining property of matroids that any two bases have the same cardinality, called the {\em rank} of $\cM$.

\subsubsection{Representing map and its kernel}
Let us fix an integer $s\geq 1$ and consider the finite ordered set $S=\{1,\ldots, s\}$. We will consider a special kind of matroid on $S$ called a {\em linear matroid}. For this, let $\K$ be $\R$ or $\C$ and let us consider an $s$-dimensional vector space $E$ over $\K$. Let $(e_1,\ldots, e_s)$ be a basis of $E$ indexed by $S$. Let $F$ be a vector space over $\K$, and let $\oR : E \to F$ be a linear map. Let us define
\[\cI=\big\{I\subseteq S : \text{the family } \{\oR(e_{i}) : i\in I\} \text{ is linearly independent in } F\big\}.\]
Then $\cM=(S,\cI)$ is a matroid on $S$, and it is said to be {\em represented} by the linear map $\oR$. Linear matroids are a fundamental example, and one of the motivations, for the notion of matroid.

A simple observation is that the subspace $\kZ=\ker \oR$ of $E$ determines entirely the matroid $\cM$. Indeed, let us define, for each subset $J\subseteq S$, the subspace $E_{J}=\Vect(e_j:j\in J)$ of $E$. Then a subset $I$ of $S$ belongs to $\cI$ if and only if $E_{I}\cap Z=\{0\}$. Moreover, a subset $T$ of $S$ is a basis of $\cM$ if and only if $E=E_{T}\oplus Z$.

The rank $n$ of the matroid $\cM$ is equal to the rank of the linear map $\oR$, so that, setting $b=\dim Z$, we have $\rank (\cM)=n=s-b=|S|-\dim \ker \oR$.

\subsubsection{Restriction of a matroid} \label{sec:restriction}
Given a subset $K$ of $S$, the matroid $\cM=(S,\cI)$ on $S$ induces a matroid $\cM_{|K}=(K,\cI_{|K})$ on~$K$ by setting $\cI_{|K}=\{I\in \cI : I\subseteq K\}$.
This construction applies to any matroid, linear or not, and is called the operation of {\em restriction}, see \cite[Section 1.3]{Oxley}.

Since the matroid $\cM$ that we consider is linear, represented by the linear map $\partial$, the restricted matroid $\cM_{|K}$ is also linear, represented by the restriction to $E_K$ of the linear map $\oR$. The kernel of this restriction is $\kZ_K=\ker (\oR_{|E_{K}})=Z  \cap E_K$.

We will exclusively apply this construction in the case where $K$ contains at least one basis of $\cM$. In this case, the set of bases of $\cM_{|K}$ is $\B(\cM_{|K})=\{T\in \B(\cM): T\subseteq K\}$. In particular, the rank of the matroid $\cM_{|K}$ is equal to $n$, the rank of $\cM$.

\subsubsection{Union of a matroid with a uniform matroid} \label{sec:matroidunion} A subset $K$ of $S$ which contains a basis of $\cM$ must have a cardinality of the form $n+k$ for some $k\in \{0,\ldots,b\}$, because $b=s-n$. For each such integer $k$, the set 
\[\B_{k}=\big\{ K \subseteq S : |K|=n+k \text{ and } \exists T\in \B(\cM), T\subseteq K\big\}\] 
of subsets of $S$ which contain a basis of $\cM$ and have cardinality $n+k$ is the set of bases of a matroid on $S$ denoted by~$\cM_k$, and that is called the {\em union} of the matroid $\cM$ with the uniform matroid of rank $k$ on $S$ \cite[Section~11.3]{Oxley}. Note that $\B_{0}=\B(\cM)$ and $\cM_{0}=\cM$.

As explained in the previous paragraph, given $k\in \{0,\ldots, s-n\}$ and $K\in \B_k$, we can consider the restricted matroid $\cM_{|K}$. In this case, the dimension of $\kZ_K$ is 
\[\dim Z_{K}=\dim E_{K}-\rank\ \cM_{|K}=(n+k)-n=k.\]

\subsubsection{Dual matroid}\label{sec:duality-mat} The set of complements in $S$ of bases of $\cM$ is the set of bases of a matroid on $S$, called the {\em dual matroid} of $\cM$, and denoted by $\cM^{*}$ (see \cite[Chapter 2]{Oxley}). The dual of a linear matroid is still linear. Starting from the matroid $\cM$, taking the dual, then the union with the uniform matroid of rank~$k$ for some $k\geq 0$ and finally taking the dual again, one finds a new matroid, namely $((\cM^{*})_{k})^{*}$, of which the bases are  the subsets of $S$ that one can obtain by removing~$k$ elements from a basis of $\cM$.

\subsubsection{A fundamental example: the circular matroid}\label{sec:dualitecirc} In the language of matroids, this paper was devoted so far to the study of the {\em circular matroid} associated to a graph $\rG=(\rV,\rE)$. An orientation of the edges being chosen, the circular matroid is a matroid on the set $\rE^{+}$, of which the independent sets are the sets of edges of acyclic subgraphs of $\rG$. It is a linear matroid, represented by the linear map $\partial : C_{1}(\rG)\to C_{0}(\rG)$. The set of bases of this matroid is the set $\cT(\rG)$ of spanning trees of $\rG$. The set of bases of the union of the circular matroid with the uniform matroid of rank $k$ is $\cC_{k}(\rG)$. The set $\cF_{k}(\rG)$ is the set of bases of the matroid obtained by the dual construction described in the previous paragraph. 

Taking advantage of this duality, we will restrict ourselves, in the present section, to the study of the operation of union of a matroid with a uniform matroid. This makes the results of this section look closer to the case of connected spanning subgraphs than to the case of acyclic spanning subgraphs, but it should be kept in mind that the acyclic case is, from the matroidal point of view, not different from the connected one.

\subsection{Fundamental bases}\label{sec:fund-bases-mat}
Let us fix once and for all in this section a linear matroid $\cM=(\cI,S)$ represented by a linear map $\partial : E \to F$ in a basis $(e_{i} : i\in S)$ of $E$.

Let us choose a basis $T$ of the matroid~$\cM$. We already said that $T$ determines a splitting 
\[E=E_{T}\oplus Z\]
and we now consider the associated a projection $\rZ_{T}:E\to Z$.

This projection can alternatively be described as follows. For any $j\in S\setminus T$, the set $T\cup \{j\}$ does not belong to $\cI$, which means by definition of $\cI$  that there exists a linear combination of $\{e_t:t\in T\}\cup\{e_j\}$ which lies in~$\kZ$. Any such linear combination must give a non-zero coefficient to~$e_{j}$, and $\rZ_{T}(e_{j})$ is the unique such linear combination for which this coefficient $1$.

The set of all $i\in S$ for which $e_i$ has a non-zero coefficient in the decomposition of $\rZ_{T}(e_{j})$ in the basis $(e_i)_{i\in S}$ is called, in the language of matroids, the \emph{fundamental circuit} associated with $T$ and~$j$, see Corollary 1.2.6 of \cite{Oxley} and the paragraph after it.

The projection $\rZ_{T}:E\to Z$ induces an isomorphism $\rZ_{T}:E_{T^{c}}\simeq E/E_{T}\to Z$, so that the family
\begin{equation} \label{eq:defZT-matroid}
\cZ_{T}=\big\{\rZ_{T}(e_{j}) : j\in T^{c}\big\}
\end{equation}
is a basis of $Z$, that we call the \emph{fundamental basis} of $\kZ$ associated with $T$.

For all $J\subseteq S$, let us denote by $\pi_{J}:E=E_{J}\oplus E_{J^{c}}\to E_{J}$ the projection parallel to $E_{J^{c}}$. Since by definition $\rZ_{T}$ vanishes on $E_{T}$, we have
\begin{equation}\label{eq:zpiidmatroid}
\rZ_{T}\circ \pi_{T^{c}}={\rm id}_{\kZ},
\end{equation}
a relation that is the analogue of \eqref{eq:zpiid}, and whose precise meaning is that the composed map
\[Z\hookrightarrow E \build{\longrightarrow}{}{\pi_{T^{c}}} E_{T^{c}} \build{\longrightarrow}{}{\rZ_{T}} Z\]
is the identity.

The basis $T$ of $\cM$ still being fixed, let us consider an integer $k\in \{0,\ldots,b\}$ and pick $K\in \B_k$, that is, a subset of $S$ of cardinality $n+k$ such that $T\subseteq K$. We can repeat the discussion above in the restricted matroid $\cM_{|K}$. The projection $E_{K}=E_{T}\oplus Z_{K}\to Z_{K}$ is simply the restriction to $E_{K}$ of the projection $\rZ_{T}:E=E_{T}\oplus Z\to Z$. Therefore, 
\begin{equation}
\cZ_{T}^{K}=\big\{\rZ_{T}(e_{j}) : j\in K\setminus T\big\}
\end{equation}
is a basis of $\kZ_K$.

\subsection{The circuit-basis identities for linear matroids}

The following proposition generalizes both Proposition \ref{lem:omid} and Proposition \ref{lem:dimo}. 
As explained at the beginning of this section, an important difference is the appearance of determinants that are no longer necessarily in $\{-1,1\}$. 

For an ordered basis $\cB$ of a vector space, we denote, as we did in Section \ref{sec:trees}, by $\det \cB$ the exterior product of the elements of this basis. For any two bases $\cB_{1}$ and $\cB_{2}$ of the same space, we denote by $\det(\cB_{1}/\cB_{2})$ the determinant of the matrix of the vectors of $\cB_{1}$ in $\cB_{2}$. Thus, the equality $\det \cB_{1}=\det(\cB_{1}/\cB_{2}) \det \cB_{2}$ holds in the top exterior power of the vector space of which~$\cB_{1}$ and~$\cB_{2}$ are bases.

For all $J=\{i_{1}<\ldots <i_{k}\} \subseteq S$, we define  $e_{J}=e_{i_{1}}\wedge \ldots \wedge e_{i_{k}}\in E^{\wedge k}$. Recall that $b$ is the dimension of $Z$.

\begin{proposition}[First circuit-basis identity]\label{lem:fundamental-global}
Let $\cZ$ be a basis of $\kZ$. Then in $\ext^{b} E$,
\begin{equation}\label{eq:omidmatroid}
\det \cZ=\sum_{T \in \B(\cM)}\det(\cZ/\cZ_T) \, e_{T^{c}}.
\end{equation}
\end{proposition}

\begin{proof} Let us write $\cZ=(z_1, \ldots, z_b)$ and decompose $z_{1}\wedge \ldots \wedge z_{b}$ on the basis $\{e_{J}: |J|=b\}$ of $E^{\wedge b}$:
\[z_{1}\wedge \ldots \wedge z_{b}=\sum_{J\subseteq S : |J|=b} a_{J} e_{J}.\]
Consider a subset $J$ of $S$ of cardinality $b$ and assume that $S\setminus J$ is not a basis of $\cM$. Then $S\setminus J$ is not independent, and $E_{S\setminus J}$ intersects $\kZ$ in a non-trivial way. This means that there exists a non-zero linear combination 
$y=u_{1}z_{1}+\ldots+u_{b}z_{b}$ that belongs to $E_{S\setminus J}$. By reordering $z_{1},\ldots,z_{b}$ if needed, let us make sure that $u_{1}\neq 0$. Then $0=(\pi_{J})^{\wedge b}(y\wedge z_{2}\wedge \ldots \wedge z_{b})=u_{1}(\pi_{J})^{\wedge b}(z_{1}\wedge \ldots \wedge z_{b})=u_{1} a_{J}e_{J}$, so that~$a_{J}=0$. 

Consider now a basis $T$ of $\cM$. Using \eqref{eq:zpiidmatroid}, we find 
\[z_{1}\wedge \ldots \wedge z_{b}=(\rZ_{T}\circ \pi_{T^{c}})^{\wedge b}(z_{1}\wedge \ldots \wedge z_{b})=\sum_{J\subseteq S : |J|=b} a_{J}\, (\rZ_{T}) ^{\wedge b} \big((\pi_{T^{c}})^{\wedge b}(e_{J})\big)\]
and the only non-zero term of the last sum is that corresponding to $J=T^{c}$, so that
\[z_{1}\wedge \ldots \wedge z_{b}=a_{T^{c}}\, (\rZ_{T})^{\wedge b}e_{T^{c}}.\]
The result follows from the observation that $(\rZ_{T})^{\wedge b}e_{T^{c}}$ is the exterior product of the elements of the basis $\cZ_{T}$. This identifies the coefficient $a_{T^{c}}$ as $\det(\cZ/\cZ_{T})$ and concludes the proof.
\end{proof}

For every basis $T$ of $\cM$, applying $(\pi_{T^{c}})^{\wedge b}$ to \eqref{eq:omidmatroid} yields
\begin{equation}\label{eq:pitcgammamatroid}
(\pi_{T^{c}})^{\wedge b} (\det\cZ)=\det(\cZ/\cZ_{T}) \, e_{T^{c}}.
\end{equation}
Since $Z^{\wedge b}$ is a line generated by $\det \cZ$, this equation and the last proposition imply that on $Z^{\wedge b}$,
\begin{equation}\label{eq:equalityextZmatroid}
\sum_{T\in \B(\cM)} (\pi_{T^{c}})^{\wedge b}={\rm id}.
\end{equation}

The following proposition generalizes Proposition~\ref{prop:sum2-sym}, in a vector space instead of a $\Z$-module. It also generalizes Proposition \ref{prop:sum2-kir}, according to the remark made in Section \ref{sec:dualitecirc}.

\begin{proposition}[Second circuit-basis identity]\label{prop:sum2-matroid-global}
Let $\cZ$ be a basis of $Z$. Then in $\big(\ext^{b}E\big)^{\otimes 2}$, we have
\[(\det \cZ)^{\otimes 2}=\sum_{T\in \B(\cM)} \det(\cZ \! /\cZ_T)^{2}\, ({\sf Z}_{T})^{\wedge k} (e_{T^{c}}) \otimes e_{T^{c}}.\]
\end{proposition}

\begin{proof} 
Let us compute the left-hand side of the equality to prove. For this, let us firstly apply \eqref{eq:equalityextZmatroid} to the second factor, which produces a sum over bases, and secondly apply  \eqref{eq:zpiidmatroid} to the first factor of each term of the sum:
\begin{align*}
(\det \cZ)^{\otimes 2}&=\sum_{T\in \B(\cM)}  \det \cZ \otimes (\pi_{T^{c}})^{\wedge b} (\det \cZ)\\
&=\sum_{T\in \B(\cM)} ({\rZ}_{T})^{\wedge b} \circ (\pi_{T^{c}})^{\wedge b} (\det \cZ) \otimes (\pi_{T^{c}})^{\wedge b}(\det \cZ).
\end{align*}
It suffices now to apply \eqref{eq:pitcgammamatroid} to each factor of each term of the sum. 
\end{proof}

Let us write down the circuit-basis identities, Propositions \ref{lem:fundamental-global} and \ref{prop:sum2-matroid-global}, applied to a restricted matroid $\cM_{|K}$.

\begin{proposition}\label{prop:sumsum2}
Let $k\in \{0,\ldots, b\}$ and $K\in \B_{k}$. Let $\cZ^K$ be a basis of $\kZ_{K}$. 

Then in $\ext^{k} E_{K}$,
\begin{equation}\label{eq:omidmatroid-local}
\det \cZ^{K}=\sum_{T \in \B(\cM), T\subseteq K}\det(\cZ^{K}/\cZ^{K}_T) \, e_{K\setminus T}.
\end{equation}

Moreover, in $\big(\ext^k E_K\big)^{\otimes 2}$, we have
\begin{equation}\label{eq:omidmatroid2-local}
(\det \cZ^{K})^{\otimes 2}=\sum_{T\in \B(\cM), T\subseteq K} \det(\cZ^K \! /\cZ^K_T)^{2}\, ({\sf Z}_{T})^{\wedge k} (e_{K\setminus T}) \otimes e_{K\setminus T}.
\end{equation}
\end{proposition}

\begin{proof} It suffices to replace $\cM$, $S$, $E$, $Z$, $b$ respectively by $\cM_{|K}$, $K$, $E_{K}$, $Z_{K}$, $k$, and to recall that the projection $\rZ_{T}:E_{K}\to Z_{K}$ associated to the restricted matroid $\cM_{|K}$ is the restriction of the projection $\rZ_{T}:E\to Z$ associated to $\cM$.
\end{proof}

\subsection{Euclidean structures}\label{sec:euclidean-case}

In the next subsection, we will define a probability measure on the set of bases of $\cM$. Recall that our fundamental data is a linear matroid $\cM=(S,\cI)$  represented by a linear map $\oR:E\to F$, where a basis $(e_{i})_{i\in S}$ of $E$ is fixed. 

This data is not sufficient to determine a probability measure, and from now on we will assume that $E^{*}$ and $F^{*}$, the dual spaces of $E$ and $F$, are inner product spaces. Moreover, we assume that the basis $(e^{\star}_{i})_{i\in S}$ of $E^{*}$, dual to the basis $(e_{i})_{i\in S}$ of $E$, is orthogonal. In other words, we assume that there exists a collection of positive {\em Euclidean weights} $\ul{x}=(\x_i:i\in S)$ such that for all $i,j\in S$,
\begin{equation}
\lanx e_{i}^{\star}, e_{j}^{\star} \ranx = x_{i} \delta_{ij}.
\end{equation}
For every subset $I$ of $S$, we will use the notation $\ul\x^{I}=\prod_{i\in I}x_{i}$ and $\ul\x^{-I}=(\ul\x^{I})^{-1}$.

Let us now consider the linear map $\oD:F^{*}\to E^{*}$, defined as the transposed map of $\oR$, and the map $\oD^{*}:E^{*}\to F^{*}$, defined as the adjoint of $\oD$.

This notation is chosen to be as close as possible to that used previously in the paper for maps between chains an cochains on graphs. For the convenience of the reader, we give a short dictionary below, see Table \ref{tab:dico}.

\begin{table}[h!]
\begin{center}
\begin{tabular}{|c|c|}
\hline
matroids & graphs \\ \hline\hline
$(S, \cI)$ & $(\rE^{+}, \cF)$\\
$\B(\cM)$, $\B_{k}(\cM)$ & $\cT(\rG)$, $\cF_{k}(\rG)$ \\
$E$, $E^{*}$ & $\Omega_{1}(\rG)$, $\Omega^{1}(\rG)$ \\
$F$, $F^{*}$ & $\Omega_{0}(\rG)$, $\Omega^{0}(\rG)$ \\
$\partial, d, d^{*}$ & $\partial, d, d^{*}$\\
\hline
\end{tabular}
\end{center}
\caption{\label{tab:dico}\small Dictionary between the graphical and matroidal situations.}
\end{table}

For each subset $J\subseteq S$, let us denote $E^{*}_{J}=\Vect(e^{\star}_i:i\in J)$. Let us observe that the annihilator of $E_{J}$ in $E^{*}$ is $E^{*}_{J^{c}}$, and that the orthogonal of $E^{*}_{J^{c}}$ is $E^{*}_{J}$. Moreover, let us note that the range of~$d$ is the annihilator of the kernel of $\partial$, that is, $\im d = Z^{\circ}$. With these observations, a subset $T$ of $S$ is a basis of $\cM$ if and only if $E=E_{T}\oplus Z$, if and only if $E^{*}=(E^{*}_{T})^{\circ}\oplus Z^{\circ}$, if and only if
\[E^{*}=((E^{*}_{T})^{\circ})^{\perp}\oplus (Z^{\circ})^{\perp}=E^{*}_{T}\oplus (\im d)^{\perp}.\]

In view of Section \ref{sec:det-toolbox}, and in particular of Proposition \ref{prop:supportdpp}, this suggests to consider the determinantal point process on $S$ associated to $\im d$ in the orthonormal basis $(e^{\star}_{i}/\sqrt{x_{i}})_{i\in S}$ of $E^{*}$. This is what we do in the next subsection.

\subsection{A determinantal probability measure on the set of bases}

The determinantal measure that we are going to construct on $\B(\cM)$ is the analogue of the `uniform' spanning tree studied in Section \ref{sec:spanningtrees}, with the important difference that in the matroidal case, the measure may not be uniform anymore, even if $x_{i}=1$ for all $i\in S$.

As a preparation, let us study how a certain determinant depends on the Euclidean weights $\ul{x}$.

\begin{lemma} \label{lem:monom-mat} Let $T\in \B(\cM)$. The determinant $\det(\oD\oD^{*})^{E^{*}_{T}}_{E^{*}_{T}}$ is equal to a constant times $\ul{x}^{T}$.
\end{lemma}

\begin{proof} Let $D$ be the matrix of $d:F^{*}\to E^{*}$ with respect to some basis of $F$ and the basis $(e^{\star}_{i} : i\in S)$ of $E^{*}$. Let $X$ be the diagonal matrix ${\rm diag}(x_{i} : i\in S)$. Let $D^{\dagger}$ be the transposed matrix of $D$. Then the matrix of $d^{*}$ with respect to these bases is $D^{\dagger}X$. Therefore,
\[\det(\oD\oD^{*})^{E^{*}_{T}}_{E^{*}_{T}}=\det (DD^{\dagger}X)_{T}^{T}= \det(DD^{\dagger})_{T}^{T}\, \ul{x}^{T},\]
as expected.
\end{proof}

To each $T\in \B(\cM)$, let us associate the weight $w(T)$ which is the real number such that
\begin{equation}\label{eq:weightsmatroid}
\det(\oD\oD^{*})^{E^{*}_{T}}_{E^{*}_{T}}=w(T)\, \ul{x}^{T}.
\end{equation}
In the case of the circular matroid on graphs, $w(T)$ is the same for all spanning trees $T$, equal to the number of vertices of the graph.

In close analogy with \eqref{eq:defTx}, let us define the generating polynomial
\begin{equation}\label{eq:mtt-matroid}
\rB(\ul{\x})=\sum_{T\in \B(\cM)} w(T)\, \ul{x}^{T}.
\end{equation}

\begin{proposition} \label{prop:wisDPP} 
Let $\X$ be a determinantal random susbset of $S$ associated with $\im \oD$ in the orthonormal basis $(e^{\star}_{i}/\sqrt{x_{i}})_{i\in S}$ of $E^{*}$. 

1. The support of the distribution of $\X$ is the set $\B(\cM)$ of bases of the matroid $\cM$.

2. For every $T\in \B(\cM)$, we have
\[\P(\X=T)=\det(\oD\oD^{*})^{E^{*}_{T}}_{E^{*}_{T}} \ \big / \det(\oD\oD^{*})^{\im d}_{\im d}=w(T)\, \ul{x}^{T} / \  \rB(\ul{x}).\]
\end{proposition}

\begin{proof} The first assertion follows from Proposition \ref{prop:supportdpp} and the fact, discussed a few lines above, that the equality $E^{*}=E^{*}_{T}\oplus \im d$ holds if and only if $T$ is a basis of $\cM$.

To prove the second assertion, let us define $\eta\in \ext^{n}E^{*}$ as the exterior product of the elements of an orthonormal basis of $\im d$. Let us consider $T\in \B(\cM)$. By definition (see \eqref{eq:defloiX}), we have
\[\P(\X=T)=\big|\big\langle \eta,e^{\star}_{T}/\sqrt{\ul x^{T}}\big\rangle\big|^{2}=\ul x^{-T} |\langle \eta,e^{\star}_{T}\rangle|^{2}.\]
Let us introduce the orthogonal projection $\proj{\im d}$ on $\im d$, and observe that $d^{*}=d^{*}\proj{\im d}$, because $\im d$ is the orthogonal of $\ker d^{*}$. We have
\[\det(\oD\oD^{*})^{E^{*}_{T}}_{E^{*}_{T}}=
\ul\x^{-T} \lanx e^{\star}_{T},(\oD\oD^{*})^{\wedge n} e^{\star}_{T}\ranx=\ul\x^{-T} \blanx(\proj{\im d})^{\wedge n} e^{\star}_{T},(\oD\oD^{*})^{\wedge n} (\proj{\im d})^{\wedge n}e^{\star}_{T}\branx.
\]
Now, $\ext^{n}\im d$ is a line, generated by the unit norm vector $\eta$, so that $(\proj{\im d})^{\wedge n}e^{\star}_{T}=\big\langle \eta, e^{\star}_{T} \big\rangle \eta$. Therefore, 
\[\det(\oD\oD^{*})^{E^{*}_{T}}_{E^{*}_{T}}=\ul x^{-T} |\langle \eta, e^{\star}_{T}\rangle|^{2} \langle \eta, (dd^{*})^{\wedge n} \eta\rangle=\big|\big\langle \eta, e^{\star}_{T}/\sqrt{\ul x^{T}}\big\rangle\big|^{2} \det(\oD\oD^{*})^{\im d}_{\im d},\]
and the result follows.
\end{proof}

It follows immediately from this proposition that
\begin{equation}\label{eq:mtt-matroid-2}
\rB(\ul{\x})=\det(\oD\oD^{*})_{\im \oD}^{\im \oD}=\det(\oD^{*}\oD)^{\im \oD^{*}}_{\im \oD^{*}},
\end{equation}
the square of the product of the non-zero singular values of $d$. In the case of the circular matroid on a graph, $\rB(\ul{x})=|\rV| \rT(\ul{x})$ and the identity above is one of the forms of the classical matrix-tree theorem.

We speak of a {\em measured matroid} to describe the situation where a linear matroid is represented by a linear map between inner product spaces, and the set of bases of this matroid is endowed with the probability measure described by Proposition \ref{prop:wisDPP}.

\subsection{A dual expression of the probability of a basis} We will now give a different expression of the distribution of the determinantal point process $\X$ introduced in the previous section. For this, let us define the antilinear isomorphism ${\sf J}_{\ul \x} : E \to E^*$ by setting, for each $i\in S$,
\[{\sf J}_{\ul \x}(e_i) =x_i^{-1} e_i^{\star}.\]
With this definition, we have for all $\alpha\in E^{*}$ and $e\in E$ the equality
\begin{equation}\label{eq:inner-product-pairing}
(\alpha, e)=\lanx {\sf J}_{\ul \x}e ,\alpha \ranx\,.
\end{equation}

Given a basis $\mathscr B=(f_{1},\ldots,f_{r})$ of a linear subspace of $E$, it will be convenient to use the notation 
\[\|\mathscr B\|^{2}=\|(\Jxm)^{\wedge r} \det \mathscr B\|^{2}=\det \big( \langle \Jxm f_{i},\Jxm f_{j}\rangle \big)_{1\leq i,j \leq r}
\]
for the Gram determinant of the image of $\mathscr B$ by $\Jxm$. The number $\|\cB\|^{2}$ can also be understood as the square of the volume of $\mathscr B$ in $E$, provided $E$ is endowed with the inner product that makes ${\sf J}_{\ul{x}}$ an isometry.

We will use the fact that if $\mathscr B_{1}$ and $\mathscr B_{2}$ are two bases of the same linear subspace of $E$, then
\begin{equation}\label{eq:detvol}
|\det (\mathscr B_{1}/\mathscr B_{2})|^{2}=\|\mathscr B_{1}\|^{2}/\|\mathscr B_{2}\|^{2}.
\end{equation}

Recall the definition \eqref{eq:defZT-matroid} of the basis $\cZ_{T}$ of $Z$ associated to a basis $T$ of $\cM$.

\begin{proposition}\label{prop:proba-det}
For all $T\in \B(\cM)$, we have
\[\P(\X=T)=\frac{\ul x^{-T^{c}}}{\|\cZ_{T}\|^{2}}.\]
\end{proposition}

\begin{proof} Let $T$ be a basis of $\cM$. 
Let $\zeta$ denote the exterior product of the elements of a basis of the orthogonal of $\im \oD$, that is, of $\ker \oD^{*}=\Jxm(Z)$. Then according to Lemma \ref{lem:orthominors},  
\[\P(\X=T)=| \lanx \zeta,e^{\star}_{S\setminus T}\ranx|^{2} \  \big/ \ \ul\x^{S\setminus T} \|\zeta\|^{2} =\ul x^{-T^{c}} \, | \lanx \zeta,e^{\star}_{T^{c}}\ranx|^{2}  \  \big/ \ \|\zeta\|^{2}.\]
Let us apply this equality with $\zeta=(\Jxm)^{\wedge b}(\det \cZ_{T})$. Using \eqref{eq:inner-product-pairing} and Proposition \ref{lem:fundamental-global}, we find
\[\blanx (\Jxm)^{\wedge b}(\det \cZ_{T}), e^{\star}_{T^{c}}\branx = (e^{\star}_{T^{c}},\det \cZ_{T})=1\]
and the result follows.
\end{proof}

This proposition, applied to the matroid $\cM$, and to a restricted matroid $\cM_{|K}$, where $K$ is a subset of $S$ containing a basis of $\cM$, implies the equalities
\begin{equation}\label{eq:functionpartition}
\sum_{T\in \B(\cM)}  \frac{\ul\x^T}{\|\cZ_{T} \|^{2}}=\ul\x^{S} \ \ \ \text{ and } \ \ \sum_{T\in \B(\cM), T\subseteq K}  \frac{\ul\x^T}{\|\cZ^{K}_{T} \|^{2}}=\ul\x^{K}.
\end{equation}
These equalities can also be deduced from the first circuit-basis identity (Propositions \ref{lem:fundamental-global} and~\ref{prop:sumsum2}) by applying a suitable exterior power of ${\sf J}_{\ul x}$ and taking square norms in the exterior algebra of $E^{*}$.

Let us conclude this section by expressing, in the present framework, the mean projection theorem (Theorem \ref{thm:projection}). For this, let us observe that for every $T\in \B(\cM)$, the endomorphism 
\begin{equation}\label{eq:projZT-matroid}
\rP_T = {\sf J}_{\ul \x} \rZ_T{\sf J}_{\ul \x}^{-1}   
\end{equation}
of $E^{*}$ is the projection on ${\Jxm}(\kZ)=\ker \oD^{*}$ parallel to ${\sf J}_{\ul \x}(E_{T})=E^{*}_{T}$.

\begin{proposition}\label{coro:mean-proj-matroid}
For each $k\in \{0,\ldots, b\}$, we have
\begin{equation}
\sum_{T\in \B(\cM)} \frac{\ul\x^{-T^{c}}}{\|\cZ_{T} \|^{2}} \; (\rP_{T})^{\wedge k}=
(\proj{\ker \oD^{*}})^{\wedge k}\,.
\end{equation}
\end{proposition}
\begin{proof}
Combine Theorem \ref{thm:projection}, Proposition \ref{prop:proba-det} and the first equality of \eqref{eq:functionpartition}.
\end{proof}

\subsection{Conditional probability measures}\label{sec:matroid-conditional}
In the case of graphs, the basis $\cZ_{T}$ associated to a spanning tree $T$ defined in Section \ref{sec:integralbases} was an integral basis of the $\Z$-module $Z_{1}(\rG)$, and the volume of this basis was the same for every $T$. This was the reason why, in the case of the circular matroid, or in general for any totally unimodular matroid, the probability measure on $\B(\cM)$ was `uniform', in the sense that it assigned to each basis $T$ a probability proportional to $\ul\x^T$.

A comparison of Proposition \ref{prop:proba-det} with its graphical analogue Proposition \ref{prop:lawUST}, or of Proposition~\ref{coro:mean-proj-matroid} with Proposition \ref{prop:projUST}, shows that the factor $\|\cZ_{T}\|^{2}$, which depends on $T$, is now ubiquitous.

Worse, in the next section, we will turn our attention to the definition of probability measures on the sets $\B_{k}$ of bases of the matroids $\cM_{k}$ defined in Section \ref{sec:matroidunion}. We will repeatedly consider subsets $K\in \B_{k}$, bases $T$ of $\cM$ included in $K$, and we will see factors $\|\cZ^{K}_{T}\|^{2}$ appear, which depend on $K$ and $T$.

A crucial help in dealing with these factors will be given by the next proposition, which computes the distribution of the determinantal subset $\X$ of $S$ conditioned on staying inside a subset $K$ of $S$, and states that it is the determinantal measure on the set of bases of the restricted matroid $\cM_{K}$ (see Section~\ref{sec:restriction}).

Note that, according to Proposition \ref{prop:wisDPP}, the condition $\P(\X\subseteq K)>0$ is equivalent to the fact that $K$ contains a basis of $\cM$, that is, that $K$ belongs to $\B_{k}$ for some integer $k$.

\begin{proposition}\label{lem:dpp-condition}
Let $\X$ be a random determinantal subset of $S$ associated with the subspace $\im \oD$ of $E^{*}$ in the basis $(e^{\star}_{i}/\sqrt{x_{i}})_{i\in S}$. Let $K$ be a subset of $S$ containing a basis of $\cM$. Then the random subset~$\X$ of $S$ conditioned on being included in $K$ is the determinantal random subset of $K$ associated with the subspace $\proj{E_{K}^*}(\im d)$ of~$E_{K}^*$ in the orthonormal basis $(e^{\star}_{i}/\sqrt{x_{i}})_{i\in K}$. 

Moreover, 
\begin{equation}\label{eq:dpp-condition-ZK}
\big(\proj{E_{K}^*}(\im d)\big)^\perp\cap E_{K}^*={\Jxm} (\kZ_K)\,.
\end{equation}
\end{proposition}

\begin{proof} The first assertion is a straightforward specialisation of Proposition \ref{lem:dpp-conditional-general}. To prove the second equality, we observe that the general identity $\proj{G}(H)=G\cap (G\cap H^\perp)^\perp$ holds for any two linear subspaces of a Euclidean space, because the adjoint of the orthogonal projection $H\to G$
is the orthogonal projection $G\to H$, and the range of the first is the orthogonal of the kernel of the second. Therefore, $(\proj{G}(H))^{\perp}\cap G=G\cap H^{\perp}$, and
$(\proj{E_{K}^*}(\im d))^\perp\cap E_{K}^*=  E_{K}^* \cap \ker d^{*}={\Jxm} (\kZ_K)$.
\end{proof}

The relation between this proposition and the volumes of bases that appear in our computations is given by the next two statements.

\begin{corollary} \label{cor:PTXK}
For all $T\in \B(\cM_K)$, we have
\begin{equation}\label{eq:dpp-condition-density}
\P\big(\X=T\, \big\vert\, \X\subseteq K\big) =\frac{\ul x^{-(K\setminus T)}}{ \| \cZ^{K}_{T}\|^{2}}.
\end{equation}
\end{corollary}

\begin{proof}
According to Proposition \ref{lem:dpp-condition} the distribution of $\X$ conditional on $\X\subseteq K$ is the distribution of the determinantal point process associated with the restricted matroid $\cM_{|K}$, which is given by Proposition \ref{prop:proba-det} applied to $\cM_{|K}$. 
\end{proof}

\begin{proposition}\label{lem:dpp-mat}
Let $K$ be a subset of $S$ containing a basis of $\cM$. Then for every $T\in \B(\cM)$ such that $T\subseteq K$,
\[\frac{\|\cZ^{K}_{T}\|^{2}}{\|\cZ_{T}\|^{2}}=\ul x^{K^{c}}\, \P(\X\subseteq K).\]
In particular, the left-hand side does not depend on $T$.
\end{proposition}

\begin{proof}
In view of Proposition \ref{prop:proba-det} and Corollary \ref{cor:PTXK}, the ratio is equal to 
\[\frac{\ul x^{S\setminus T} \, \P(\X=T)}{\ul x^{K\setminus T} \, \P(\X=T \, |\, \X\subseteq K)}
=\ul x^{S\setminus K} \, \P(\X\subseteq K),\]
as expected.
\end{proof}

Loosely speaking, the crucial point here is that the operations of restriction and conditioning (represented by the vertical arrows in the diagram below) commute with the operation of 'probabilisation', which assigns to a matroid a measured matroid (represented by the horizontal arrows).
\[
\xymatrix@=1.2cm{
\cM \ar[r]^{\text{prob.}} \ar[d]_{\rotatebox{90}{\text{\footnotesize rest.}}}& \X_{\cM} \ar[d]^{\rotatebox{270}{\text{\footnotesize cond.}}} \\
\cM_{|K} \ar[r]^{\text{prob.}} & \X_{\cM_{|K}}
}
\]

\subsection{A geometric multilinear identity}

We now combine the circuit-bases identities (Proposition \ref{prop:sumsum2}) and the mean projection theorem (Proposition \ref{coro:mean-proj-matroid}) to prove an identity which will be instrumental in our description of the determinantal probability measures of the matroids $\cM_{k}$.

\begin{proposition}\label{prop:gmi-mat}
 In $\ext^{k}E\otimes_{\R}\ext^{k} E$, the following equality holds:
\[ \sum_{\substack{ T\in \B_{0},K\in \B_{k}\\ T\subseteq  K}} \frac{\ul\x^{-T^{c}}}{\|\cZ_{T}\|^{2}} \frac{\big(({\sf J}_{\ul x})^{\wedge k} \det \mathscr Z^{K}\big)^{\otimes 2}}{\|\cZ^{K}\|^{2}}= \sum_{I\subseteq S, |I|=k} (\proj{\ker d^{*}})^{\wedge k}\big(e^{\star}_{I}/\sqrt{\ul\x^{I}}\big)\otimes e^{\star}_{I}/\sqrt{\ul\x^{I}}.\]
\end{proposition}

\begin{proof} 
Let us start from \eqref{eq:omidmatroid2-local} and express the square of the determinant of the change of bases as a quotient of squares of volumes, according to \eqref{eq:detvol}. Then, let us use Proposition \ref{lem:dpp-mat} to express the volume $\|\cZ^{K}_{T}\|^{2}$. We find
\[ \P(\X\subseteq K) \, \frac{(\det \cZ^{K})^{\otimes 2}}{\|\cZ^{K}\|^{2}}= \sum_{\substack{T\in \B(\cM), T\subseteq K}}  \frac{\ul x^{-K^{c}}  }{\| \cZ_{T} \|^{2}} \;  ({\rZ}_{T})^{\wedge k} (e_{K\setminus T}) \otimes e_{K\setminus T}.
\]

Let us apply $((\Jx)^{\wedge k})^{\otimes 2}$ to both sides of this equality and sum over all  $K\in \B_{k}$. On the left-hand side, we find, by Proposition \ref{prop:proba-det}, the left-hand side of the equality that we want to prove. On the right-hand side, thanks to \eqref{eq:projZT-matroid}, we find
\[\sum_{T\subseteq K}  \frac{\ul x^{-T^{c}}  }{\| \cZ_{T} \|^{2}} ({\rP}_{T})^{\wedge k} (e^{\star}_{K\setminus T}/\sqrt{\ul\x^{{\scriptscriptstyle K\setminus T}}}) \otimes e^{\star}_{K\setminus T} / \sqrt{\ul\x^{{\scriptscriptstyle K\setminus T}}},\]
where the sum is over all pairs formed by a  subset $K\in \B_{k}$ and a basis $T$ of $\cM_{|K}$. Let us re-index this sum as a sum over pairs formed by a basis $T$ of $\cM$ and a $k$-subset $I$ of $S$ disjoint from $T$, this set~$I$ playing the role of $K\setminus T$. We find
\[\sum_{T,I} \frac{\ul x^{-T^{c}}  }{\| \cZ_{T} \|^{2}} ({\rP}_{T})^{\wedge k}(e^{\star}_{I}/\sqrt{\ul\x^{I}}) \otimes e^{\star}_{I}/\sqrt{\ul\x^{I}}\]
and since $({\rP}_{T})^{\wedge k}(e^{\star}_{I})=0$ whenever the subset $I$ is not disjoint from the basis $T$, we can read the last sum as a double sum over all bases $T$ of $\cM$ and all $k$-subsets $I$ of $S$. An application of the mean projection theorem (Proposition \ref{coro:mean-proj-matroid}) concludes the proof.
\end{proof}

\subsection{Determinantal probability measures}
We will now construct determinantal probability measures on the sets $\B_{k}(\cM)$ of bases of the union matroid $\cM_{k}$. This construction is parallel to that performed in Section \ref{sec:connected}.

Let us choose an integer $k\geq 0$ and a linear subspace $H$ of $E^{*}$ such that 
\[\im d \subseteq H \ \text{ and } \ \dim H = \rank(d)+k.\]
Let us choose $k$ vectors $\theta_{1},\ldots,\theta_{k}$ of $E^{*}$ such that
\[H=\im d \oplus \Vect(\theta_{1},\ldots,\theta_{k})\]
and set
\[\vartheta=\theta_{1}\wedge \ldots \wedge \theta_{k} \in \ext^{k}E^{*}.\]

We want to use $\vartheta$ to associate to every $K\in \B_{k}(\cM)$ a non-negative weight, but we do not have any preferred basis of the vector space $Z_{K}$. 

We can nevertheless consider a quantity such as $\|  (\proj{{\sf J}_{\ul x} (Z_{K})})^{\wedge k} \vartheta \|^{2} $, which is independent of the choice of a basis of $Z_{K}$. If a basis $\cZ^{K}=(z_{1}, \ldots, z_{k})$ of $Z_{K}$ is chosen anyway, this number can be written more concretely as
\begin{equation}\label{eq:weightSym-mat}
\big\|  \big(\proj{{\sf J}_{\ul x} (Z_{K})}\big)^{\wedge k} \vartheta \big\|^{2}=\frac{\big| \big(\vartheta,\det \cZ^{K}\big)\big|^{2}}{\|\cZ^{K}\|^{2}}=\frac{\big\vert\det\big(\theta_i(z_{j})\big)_{1\le i,j\le k}\big\vert^2}{\|\cZ^{K}\|^{2}}.
\end{equation}

We can also, in a less canonical way, choose a basis $T$ of $\cM$ contained in $K$, and consider the basis $\cZ^{K}_{T}$ of $Z_{K}$. The next proposition relates these two approaches. The random subset $\X$ of $S$ is that described by Proposition \ref{prop:wisDPP}.

\begin{proposition}\label{prop:weight-matroid}
Let $K\in \B_{k}$ and $T\in \B_{0}$ such that $T\subseteq K$. Then 
\[w(T) \, \big| \big(\vartheta,\det \cZ^{K}_{T}\big)\big|^{2}\, \ul{x}^{K}=\rB(\ul{x}) \, 
\P(\X\subseteq K)\, \big\|  \big(\proj{{\sf J}_{\ul x} (Z_{K})}\big)^{\wedge k} \vartheta \big\|^{2}.\]
In particular, the left-hand side does not depend on the basis $T$ chosen within $K$, and the right-hand side is, as a function of the Euclidean weights $\ul x$, a constant times $\ul{x}^{K}$.
\end{proposition}

\begin{proof} Let us start from the right-hand side. Let us apply \eqref{eq:weightSym-mat} with the basis $\cZ^{K}_{T}$ of $Z_{K}$ to compute the middle term, and Proposition \ref{lem:dpp-mat} to compute the probability $\P(\X\subseteq K)$, with the basis $T$. A simplification of the squared volume of $\cZ^{K}_{T}$ occurs and we find that the right-hand side is equal to
\[\ul{x}^{-K^{c}} \rB(\ul{x})\frac{\big| \big(\vartheta,\det \cZ^{K}_{T}\big)\big|^{2}}{\|\cZ_{T}\|^{2}}=
\frac{\ul{x}^{-T^{c}}}{\|\cZ_{T}\|^{2}}\,  \rB(\ul{x})\, \big| \big(\vartheta,\det \cZ^{K}_{T}\big)\big|^{2}
 \, \ul{x}^{K\setminus T}.\]
 By Proposition \ref{prop:proba-det}, the first term is equal to $\P(\X=T)$, and by Proposition \ref{prop:wisDPP}, the product of this term with $\rB(\ul{x})$ is equal to $w(T)\, \ul{x}^{T}$. The result follows.
\end{proof}

To each $K\in \B_{k}$, let us associate the weight
\begin{equation}\label{eq:defweightK}
w(K)=w(T) \, \big| \big(\vartheta,\det \cZ^{K}_{T}\big)\big|^{2},
\end{equation}
where $T$ is any basis of $\cM$ contained in $K$. Let us emphasize that this number does not depend on the Euclidean weights $\ul{x}$.

Let us also define the generating polynomial 
\begin{equation}\label{eq:defLpol}
{\rL}_{k}(\vartheta,\ul\x)=\sum_{K\in \B_{k}(\cM)} w(K)\, \ul{x}^{K}=\rB(\ul{x}) \sum_{K\in \B_{k}(\cM)} \P(\X\subseteq K)\, \big\|  \big(\proj{{\sf J}_{\ul x} (Z_{K})}\big)^{\wedge k} \vartheta \big\|^{2} .
\end{equation}

The main theorem of this section is the following.

\begin{theorem}\label{thm:DPPmat} \label{thm:matroid}
Let $\Y$ be a determinantal random subset of $S$ associated to $H$ in the orthonormal basis $(e^{\star}/\sqrt{x_{e}})_{e\in S}$ of  $E^{*}$. Then $\Y$ belongs to $\B_{k}(\cM)$ almost surely and for every $K\in \B_{k}(\cM)$, 
\begin{equation}\label{eq:distribK-mat}
\P(\Y=K)=w(K)\, \ul{x}^{K} \ \big/ \ \rL_{k}(\vartheta,\ul\x)=w(T) \, \big| \big(\vartheta,\det \cZ^{K}_{T}\big)\big|^{2} \ \big/ \ \rL_{k}(\vartheta,\ul\x),
\end{equation}
where $T$ is any basis of $\cM$ contained in $K$.
\end{theorem}

The proof relies on Proposition \ref{prop:routine}, for the application of which we need to introduce a linear map with range equal to $H$. To this end, let us consider the map 
\[\omega_{\vartheta} : \C^k \to E^{*}, \ \ (t_{1},\ldots, t_{k})\mapsto t_{1} \theta_1+\ldots+ t_{k} \theta_k.\]
This map does not only depend on the tensor $\vartheta$, but on the whole family $(\theta_{1},\ldots,\theta_{k})$, so that our notation is slightly abusive. 
The space $\C^{k}$ being endowed with the usual Hermitian inner product, let us consider the orthogonal direct sum $(\ker d)^{\perp} \oplus \C^{k}$ and the linear maps
\[d\oplus \omega_{\vartheta}:(\ker d)^{\perp}\oplus \C^{k} \to E^{*} \ \text{ and } \ \Delta_{\vartheta}=(d\oplus \omega_{\vartheta})^*(d\oplus \omega_{\vartheta}).\]

Theorem \ref{thm:DPPmat} will be a consequence of the next proposition.

\begin{proposition} \label{prop:egalitesmat} One has the equalities
\begin{equation}\label{eq:==mat}
\rL_{k}(\vartheta,\ul\x)=\rB(\ul{x}) \big\| (\proj{\ker d^{*}})^{\wedge k}(\vartheta)\big\|^{2}= \det \Delta_{\vartheta}.
\end{equation}
\end{proposition}

\begin{proof} Let us start by observing Proposition \ref{prop:gmi-mat} can be written under the form 
\[\sum_{K\in \B_{k}} \P(\X\subseteq K) \frac{\big(({\sf J}_{\ul x})^{\wedge k} \det \mathscr Z^{K}\big)^{\otimes 2}}{\|\cZ^{K}\|^{2}}=\sum_{I\subseteq S, |I|=k} (\proj{\ker d^{*}})^{\wedge k}\big(e^{\star}_{I}/\sqrt{\ul\x^{I}}\big)\otimes e^{\star}_{I}/\sqrt{\ul\x^{I}}.\]
Applying the sesquilinear form $\langle \vartheta, \cdot \rangle \otimes \langle \cdot, \vartheta \rangle$ to both sides of this equality, as we did in the proof of Proposition \ref{prop:egalitesconnexes}, and multiplying by $\rB(\ul{x})$, we find
\[\rB(\ul{x})\sum_{K\in \B_{k}(\cM)} \P(\X\subseteq K)\,  \big\|  \big(\proj{{\sf J}_{\ul x} (Z_{K})}\big)^{\wedge k} \vartheta \big\|^{2} =\rB(\ul{x})\big\| (\proj{\ker d^{*}})^{\wedge k}(\vartheta)\big\|^{2},\]
recognise ${\rL}_{k}(\vartheta,\ul\x)$ in the left-hand side, and the first equality is proved.

To prove the second equality, let us apply the Schur complement formula under the form given by Lemma \ref{lem:Schurcomp}, with $E_{0}=(\ker d)^{\perp}\oplus \C^{k}$, $E_{1}=E^{*}$, $G=\C^{k}$ and $a=d\oplus \omega_{\vartheta}$. We find
\[\det \Delta_{\vartheta} =\det(d^*d)^{\im d^{*}}_{\im d^{*}}\  \det\big(\omega_{\vartheta}^{*}\proj{\ker d^{*}} \omega_{\vartheta}\big).
\]
The first factor is equal to $ \rB(\ul\x)$, by definition (see \eqref{eq:mtt-matroid}). Let us compute the second. For this, let us observe that the adjoint of $\omega_{\vartheta}$ is given, for all $\alpha\in E^{*}$, by $\omega_{\vartheta}^{*}(\alpha)= (\lanx\theta_{1},\alpha\ranx,\ldots,\lanx\theta_{k},\alpha\ranx)$, so that the matrix in the canonical basis of $\C^{k}$ of $\omega_{\vartheta}^{*}\proj{\ker d^{*}} \omega_{\vartheta}$ is
$ \big(\lanx \theta_{i},\proj{\ker d^{*}}\theta_{j} \ranx\big)_{1\leq i,j \leq k}$
and the equality
\[\det\big(\omega_{\vartheta}^{*}\proj{\ker d^{*}} \omega_{\vartheta}\big)= \big\| (\proj{\ker d^{*}})^{\wedge k} \vartheta \big\|^{2}\]
holds, concluding the proof.
\end{proof}

We can now conclude the proof of the main theorem. 

\begin{proof}[First proof of Theorem \ref{thm:DPPmat}] In view of the equality of the first and third terms of \eqref{eq:==mat}, an application of the first part of Proposition \ref{prop:routine} to the operator 
 $d\oplus \omega_{\vartheta}$, whose range is the subspace~$H$ of $E^{*}$, implies that the random subset $\Y$ is determinantal with distribution given by \eqref{eq:distribK-mat}.
\end{proof}

Just as Theorems \ref{thm:DPP1} and \ref{thm:DPP2}, we can prove Theorem \ref{thm:DPPmat} via a second method, which we described as going from local to global in the introduction. We will only give a sketch of this proof and refer the interested reader to the second proof of Theorem \ref{thm:DPP1}.

\begin{proof}[Second proof of Theorem \ref{thm:DPPmat}] We first make sure that $\Y$ almost surely belongs to $\B_{k}$. Once this is done, we choose $K\in \B_{k}$ and compute $\P(\Y=K)$ using \eqref{eq:defloiX}.

 For this, we construct a basis of $H$ by choosing an arbitrary basis of $\im d$ and completing it with $\theta_{1},\ldots,\theta_{k}$. We also choose a basis $T$ of $\cM$ contained in $K$. Then the matrix of which the square of the determinant computes $\P(\Y=K)$ has a $2\times 2$ block structure corresponding to the decompositions $H=\im d \oplus {\rm Vect}(\theta_{1},\ldots,\theta_{k})$ and $E^{*}_{K}=E^{*}_{T}\oplus E^{*}_{K\setminus T}$. 

Performing elementary row or column operations on this matrix, in a way which amounts to replacing the basis vectors $(e^{\star}_{i} : i\in K\setminus T)$ by the basis ${\sf J}_{\ul x}(\cZ^{K}_{T})$, makes it block triangular, and allows us to write its determinant as the product of two terms.

The first term corresponds to $\im d$ and $E^{*}_{T}$, and is, up to an explicit constant, equal to the weight $w(T) \, \ul{x}^{T}$. The second term corresponds to ${\rm Vect}(\theta_{1},\ldots,\theta_{k})$ and $\cZ^{K}_{T}$ and produces the factor~$|(\vartheta,\det \cZ^{K}_{T})|^{2}\, \ul{x}^{K\setminus T}$.
\end{proof}


\section{Further examples}\label{sec:examples}

We have seen in Section \ref{sec:dss} how, starting from the uniform spanning tree, we could generate new interesting random spanning subgraphs. In the light of Section \ref{sec:matroid-setting}, we see that this case corresponded to the circular matroid being extended by $k$ points. With the general result Theorem \ref{thm:matroid} at hand, we see several avenues of generalisation:
\begin{itemize}
\item apply the construction to another matroid on a graph, for instance the bicircular matroid; see Section~\ref{sec:bicircular-case};
\item apply the construction iteratively, by adding or removing points to a given matroid, for instance starting with a seed given by the circular matroid or the bicircular matroid; see Section~\ref{sec:iterated-case};
\item apply the theorem to a `circular' or `bicircular' matroid on higher rank vector bundles on graphs or complexes; see Section~\ref{sec:higher-rank-bundles};
\item apply the theorem to a `circular' or `bicircular' matroid on simplicial complexes; see Section~\ref{sec:higher-dim-cells};
\item combine the last two generalizations to higher dimension, and higher rank.
\end{itemize}

\subsection{Bicircular case}\label{sec:bicircular-case}

The bicircular matroid of a graph $\rG$ is the matroid on its set of edges, the set of bases of which is the set $\cU(\rG)$ of cycle-rooted spanning forests. Its collection of circuits is the set of connected subgraphs with one more edge than vertices, that is the subgraphs with Betti numbers $(b_0,b_1)=(1,2)$. 

\subsubsection{Linear representation of the bicircular matroid}

For more details about the content of this section, see~\cite{Kenyon, Kassel-ESAIM, KL3, KL6}.

We consider a connection, that is the data, for each half-edge $(\ul{e},e)$ or $(e,\ol{e})$, of a unitary complex number $h_{\ul{e},e}$ or $h_{\ol{e},e}$. We set $h_{e,\ul{e}}=h_{\ul{e},e}^{-1}$ and $h_{e,\ol{e}}=h_{\ol{e},e}^{-1}$. Given a connection, we define a twisted coboundary map $d_{h}:\Omega^0(\rG)\to \Omega^1(\rG)$ by 
\[d_hf(e)= h_{e,\ol{e}} f(\ol{e}) - h_{e,\ul{e}}f(\ul{e}),\]
and let $d_h^*$ be the adjoint linear map. The twisted Laplacian is defined by 
$\Delta_h=d_h^* d_h$. Then for all $f\in \Omega^0(\rG)$, and all $v\in \rV$, we have $\Delta_h f(v)=\sum_{e\in \rE: \ol{e}=v} x_e (f(v)-h_{v,\ul{e}} f(\ul{e}))$ where $h_{v,\ul{e}}=h_{v,e}h_{e,\ul{e}}$.
We still denote by $\Delta_h$ the matrix of this operator in the canonical basis of $\Omega^0$ indexed by vertices. 

By a theorem of Forman \cite{Forman}, rediscovered and extended to quaternions by Kenyon \cite{Kenyon} (see \cite{KL2, KL2C} for a short proof and references therein), for each $r$, the $r\times r$ principal minors of $\Delta_h$ are a weighted sum over bases of the $r$-truncation of the bicircular matroid. We denote by $\B_{-k}(\rG)$ the set of bases of this matroid, with $r+k=|\rV|$. In particular, $\B_{0}(\rG)=\cU(\rG)$.  

\begin{theorem}[Forman, Kenyon]\label{thm:forman-kenyon-minors}
Let $W\subseteq \rV$ be a subset of size $k$. Then
\begin{equation}
\det (\Delta_h)_{\widehat{W}}^{\widehat{W}}=\sum_{F\in \B_{-k}(\rG)} \prod_{c\, \text{\rm cycle}} \vert 1-\hol_h(c)\vert^2 \; \ul{x}^F \1_{\{F\, \text{\rm separates points of}\; W\sqcup \{\text{\rm cycles}\}\}}\,
\end{equation}
where for each $F\in \B_{-k}(\rG)$, the product is over the set of cycles, each with a fixed chosen orientation.
\end{theorem}

\subsubsection{Determinantal cycle-rooted spanning forests}

This twisted matrix-tree formula (Theorem~\ref{thm:forman-kenyon-minors} in the case $W=\varnothing$) implies the existence of a determinantal probability measure on~$\cU(\rG)$ according to Proposition~\ref{prop:routine}, a fact first proved by Kenyon \cite{Kenyon} who rediscovered Forman's result, and extended it to the quaternion case. See Figure~\ref{fig:crsf} for a sample. This is the analogue, in this twisted setting, of the Burton--Pemantle theorem (Proposition \ref{prop:lawUST}) and its first proof. 

\begin{theorem}[Kenyon]
The probability distribution on $\cU(\rG)$ which assigns any $F$ a weight proportional to 
\begin{equation}\label{eq:wF-crsf}
w(F)\,\ul{x}^F=\ul\x^F \prod_{c\,\text{\rm cycle}} \vert 1-\hol_h(c)\vert^2\,,
\end{equation}
is determinantal associated with the subspace $\im d_h$ of $\Omega^1(\rG)$ in the orthonormal basis $(\frac{e^\star}{\sqrt{x_e}}:e\in \rE^+)$.
\end{theorem}

\begin{figure}[!ht]
\centering
\includegraphics[width=6cm]{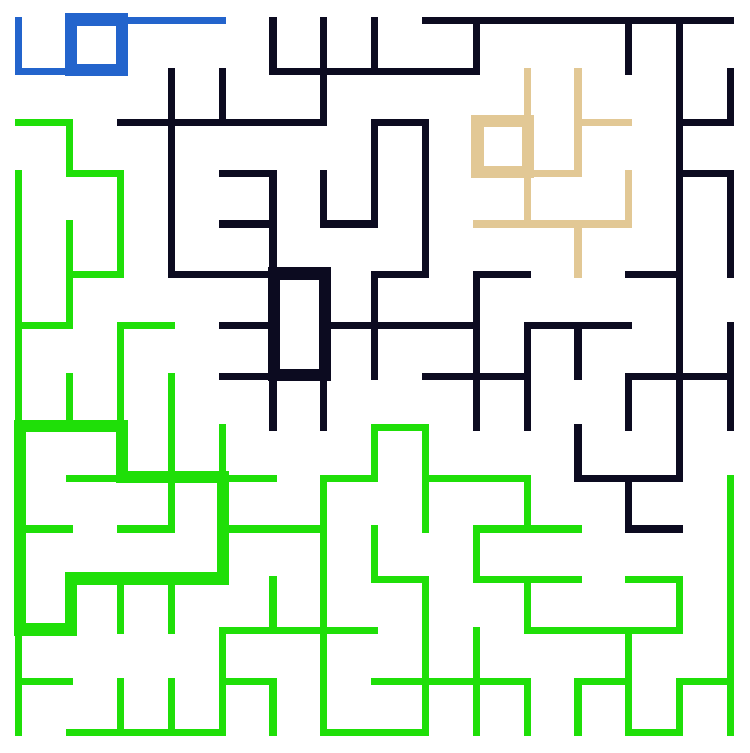}
\caption{\small A determinantal random cycle-rooted spanning forest, that is a random basis of the bicircular matroid of a $15\times 15$ grid. Its cycles are represented by thickened edges. }\label{fig:crsf}
\end{figure}

We see from \eqref{eq:wF-crsf} that the support of the measure is a certain subset of $\cU(\rG)$ corresponding to those $F$ all of whose cycles have non-trivial holonomy. This defines a submatroid of the bicircular matroid.

One sees from \eqref{eq:wF-crsf} that the determinantal distribution on $\cU(\rG)$ only depends on the holonomy representation $\rho$ of the fundamental group defined by the connection $h$ (see e.g.~\cite{CK} for a review of these notions). 

\subsubsection{Symanzik polynomials in the bicircular case}\label{sec:bicircular-Symanzik-polynomials}

Let us define a family of multivariate polynomials in the Euclidean weights $x_e$ associated to the linear representation of the bicircular matroid we just introduced, and more precisely to its holonomy representation~$\rho$.

First, we define the weighted basis generating polynomial for the bicircular matroid `measured' by its representing map $d_h^*$. It reads
\begin{equation}
\rB^{\rho}(\ul{x})=\sum_{F\in \cU(\rG)} \prod_{c\, \text{cycle}} \vert 1-\hol_h(c)\vert^2 \ul{x}^F\,.
\end{equation}
This is the special case of the polynomial $\rB(\ul{x})$ defined in \eqref{eq:mtt-matroid} in Section \ref{sec:matroid-setting}.
From Theorem \ref{thm:forman-kenyon-minors}, we have 
\begin{equation}\label{eq:Brho-det}
\rB^\rho(\ul{x})=\det\Delta_h\,.
\end{equation}

Next, we define a family of polynomials indexed by positive integers $k$, corresponding to weighted basis generating polynomials of truncations of the bicircular matroid. The bases of the $r$-truncation of the bicircular matroid are multi-type spanning forests, that is spanning subgraphs such that each connected component is either a tree or a unicycle (a graph with first Betti number equal to $1$), and such that there are exactly $k=|\rV|-r$ trees (the number of unicycles is not fixed)\footnote{Note that such multi-type spanning forests have been considered by Kenyon \cite{Kenyon,Kenyon-massive} and more recently in a different context motivated by applications by \cite{Bardenet-Fanuel,Jaquard})}. 
Recall that we denote $\B_{-k}(\rG)$ this set of bases, and further note that $\cF_{k-1}(\rG)\subset\B_{-k}(\rG)$.

For all positiver integer $k$ and all $k$-tuple ${\sf q}=(q_1,\ldots, q_k)$  of elements of $\Omega^0(\rG)$, we set 
\begin{equation}\label{eq:Arho}
\rA^{\rho}_k({\sf q},\ul{x})=\sum_{F\in \B_{-k}(\rG)} \vert\det\left( q_i (\rV(T_j)) \right)_{1\le i,j\le k}\vert^2 \;  \prod_{c \, \text{cycle}} \vert 1-\hol_h(c)\vert^2 \; \ul{x}^F\,,
\end{equation}
where for each $F\in \B_{-k}(\rG)$, we write $\{T_1,\ldots, T_k\}$ for its connected components which are trees, and the product is over the simple cycles of the unicyclic connected components. These polynomials are the analogues of the Symanzik polynomials ${\rA}_{k}(\calc,\ul\x)$ defined in \eqref{eq:defApol} of Section~\ref{sec:dss}.

\subsubsection{Twisted Green function minors}

Let us assume that there is at least one cycle $c$ in $\rG$ such that $\hol_h(c)\ne 1$ (this is equivalent to $\rho\ne 1$). Otherwise, we are back in the circular case of Section~\ref{sec:dss}. This implies that $\Delta_h$ is invertible, and we let~$G_h$ be the inverse of $\Delta_h$. The following proposition is the analogue, in the bicircular matroid case, of Proposition~\ref{prop:A-T-G} which holds for the circular matroid case.

\begin{proposition}\label{prop:bicircular-green}
For all ${\sf q}=(q_1,\ldots, q_k)$ $k$-tuple of elements of $\Omega^0(\rG)$, we have
\begin{equation}\label{eq:height-pairing}
\frac{\rA^\rho_k({\sf q},\ul{x})}{\rB^\rho(\ul{x})} = \det \left(\langle q_i, G_h q_j \rangle \right)_{1\le i,j\le k} \,.
\end{equation}
\end{proposition}

Note that the right-hand side of \eqref{eq:height-pairing} can be written in the language of the exterior algebra as $\langle q_1\wedge \ldots \wedge q_k ,\ext^kG_h (q_1\wedge\ldots\wedge q_k) \rangle$.

\begin{proof}
We first note that each side of the equation we want to prove is a quadratic form in $q_1\wedge \ldots \wedge q_k$. Hence it is enough to prove it for an orthonormal basis of $\ext^k \Omega^0(\rG)$, and then conclude by bilinearity. We will thus prove the formula for each collection of disjoint $k$ vertices~$v_i$, letting $q_i=\1_{v_i}$.

In that case, the right-hand side is the principal minor of $G_h$ indexed by $W=\{v_1,\ldots, v_k\}$, which by Jacobi's complementary minors formula is the quotient of $\det(\Delta_h)_{\widehat{W}}^{\widehat{W}}$ by $\det \Delta_h$.

Now, we turn to the numerator of the left-hand side. By \eqref{eq:Arho}, this is the sum over $F$ of a weight which we now compute. We observe that 
\[\det\left(q_i(\rV(T_j))\right)_{1\le i,j\le k}\] 
is nonzero if and only if there exists a bijection $\sigma$ of $\{1,\ldots, k\}$ such that $v_i$ belongs to $T_{\sigma(i)}$, in which case, the determinant is the signature of this permutation matrix, equal to $\pm 1$. Hence in that case, $\rA^\rho_k({\sf q},\ul{x})=\det(\Delta_h)_{\widehat{W}}^{\widehat{W}}$. The result now follows by invoking \eqref{eq:Brho-det}.
\end{proof}

We note that since $d_h^*$ is onto, we may write for each $i\in \{1,\ldots, k\}$, that $q_i=d_h^* \phi_i$, for $\phi_i \in \Omega^1(\rG)$. Then, letting $\varphi=\phi_1\wedge\ldots\wedge\phi_k$, Equation~\eqref{eq:height-pairing} reads 
\begin{equation}\label{eq:abphi}
\lVert \ext^k \proj{\im d_h} (\varphi) \rVert^2 = \frac{\rA^\rho_k({\sf q},\ul{x})}{\rB_\rho(\ul{x})}\,.
\end{equation}

\subsubsection{Determinantal multi-type spanning forests}

Given $k$ functions $q_1,\ldots, q_k\in \Omega^0(\rG)$, let 
\begin{equation}
H_{{\sf q}}=\{ d_h f: f\in \Omega^0(\rG), \forall 1\le i\le k, \langle f, q_i\rangle=0\}\,.
\end{equation}
For any choice of $\phi_1,\ldots,\phi_k\in \Omega^1(\rG)$ such that $q_i=d_h^* \phi_i$ for all $i\in\{1,\ldots, k\}$, we may rewrite this space as
\begin{equation}
H_{{\sf q}}=\im d_h \cap \Vect(\phi_1, \ldots, \phi_k)^\perp\,.
\end{equation}

The following theorem is the analogue, in the bicircular matroid case, of Theorem~\ref{thm:sym-forests} which holds for the circular matroid case. See Figure \ref{fig:bicircular-truncation} for an example.

\begin{figure}[!ht]
\centering
\includegraphics[width=6cm]{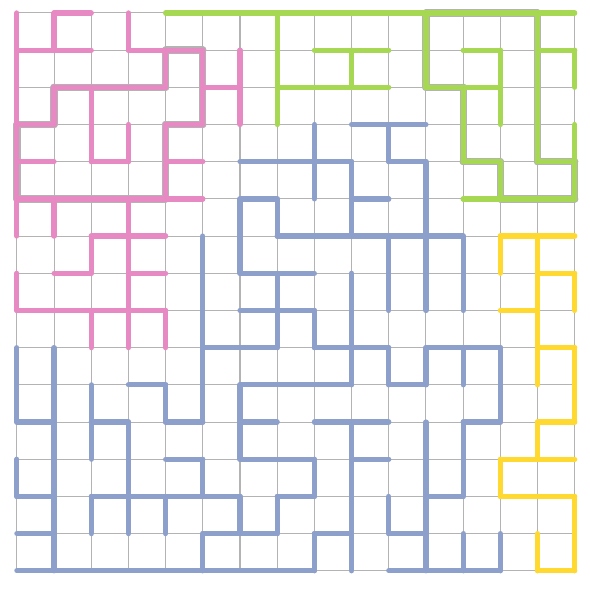}
\caption{\small A determinantal random element of $\B_{-2}(\rG)$.}\label{fig:bicircular-truncation}
\end{figure}

\begin{theorem}
The probability distribution on $\B_{-k}(\rG)$ associated to the polynomial $\rA^{\rho}_k({\sf q},\ul{x})$, seen as the generating function of a positive measure, is determinantal associated to the orthogonal projection on $H_{{\sf q}}$ in the canonical orthonormal basis of $\Omega^1(\rG)$ indexed by $\rE^+$. 
\end{theorem}

\begin{proof}
We follow the same strategy as in the first proof of Theorem \ref{thm:DPP2}. Hence we want to apply the second part of Proposition \ref{prop:routine} and for that we need an expression for the partition function of our probability measure of the form $\ul{x}^{\rE}\det bb^*$ where $\ker b=H_{{\mathsf q}}$. 

For that, we first need to prove an analogue of the second equality in Proposition \ref{prop:egalitesacycliques}. We will do so using Lemma~\ref{lem:fundamental-global} in place of \eqref{eq:sym1-omid}.

Let us start by defining the map $b$. We follow word for word the construction in the circular matroid case. The only difference is that we do not use $Z_1(\rG,\K)$. Instead we replace that space with $Z=({\sf J}^{1}_{\ul 1})^{-1}\ker d_h^\dagger\subset \Omega_1(\rG)$. We let $c_i$ be the preimages by $\Jx$ of $\phi_i$. 

The key point is that we will observe how changing the Euclidean weights from $1$ to $\ul{x}$ changes the volume of an \emph{arbitrary} basis $\cZ$ of $Z$ (we have no integral structure here unlike in the circular matroid case). This is where Lemma \ref{lem:fundamental-global} comes into play. We find that 
\[
\lVert \det \cZ \lVert^2=\sum_{T \in \cU(\rG)}\det(\cZ/\cZ_T)^2 \, \frac{\ul{x}^T}{\ul{x}^\rE}\,.
\]
Let us call $\cV(\ul{x)}$ this function of $\ul{x}$. By the same argument as in the proof of Proposition \ref{prop:egalitesacycliques},
we find that $\det(1^Q bb^*1_Q)$ is equal to $\cV(\ul{x})/\cV(\ul{1})$.

But now, combining Proposition \ref{prop:proba-det} and \eqref{eq:wF-crsf}, we find that 
\[T\mapsto \frac{w(T)}{(\det\cZ/\cZ_T)^2}\] 
is a constant function of $T\in \cU(\rG)$. 
Hence, the above ratio of volume $V(\ul{x})/V(\ul{1})$ is proportional to $\rB_\rho(\ul{x})/\ul{x}^{\rE}$ and evaluating it at $\ul{x}=\ul{1}$ where it is equal to $1$, we find that the volume-ratio is thus equal to 
\[
\frac{\cV(\ul{x})}{\cV(\ul{1})}=\frac{\rB_\rho(\ul{x})}{\rB_\rho(\ul{1})\ul{x}^{\rE}}\,.
\]

Hence, combining \eqref{eq:abphi}, and the same Schur lemma (Lemma \ref{lem:Schurcomp}) argument as in the proof of Proposition \ref{prop:egalitesacycliques} we are following step by step, we obtain
\begin{equation}
\rA^\rho_k({\mathsf q}, \ul{x})=\rB_\rho(\ul{x}) \big\| \ext^k\proj{\im d_h} (\varphi)\big\|^2=  \rB_\rho(\ul{1}) \ul{x}^\rE \det(bb^*)\,.
\end{equation}
The result now follows by an application of the second determinantality criterion (second assertion of Proposition \ref{prop:routine}).
\end{proof}

\subsubsection{An example on the `rank $1$ extension' of the bicircular matroid}

Let us consider a unitary rank~$1$ connection $h$ on $\rG$. As we have just seen, the subspace $\im d_h$ of~$\Omega^1(\rG)$ defines a determinantal probability measure on $\cU(\rG)$. And by considering subspaces $H\subseteq \im d_h$, we defined determinantal multi-type spanning forests. Let us now take a look at what happens if we take a subspace $H$ containing $\im d_h$, in the simplest case where the codimension of $\im d_h$ in $H$ is $1$.

Hence, let us consider a $1$-form $\theta$ outside $\im d_h$ and consider the subspace 
\[H_\theta=\im d_h \oplus \C \theta\subseteq \Omega^1(\rG)\,.\] 
By Theorem \ref{thm:matroid}, the corresponding determinantal measure is supported on the set of spanning subgraphs $K$ with connected components $K_{1},\ldots,K_{m}$, all of which have first Betti number equal to~$1$, but $K_{1}$, which has first Betti number equal to $2$.

Given such a $K$, we pick a cycle-rooted spanning forest $F\subseteq K$. The unique edge in $K\setminus F$ belongs to $K_{1}$, indeed to the $2$-core of $K_{1}$. This edge determines a fundamental $1$-chain $z^K_F$, as explained in Section~\ref{sec:fund-bases-mat}. We observe that it only depends on~$K_1$ and~$e_1$ and leave it to the interested reader to compute $z^K_F$ explicitly. Then, by Theorem \ref{thm:matroid}, and using \eqref{eq:wF-crsf}, the weight of $K$ is 
\begin{equation}
w(K)= \vert \theta(z^K_F)\vert ^2\; \prod_{c} |1-h_{c}|^{2} \,,
\end{equation}
where the product is over the cycles of $F$.

\subsection{Iteration of the construction}\label{sec:iterated-case}

Since our construction yields a new measured matroid from a given one, we can apply this construction iteratively, in particular intertwining duality at each step to shuffle the matroid in two `opposite directions'.

\subsubsection{Perturbation of spanning trees: subgraphs with fixed Euler characteristic}\label{sec:betti-fixed}

For each pair of non-negative integers $(k,\ell)\in \{0,\ldots, \vert \rV\vert-1\}\times \{0, \ldots, b_1(\rG)\}$, there is a matroid $\cM_{k,\ell}(\rG)$ on the set of edges of $\rG$ whose set of bases $\B_{k,\ell}(\rG)$ is the collection of subgraphs $B$ of $\rG$ satisfying
\begin{itemize}
\item $\chi(B)=k-\ell+1$
\item $\max(0,\ell-k)\le b_1(B)\le \ell$
\end{itemize}
where $\chi(B)=b_0(B)-b_1(B)=\vert \rV(B)\vert-\vert \rE(B)\vert$ is the Euler characteristic of $B$. We recover the previously considered families of subgraphs: $\B_{0,0}(\rG)=\cT(\rG)$, $\B_{0,\ell}(\rG)=\cC_\ell(\rG)$ and $\B_{k,0}(\rG)=\cF_{k}(\rG)$. The elements of $\B_{k,\ell}(\rG)$ are the subgraphs obtained by taking any spanning tree, adding $\ell$ edges, and then removing $k$ edges. See Figure~\ref{fig:bkl}.

\begin{figure}[!ht]
\centering
\includegraphics[width=6cm]{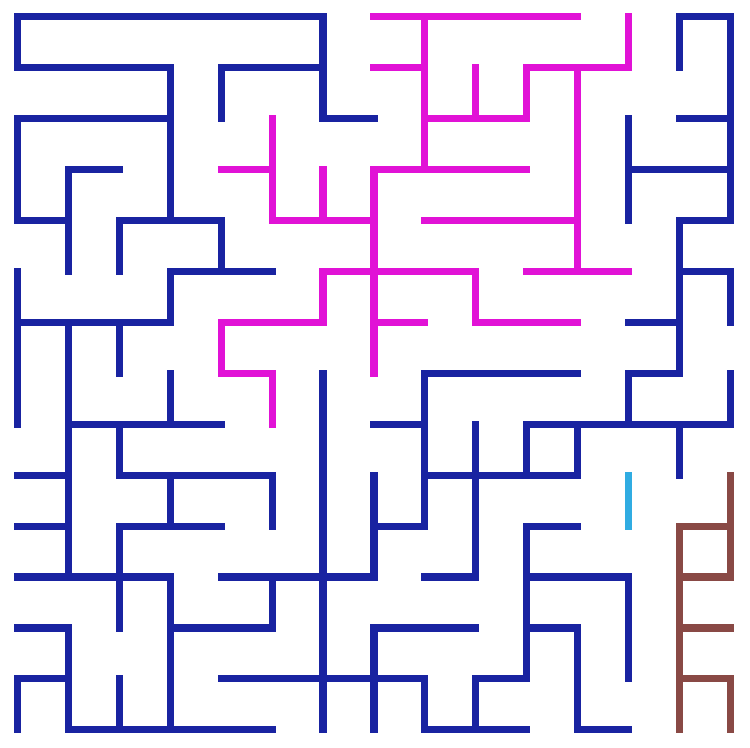}
\caption{\small A determinantal random element of $\B_{4,4}(\rG)$ on a $15\times 15$ grid. Note that this subgraph also belongs to $\B_{3,3}(\rG)$, and indeed to $\B_{k,k}(\rG)$ for any $k\geq 3$.}\label{fig:bkl}
\end{figure}

The existence of these matroids, simply obtained by taking unions with the uniform matroid, or duals, starting from the circular matroid, does not contradict the result of \cite{Simoes-Pereira} mentioned in the introduction: the circular and bicircular matroids of a graph are the only two matroids such that their sets of circuits are homeomorphism classes of connected graphs. Indeed the circuits of $\cM_{k,0}$ consist in simple cycles and spanning forests with $k$ components; because of the \emph{spanning} assumption this is not the class of subgraphs \emph{homeomorphic} to a fixed subclass of connected graphs. Similarly, the circuits of $\cM_{0,\ell}(\rG)$ are the minimal subgraphs with $\ell+1$ independent cycles; because these subgraphs are not necessarily \emph{connected}, this is not the class of subgraphs homeomorphic to a fixed subclass of \emph{connected} graphs.

For all $(k,\ell)$ in the above range, one can define natural determinantal probability measures whose supports are included in $\B_{k,\ell}(\rG)$ by taking a subspace of the form $\im d \cap ({\Jx}(\Phi))^{\perp} \oplus \Theta$, with $\Phi\subset \Omega_1(\rG)$ a $k$-dimensional subspace such that $\Phi\cap Z_1(\rG,\K)=\{0\}$ and $\Theta\subset \Omega^1(\rG)$ an~$\ell$-dimensional subspace such that $\Theta\cap\im d \cap ({\Jx}(\Phi))^{\perp}=\{0\}$. 

Repeated uses of Theorem \ref{thm:matroid} can in principle allow us to describe explicitly these measures via geometric-topological weights on subgraphs.

\subsubsection{Perturbation of cycle-rooted spanning forests}\label{sec:crsf}

As in the circular case above, we may define, for all $(k,\ell)\in \{0, \ldots, \vert\rV\vert\}\times \{0,\ldots, b_1(\rG)-1\}$, variants $\cM^{\text{bicirc}}_{k,l}(\rG)$, where $\cM^{\text{bicirc}}_{0,0}(\rG)$ is the bicircular matroid of~$\rG$. The collection of bases of~$\cM^{\text{bicirc}}_{k,l}(\rG)$ is obtained as the collection of subgraphs of $\rG$ built by adding $k$ edges to any element of $\cU(\rG)$, and then removing $\ell$ edges. 

For all $(k,\ell)$ in the above range, we can define natural determinantal probability measures whose support is included in $\B_{k,\ell}^{\text{bicirc}}(\rG)$, by considering a subspace of the form $\im d_h \cap ({\Jx}(\Phi))^{\perp} \oplus \Theta$ with assumptions similar to those in the circular case. The description of weights can in principle be obtained using Theorem~\ref{thm:matroid} but we have not worked this out (see Section \ref{sec:bicircular-case} however for the case $k=0$, $\ell=1$).

As proved by Kenyon in \cite[Theorem 3]{Kenyon}, when we consider $\varepsilon \theta$ in place of $\theta$, then, letting~$\varepsilon$ tend to $0$, we have a family of determinantal measures on $\cU(\rG)$ converging to a determinantal measure on~$\cC_1(\rG)$, which is precisely the measure described in Theorem \ref{thm:DPP1} for the line~$\Theta=\K\theta$. This convergence result is generalised to higher rank in \cite{KL6}; see Section \ref{sec:higher-rank-bundles} for a preview.

\subsection{Quantum spanning forests}\label{sec:higher-rank-bundles}

In \cite[Section 1.5]{KL3} and \cite{KL6}, we consider higher rank vector bundles on graphs, following our work \cite{KL1}. We consider a subspace $\im d_h$ of $\Omega^1(\rG, \K^N)$ where $d_h$ is a $\rU(N,\K)$-twisted discrete covariant derivative. Here, $\rU(N,\K)$ is the unitary group of $\K^{N}$. The quantum spanning forest is the determinantal linear process \cite[Definition~3.1]{KL3} associated with the subspace $\im d_h$ and the natural splitting of $\Omega^1(\rG,\K^N)$ as sums of blocks $\K^N$, where the sum is over~$\dE$. It is a certain random subspace $\Q$ of the form $\Q=\oplus_{e\in \dE}\Q_e$ which is $h$-acyclic in the sense that $\Q\cap \ker d_h^*=\{0\}$ and which is maximal for these properties.

Viewing $\Omega^1(\rG)$ as a linear subspace of $\Omega^1(\rG, \K^N)$ by picking a line over each edge, we then consider the compression on that subspace of the orthogonal projection onto $\im d_h$. It is an element of $\End(\Omega^1(\rG))$, which defines a determinantal random subgraph, which we call a marginal of the quantum spanning forest. See Figure~\ref{fig:qsf}. Taking a full orthogonal basis of each block~$\K^N$ over each edge, one then obtains a collection of $N$ correlated marginal subgraphs, which are the marginals of the quantum spanning forest.

\begin{figure}[!ht]
\centering
\includegraphics[width=6cm]{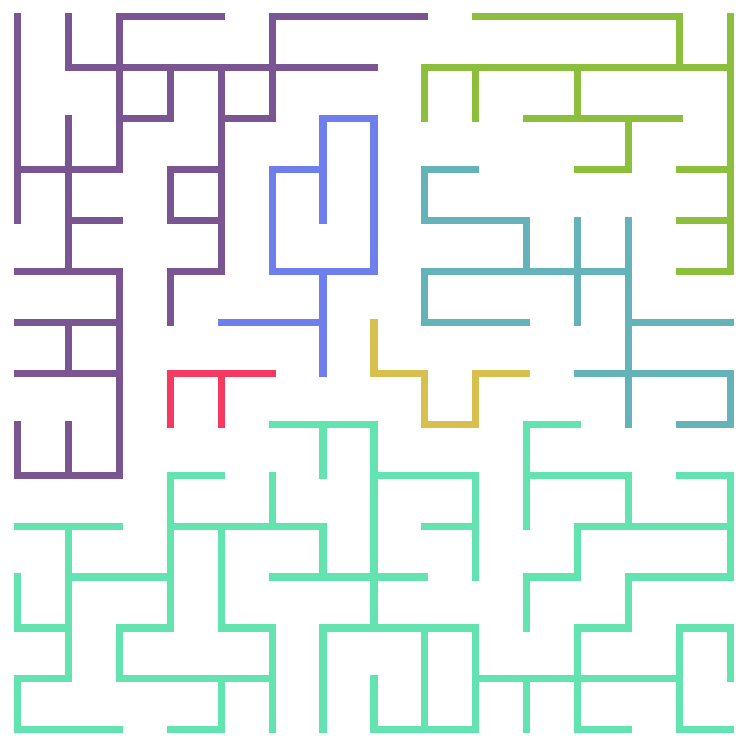}\hspace{1cm}\includegraphics[width=6cm]{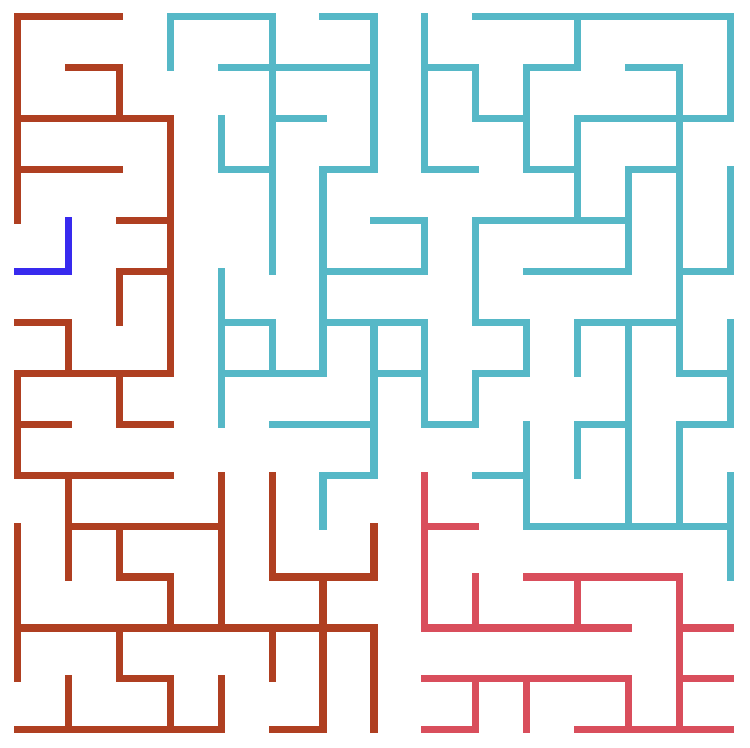}
\caption{\small The marginals of a rank-$2$ quantum spanning forest on a $15\times 15$ grid.}\label{fig:qsf}
\end{figure}

By \cite[Proposition 6.13]{KL3}, in the case where holonomies of loops are in $\SU(2)$, the law of the total occupation number of the marginals of the quantum spanning forest (like those in Figure~\ref{fig:qsf}) is equal to the occupation number of the union of two independent samples of the associated Q-determinantal measure on~$\cU(\rG)$ (like that in Figure \ref{fig:crsf}).

In the case where the connection $h$ tends to a trivial one, we obtain a determinantal probability measure on $N$-tuples of connected subgraphs whose total first Betti number is $N$. The partition function of this probability measure is quite analogous to the one in the Symanzik case considered above in the case $k=N$. We wonder whether one can obtain the above measure as the limit of a richer probability distribution, or even as a certain marginal or conditioning of the quantum spanning forest.

An easy case is when the connection is `diagonal' (see \cite[Section 1]{KL6}) in which case we have a collection of independent rank $1$ random subgraphs. In that case the probability measure on $N$-tuples of connected subgraphs is just the collection of independent elements of $\cC_1(\rG)$ chosen according to the corresponding determinantal probability measure.

\subsection{Higher dimensional random simplicial complexes}\label{sec:higher-dim-cells}

\subsubsection{Three-dimensional case}

Consider a $3$-dimensional cube, and a simplicial complex discretization $X=(X_0,X_1,X_2,X_3)$ of it. Let $\rG=X_1$ be is $1$-skeleton and let $v_0$ be a vertex lying in the bulk of the cube. Let $W$ be the set of boundary vertices, and set $q_w=\1_{v_0}-\1_w$ for each $w\in W$. By Theorem \ref{thm:DPP2}, the random element of $\cF_1(\rG)$ associated with $H=\im d\oplus (\bigoplus_{w\in W}\R q_w)$ is the uniform measure on ($\vert W\vert+1$)-component spanning forests where each connected component contains exactly one vertex $w$ of $W$ and not $v_0$, or the vertex $v_0$ and no vertex of $W$. The Poincaré dual of the connected component containing $v_0$ is a two-dimensional sphere.

\subsubsection{General case}

Instead of the circular or bicircular matroids and their variants on graphs, one can consider matroids on the cells of higher dimensional simplicial complexes such as the ones defined in \cite{Lyons-Betti} (see also \cite{Kalai, CCK-complexes, Duval-Klivans-Martin}), which could be called `circular' or `co-circular', and those mentioned in \cite[Section 1.5]{KL3} which are associated to twisted coboundary maps, and could be called `bicircular'. 
The corresponding partition functions would generalize the Kirchhoff and Symanzik polynomials, and might be related to those of \cite{Piquerez}.

For example, in the circular case in dimension $2$, one may consider a two-dimensional simplicial complex $X=(X_0,X_1,X_2)$ with chain complex map $d$, and in $\Omega^2(X)$ the space $H=\im d \oplus \R \omega$ by adding a $2$-form $\omega$ to the space of exact $2$-forms. The determinantal subset of $X_2$ associated to $H$ induces a cycle in $X_2$, which is a random two dimensional complex, and can be seen as a random surface (whose law depends on $\omega$). We can ask about the genus of this random surface model.

In the case where the two form chosen is $\omega=d^*(\1_{v_1^*}-\1_{v^*_0})$ for two $3$-faces $v_1^*$ and $v_0^*$, we simply recover the random sphere separating the two components of a uniform two-component spanning forest rooted at $v_1^*1$ and $v_0^*$ in the Poincaré dual $X^*=(X_0^*\simeq X_3,X_1^*\simeq X_1,X_2^*\simeq X_0)$ of $X$.


\section{Concluding remarks}\label{sec:conclusion}

In this concluding section, we record two observations that we came upon in  the course of writing this paper, and make brief remarks on infinite volume limits and putative scaling limits.

\subsection{Counting complexity, determinantal measures, and uniform sampling}\label{sec:complexity}

Finding the number of elements of~$\cU(\rG)$ is known to be~$\#P$-hard~\cite[Section~3]{Gimenez-Noy}\footnote{The authors of that paper use a short two-step reduction to the counting problem of perfect matchings of a graph, which is known to be $\#\mathsf{P}$-hard by a celebrated work of Valiant \cite{Valiant}, where this computational complexity class was in fact introduced.}, and thus there can be no polynomial time computable formula for it.\footnote{Bounds on the cardinality of $\cU(\rG)$ in terms of that of $\cT(\rG)$ were obtained in \cite{Gimenez-deMier-Noy}.} In particular, the uniform measure on~$\cU(\rG)$ cannot be determinantal, by \cite[Corollary 5.5]{Lyons-DPP}. However, the uniform measure on $\cU(\rG)$ may be sampled exactly in polynomial time, see \cite[Theorem~1]{Kassel-Kenyon}, \cite[Section~2.4]{Kassel-ESAIM}, and \cite{Guo-Jerrum}.\footnote{Similar variations on Wilson's algorithm \cite{Wilson} were proposed in \cite{Bouttier-Bowick-Guitter-Jeng, Gorodezky-Pak}. A general theory of partial rejection sampling was developped recently in \cite{Jerrum} which encompasses these as special cases.} This yields a fully-polynomial approximation scheme (FPRAS) for enumerating $\cU(\rG)$, as shown in~\cite{Guo-Jerrum}. 

Further note, that the bicircular is not unimodular (otherwise there would be a determinantal expression for its number of bases \cite{Maurer}), hence it is not regular \cite[Theorem 3.1.1]{White}, and hence, by a theorem of Tutte, it is not representable both on $\mathbb{F}_2$ and $\mathbb{F}_3$. However as any transversal matroid, bicircular matroids are representable over any infinite field. The question of finding over which finite fields they are representable has been studied partially by Zaslavsky.

Enumerating elements in $\cF_k(\rG)$ can be done in polynomial time, using a combination of determinants (see \cite{Liu-Chow}, further simplified by \cite{Myrvold} and \cite{Kassel-Wilson-sand}). Similarly, according to \cite{Chow}, the set $\cC_k(\rG)$ can be enumerated in polynomial time. However, the enumerations of~$\cup_{k\geq 0} \cF_{k}(\rG)$
and $\cup_{k\geq 0} \cC_{k}(\rG)$ are known to be impossible in polynomial time, as they are $\# P$-hard evaluations of the Tutte polynomial of $\rG$, see~\cite{Welsh,Provan-Ball}. In particular, the uniform measure on these sets is not determinantal in general. It is nevertheless conjectured that they satisfy a form of negative dependence, see \cite{Grimmett-Winkler}.

\subsection{Probabilistic interpretation of the matroid stratification of the Grassmannian}

As is well known \cite{Lyons-DPP}, and as we have recalled in Section \ref{sec:dpp}, any subspace $H$ of a finite-dimensional Euclidean $E$, coming with an orthogonal basis~$(e_i)_{i\in S}$, defines a determinantal probability measure on $S$ whose support is the set of bases of a matroid.

This probabilistic point of view implies the following interpretation of the matroid stratification of the Grassmannian of \cite{GGMS}. It can be understood as partitioning $\Gr_n(E)$ by assigning to each matroid $\cM$ on $S$ of rank $n$, the set $\mathcal{G}_{\cM}$ of subspaces $H\in\Gr_n(E)$ whose associated determinantal measure $\P_H$ in the basis $(e_i)_{i\in S}$ has support equal to the set of bases~$\B(\cM)$ of $\cM$. 

On the complex Grassmannian, there is a natural action of the torus $(\C^*)^d/\C^*$ on~$\Gr_n(E)$ by scaling in each direction of the basis~$(e_i)_{i\in S}$, modulo global scaling. This action is Hamiltonian with respect to the natural symplectic structure of the Grassmannian and gives rise to a moment map, which turns out to be a vectorial form of the incidence measure $\mu_H$, restricted to singletons, of the corresponding determinantal measure: 
\[\mu:\Gr_n(E)\to \R^d, \quad H\mapsto \sum_{T\in \B(\cM)} \P_H(\X=T) \big(\sum_{i\in T} e_i\big) =\sum_{i=1}^d \P_{H}(i\in \X) e_i=\sum_{i=1}^d \mu_H(\{i\}) e_i\,.\]

Furthermore, for each matroid $\cM$ of rank $n$, the closure of the image of $\mathcal{G}_{\cM}$ by the moment map~$\mu$ is the matroid polytope of~$\cM$, defined to be the convex hull in $\R^n$ of the collection $\{\sum_{i\in T} e_i : T\in \B(\cM)\}$, see \cite{GGMS, Ardila-ICM}.

\subsection{Infinite volume limits}

Famously, there is an infinite volume limit for uniform spanning trees, due to Burton and Pemantle and further studied in \cite{BLPS} and later works. 

In our setup, we cannot ensure translation invariance. We could consider the model on a sequence of tori for instance, but we would probably find in the limit the usual uniform spanning forest measure.

However, we can define a limiting object corresponding to the closed subspace of $\Omega^1_{\ell^2}(\rG)$ defined by 
\begin{equation}
\bigstar_{\ell^2}(\rG)\oplus \Vect(\theta_1,\ldots,\theta_k)\,,
\end{equation}
where $\bigstar_{\ell^2}(\rG)$ is the $\ell^2$-closure of the image by the discrete derivative of functions on vertices with finite support, and $\theta_1,\ldots,\theta_k$ are linearly independent elements of $\Omega^1_{\ell^2}(\rG)$. Similarly, we can consider the process generated by the closed subspace
\begin{equation}
\bigstar_{\ell^2}(\rG)\cap \left(\Vect(\phi_1,\ldots,\phi_k)\right)^{\perp}\,,
\end{equation}
where $\phi_1,\ldots,\phi_k$ are linearly independent elements of $\Omega^1_{\ell^2}(\rG)$.
Equivalently, and in view of the reparametrization (see Section \ref{sec:reparametrization} and Theorem \ref{thm:sym-forests}), we may also consider $k$ linearly independent zero-mean functions $q_1,\ldots, q_k$ in $\Omega^0_{\ell^2}(\rG)$. We then consider the closed subspace of $\Omega^1_{\ell^2}(\rG)$ defined by 
\[H_{{\sf q}}=\{df: f\in \Omega_{\ell^2}^0(\rG), \langle f,q_i\rangle=0, \forall i\in\{1,\ldots, k\}\}\,.\]

Because of the lack of translation invariance, it does not seem clear to us how to give an easy expression for the correlation kernel, like is the case for spanning trees \cite{Burton-Pemantle} (or its twisted generalizations~\cite{Constantin}). 

Can \eqref{eq:sym-green-intro} be used to say something about the limiting measure? This equation looks like a `finite size' correction to the entropy, expressed in terms of the Green function on the infinite graph: does this have a probabilistic meaning?

Interestingly, eigenvectors of large regular graphs, or random graphs, have well-studied properties of delocalization. This intuitively means that the added or removed edges would be quite delocalized on the graph. Is there a way to translate this into information about the probability measure on the spanning subgraphs?

\subsection{Scaling limits}
The uniform spanning tree (UST) on fine mesh approximations of simply connected planar domains has a well-understood scaling limit described by the Schramm--Loewner evolution (SLE)~\cite{LSW}. Do higher-genus analogues, such as $k$-quasitrees (for $k\le 2g$), admit scaling limits as well? Simulations such as those in Figures~\ref{fig:torus} and~\ref{fig:torus2} (which show uniform $2g$-quasitrees, which are Poincaré duals to uniform spanning trees) suggest this possibility, and this natural question already arises from~\cite{Lyons-Betti}.

Motivated by the construction in Section~\ref{sec:random-g-quasitrees-periods}, one may ask whether random $g$-quasitrees associated with a distinguished Lagrangian subspace -- chosen to ensure self-duality under Poincaré duality and conjectured criticality -- also admit a scaling limit. Since the corresponding partition function can be expressed as a deformation of the spanning tree partition function by the determinant of the discrete normalized period matrix, and that the latter converges to its continuous counterpart under circle pattern approximations~\cite{Mercat, Bobenko-Skopenkov}, this suggests that proving scaling limit might be accessible when the Laplacian determinant is known to converge upon rescaling.

In genus $g=1$, simple closed curves correspond to primitive homology classes $(p,q)\in H_1(\T^2,\Z)\simeq\Z^2$ with $\gcd(p,q) = 1$. Can we characterize the distribution on primitive classes induced by the random simple loop obtained from the determinantal $g$-quasitree? Does this distribution admit a universal scaling limit, depending only on the modulus of the torus?

We could also consider the case of a discretization of a curved disk by a fine graph: for example, take a portion of the round sphere, discretize the $1$-form $dF$, where $F$ is the curvature $2$-form and consider the determinantal random subgraph in $\cC_1(\rG)$ defined by $\im d \oplus \K F$. Does this random spanning subgraph, or the unique random simple loop it defines, have a scaling limit?

Before tackling two-dimensional scaling limits, one could also consider simpler one-dimensional limits. For instance, we could consider the Laplacian spanning forest on a cycle graph with~$n$ edges. Does the process associated with the kernel $\k$, or more likely, the process associated with the kernel $\k^a$ for $0<a<1$ (where $a$ is chosen appropriately with $n$), have a scaling limit on the unit circle when $n$ tends to infinity?



\def\@rst #1 #2other{#1}
\renewcommand\MR[1]{\relax\ifhmode\unskip\spacefactor3000 \space\fi
  \MRhref{\expandafter\@rst #1 other}{#1}}
\renewcommand{\MRhref}[2]{\href{http://www.ams.org/mathscinet-getitem?mr=#1}{MR#1}}
\newcommand{\arXiv}[1]{\href{http://arxiv.org/abs/#1}{arXiv:#1}}
\bibliographystyle{alphaurl}

\bibliography{determinantal-subgraphs}

\begin{thebibliography}{ABBGF16}

\bibitem[ABBGF16]{ABBGF}
Omid Amini, Spencer Bloch, José~I. Burgos~Gil, and Javier Fres\'{a}n.
\newblock Feynman amplitudes and limits of heights.
\newblock {\em Izv. Ross. Akad. Nauk Ser. Mat.}, 80(5):5--40, 2016.
\newblock \href {https://doi.org/10.4213/im8492} {\path{doi:10.4213/im8492}}.

\bibitem[ABKS14]{Baker}
Yang An, Matthew Baker, Greg Kuperberg, and Farbod Shokrieh.
\newblock Canonical representatives for divisor classes on tropical curves and
  the matrix-tree theorem.
\newblock {\em Forum Math. Sigma}, 2:Paper No. e24, 25, 2014.
\newblock \href {https://doi.org/10.1017/fms.2014.25}
  {\path{doi:10.1017/fms.2014.25}}.

\bibitem[AG18]{Avena-Gaudilliere}
Luca Avena and Alexandre Gaudilli\`ere.
\newblock Two applications of random spanning forests.
\newblock {\em J. Theoret. Probab.}, 31(4):1975--2004, 2018.
\newblock \href {https://doi.org/10.1007/s10959-017-0771-3}
  {\path{doi:10.1007/s10959-017-0771-3}}.

\bibitem[Ald90]{Aldous}
David~J. Aldous.
\newblock The random walk construction of uniform spanning trees and uniform
  labelled trees.
\newblock {\em SIAM J. Discrete Math.}, 3(4):450--465, 1990.
\newblock \href {https://doi.org/10.1137/0403039} {\path{doi:10.1137/0403039}}.

\bibitem[Ami19]{Amini-exchange}
Omid Amini.
\newblock The exchange graph and variations of the ratio of the two {S}ymanzik
  polynomials.
\newblock {\em Ann. Inst. Henri Poincar\'{e} D}, 6(2):155--197, 2019.
\newblock \href {https://doi.org/10.4171/AIHPD/68}
  {\path{doi:10.4171/AIHPD/68}}.

\bibitem[AN22]{Amini-Nicolussi}
Omid Amini and Noema Nicolussi.
\newblock {Moduli of hybrid curves II: Tropical and hybrid Laplacians}.
\newblock 2022.
\newblock \arXiv{2203.12785}.

\bibitem[Ard18]{Ardila}
Federico Ardila.
\newblock The geometry of matroids.
\newblock {\em Notices Amer. Math. Soc.}, 65(8):902--908, 2018.

\bibitem[Ard21]{Ardila-ICM}
Federico Ardila.
\newblock The geometry of geometries: matroid theory, old and new.
\newblock {\em ICM 2022 proceedings (submitted)}, 2021.
\newblock \arXiv{2111.08726}.

\bibitem[BBGJ07]{Bouttier-Bowick-Guitter-Jeng}
J{é}r{é}mie Bouttier, Mark Bowick, Emmanuel Guitter, and Monwhea Jeng.
\newblock Vacancy localization in the square dimer model.
\newblock {\em Phys. Rev. E}, 76:041140, Oct 2007.
\newblock \href {https://doi.org/10.1103/PhysRevE.76.041140}
  {\path{doi:10.1103/PhysRevE.76.041140}}.

\bibitem[BBL09]{Borcea-Branden-Liggett}
Julius Borcea, Petter Br\"{a}nd\'{e}n, and Thomas~M. Liggett.
\newblock Negative dependence and the geometry of polynomials.
\newblock {\em J. Amer. Math. Soc.}, 22(2):521--567, 2009.
\newblock \href {https://doi.org/10.1090/S0894-0347-08-00618-8}
  {\path{doi:10.1090/S0894-0347-08-00618-8}}.

\bibitem[BdlHN97]{Nagnibeda}
Roland Bacher, Pierre de~la Harpe, and Tatiana Nagnibeda.
\newblock The lattice of integral flows and the lattice of integral cuts on a
  finite graph.
\newblock {\em Bull. Soc. Math. France}, 125(2):167--198, 1997.
\newblock URL: \url{http://www.numdam.org/item?id=BSMF_1997__125_2_167_0}.

\bibitem[BdTR17]{BdTR}
C\'edric Boutillier, B\'eatrice de~Tili\`ere, and Kilian Raschel.
\newblock The {$Z$}-invariant massive {L}aplacian on isoradial graphs.
\newblock {\em Invent. Math.}, 208(1):109--189, 2017.
\newblock \href {https://doi.org/10.1007/s00222-016-0687-z}
  {\path{doi:10.1007/s00222-016-0687-z}}.

\bibitem[BH20]{Branden-Huh}
Petter Br\"{a}nd\'{e}n and June Huh.
\newblock Lorentzian polynomials.
\newblock {\em Ann. of Math. (2)}, 192(3):821--891, 2020.
\newblock \href {https://doi.org/10.4007/annals.2020.192.3.4}
  {\path{doi:10.4007/annals.2020.192.3.4}}.

\bibitem[Big74]{Biggs-book}
Norman~L. Biggs.
\newblock {\em Algebraic graph theory}.
\newblock Cambridge Univ.\ Press, 1974.
\newblock Cambridge Tracts in Mathematics, No.~67.

\bibitem[Big97]{Biggs-graph}
Norman Biggs.
\newblock Algebraic potential theory on graphs.
\newblock {\em Bull. London Math. Soc.}, 29(6):641--682, 1997.
\newblock \href {https://doi.org/10.1112/S0024609397003305}
  {\path{doi:10.1112/S0024609397003305}}.

\bibitem[Big99]{Biggs-sand}
Norman~L. Biggs.
\newblock Chip-firing and the critical group of a graph.
\newblock {\em J. Algebraic Combin.}, 9(1):25--45, 1999.
\newblock \href {https://doi.org/10.1023/A:1018611014097}
  {\path{doi:10.1023/A:1018611014097}}.

\bibitem[BLPS01]{BLPS}
Itai Benjamini, Russell Lyons, Yuval Peres, and Oded Schramm.
\newblock Uniform spanning forests.
\newblock {\em Ann. Probab.}, 29(1):1--65, 2001.
\newblock \href {https://doi.org/10.1214/aop/1008956321}
  {\path{doi:10.1214/aop/1008956321}}.

\bibitem[Bor11]{Borodin}
Alexei Borodin.
\newblock Determinantal point processes.
\newblock In {\em The {O}xford handbook of random matrix theory}, pages
  231--249. Oxford Univ. Press, Oxford, 2011.

\bibitem[BP93]{Burton-Pemantle}
Robert Burton and Robin Pemantle.
\newblock Local characteristics, entropy and limit theorems for spanning trees
  and domino tilings via transfer-impedances.
\newblock {\em Ann. Probab.}, 21(3):1329--1371, 1993.
\newblock URL:
  \url{http://links.jstor.org/sici?sici=0091-1798(199307)21:3<1329:LCEALT>2.0.CO;2-L&origin=MSN}.

\bibitem[BS16]{Bobenko-Skopenkov}
Alexander Bobenko and Mikhail Skopenkov.
\newblock Discrete {Riemann} surfaces: linear discretization and its
  convergence.
\newblock {\em J. Reine Angew. Math.}, 720:217--250, 2016.
\newblock URL: \url{https://doi.org/10.1515/crelle-2014-0065}.

\bibitem[BTA23]{Barthelme}
Simon Barthelme, Nicolas Tremblay, and Pierre-Olivier Amblard.
\newblock A faster sampler for discrete determinantal point processes.
\newblock In {\em International Conference on Artificial Intelligence and
  Statistics}, pages 5582--5592. PMLR, 2023.

\bibitem[BW10]{Bogner-Weinzierl}
Christian Bogner and Stefan Weinzierl.
\newblock Feynman graph polynomials.
\newblock {\em Internat. J. Modern Phys. A}, 25(13):2585--2618, 2010.
\newblock \href {https://doi.org/10.1142/S0217751X10049438}
  {\path{doi:10.1142/S0217751X10049438}}.

\bibitem[CCK13]{CCK-line}
Michael~J. Catanzaro, Vladimir~Y. Chernyak, and John~R. Klein.
\newblock On {K}irchhoff's theorems with coefficients in a line bundle.
\newblock {\em Homology Homotopy Appl.}, 15(2):267--280, 2013.
\newblock \href {https://doi.org/10.4310/HHA.2013.v15.n2.a16}
  {\path{doi:10.4310/HHA.2013.v15.n2.a16}}.

\bibitem[CCK15]{CCK-complexes}
Michael~J. Catanzaro, Vladimir~Y. Chernyak, and John~R. Klein.
\newblock Kirchhoff's theorems in higher dimensions and {R}eidemeister torsion.
\newblock {\em Homology Homotopy Appl.}, 17(1):165--189, 2015.
\newblock \href {https://doi.org/10.4310/HHA.2015.v17.n1.a8}
  {\path{doi:10.4310/HHA.2015.v17.n1.a8}}.

\bibitem[Cho92]{Chow}
Yutze Chow.
\newblock On enumeration of spanning subgraphs with a preassigned cyclomatic
  number in a graph.
\newblock {\em Acta Math. Hungar.}, 60(1-2):81--91, 1992.
\newblock \href {https://doi.org/10.1007/BF00051759}
  {\path{doi:10.1007/BF00051759}}.

\bibitem[CK23]{CK}
David Cimasoni and Adrien Kassel.
\newblock Graph coverings and twisted operators.
\newblock {\em Algebr. Comb.}, 6(1):75--94, 2023.
\newblock \href {https://doi.org/10.5802/alco.258}
  {\path{doi:10.5802/alco.258}}.

\bibitem[CKS11]{CKS}
Abhijit Champanerkar, Ilya Kofman, and Neal Stoltzfus.
\newblock Quasi-tree expansion for the {B}ollob\'{a}s-{R}iordan-{T}utte
  polynomial.
\newblock {\em Bull. Lond. Math. Soc.}, 43(5):972--984, 2011.
\newblock \href {https://doi.org/10.1112/blms/bdr034}
  {\path{doi:10.1112/blms/bdr034}}.

\bibitem[Con23]{Constantin}
H{\'e}lo{\"i}se Constantin.
\newblock {\em {Spanning Forests and phase transition}}.
\newblock Ph{D} {T}hesis, {\'Ecole normale sup{\'e}rieure de {L}yon}, 2023.
\newblock URL: \url{https://theses.hal.science/tel-04197144}.

\bibitem[CP18]{Corry-Perkinson}
Scott Corry and David Perkinson.
\newblock {\em Divisors and sandpiles}.
\newblock American Mathematical Society, Providence, RI, 2018.
\newblock An introduction to chip-firing.
\newblock \href {https://doi.org/10.1090/mbk/114} {\path{doi:10.1090/mbk/114}}.

\bibitem[DKM15]{Duval-Klivans-Martin}
Art~M. Duval, Caroline~J. Klivans, and Jeremy~L. Martin.
\newblock Cuts and flows of cell complexes.
\newblock {\em J. Algebraic Combin.}, 41(4):969--999, 2015.
\newblock \href {https://doi.org/10.1007/s10801-014-0561-2}
  {\path{doi:10.1007/s10801-014-0561-2}}.

\bibitem[DM21]{Derezinski}
Micha{\l} Derezi{\'n}ski and Michael~W. Mahoney.
\newblock Determinantal point processes in randomized numerical linear algebra.
\newblock {\em Notices Am. Math. Soc.}, 68(1):34--45, 2021.
\newblock \href {https://doi.org/10.1090/noti2202}
  {\path{doi:10.1090/noti2202}}.

\bibitem[FB25]{Bardenet-Fanuel}
Michaël Fanuel and Rémi Bardenet.
\newblock Sparsification of the regularized magnetic {Laplacian} with
  multi-type spanning forests.
\newblock {\em Appl. Comput. Harmon. Anal.}, 78:42, 2025.
\newblock Id/No 101766.
\newblock \href {https://doi.org/10.1016/j.acha.2025.101766}
  {\path{doi:10.1016/j.acha.2025.101766}}.

\bibitem[For93]{Forman}
Robin Forman.
\newblock Determinants of {L}aplacians on graphs.
\newblock {\em Topology}, 32(1):35--46, 1993.
\newblock \href {https://doi.org/10.1016/0040-9383(93)90035-T}
  {\path{doi:10.1016/0040-9383(93)90035-T}}.

\bibitem[GdMN05]{Gimenez-deMier-Noy}
Omer Gim\'{e}nez, Anna de~Mier, and Marc Noy.
\newblock On the number of bases of bicircular matroids.
\newblock {\em Ann. Comb.}, 9(1):35--45, 2005.
\newblock \href {https://doi.org/10.1007/s00026-005-0239-x}
  {\path{doi:10.1007/s00026-005-0239-x}}.

\bibitem[GGMS87]{GGMS}
Israel~M. Gel'fand, R.~Mark Goresky, Robert~D. MacPherson, and Vera~V.
  Serganova.
\newblock Combinatorial geometries, convex polyhedra, and {S}chubert cells.
\newblock {\em Adv. in Math.}, 63(3):301--316, 1987.
\newblock \href {https://doi.org/10.1016/0001-8708(87)90059-4}
  {\path{doi:10.1016/0001-8708(87)90059-4}}.

\bibitem[GJ21]{Guo-Jerrum}
Heng Guo and Mark Jerrum.
\newblock Approximately counting bases of bicircular matroids.
\newblock {\em Combin. Probab. Comput.}, 30(1):124--135, 2021.
\newblock \href {https://doi.org/10.1017/s0963548320000292}
  {\path{doi:10.1017/s0963548320000292}}.

\bibitem[GN06]{Gimenez-Noy}
Omer Gim\'{e}nez and Marc Noy.
\newblock On the complexity of computing the {T}utte polynomial of bicircular
  matroids.
\newblock {\em Combin. Probab. Comput.}, 15(3):385--395, 2006.
\newblock \href {https://doi.org/10.1017/S0963548305007327}
  {\path{doi:10.1017/S0963548305007327}}.

\bibitem[GP14]{Gorodezky-Pak}
Igor Gorodezky and Igor Pak.
\newblock Generalized loop-erased random walks and approximate reachability.
\newblock {\em Random Structures \& Algorithms}, 44(2):201--223, 2014.

\bibitem[GW04]{Grimmett-Winkler}
Geoffrey~R. Grimmett and S~N. Winkler.
\newblock Negative association in uniform forests and connected graphs.
\newblock {\em Random Structures \& Algorithms}, 24(4):444--460, 2004.
\newblock \href {https://doi.org/10.1002/rsa.20012}
  {\path{doi:10.1002/rsa.20012}}.

\bibitem[HH24]{HH}
Noah Halberstam and Tom Hutchcroft.
\newblock Uniqueness of the infinite tree in low-dimensional random forests.
\newblock {\em Probab. Math. Phys.}, 5(4):1185--1216, 2024.
\newblock \href {https://doi.org/10.2140/pmp.2024.5.1185}
  {\path{doi:10.2140/pmp.2024.5.1185}}.

\bibitem[HKPV06]{HKPV}
J.~Ben Hough, Manjunath Krishnapur, Yuval Peres, and B\'alint Vir\'ag.
\newblock Determinantal processes and independence.
\newblock {\em Probab. Surv.}, 3:206--229, 2006.
\newblock \href {https://doi.org/10.1214/154957806000000078}
  {\path{doi:10.1214/154957806000000078}}.

\bibitem[JABT25]{Jaquard}
Hugo Jaquard, Pierre-Olivier Amblard, Simon Barthelm{\'e}, and Nicolas
  Tremblay.
\newblock Random multitype spanning forests for synchronization on sparse
  graphs.
\newblock {\em SIAM Journal on Mathematics of Data Science}, 7(3):1123--1153,
  2025.
\newblock \href {https://doi.org/10.1137/24M1649563}
  {\path{doi:10.1137/24M1649563}}.

\bibitem[Jer24]{Jerrum}
Mark Jerrum.
\newblock {Fundamentals of partial rejection sampling}.
\newblock {\em Probability Surveys}, 21(none):171 -- 199, 2024.
\newblock \href {https://doi.org/10.1214/24-PS29} {\path{doi:10.1214/24-PS29}}.

\bibitem[Kal83]{Kalai}
Gil Kalai.
\newblock Enumeration of {${\bf Q}$}-acyclic simplicial complexes.
\newblock {\em Israel J. Math.}, 45(4):337--351, 1983.
\newblock \href {https://doi.org/10.1007/BF02804017}
  {\path{doi:10.1007/BF02804017}}.

\bibitem[Kas15]{Kassel-ESAIM}
Adrien Kassel.
\newblock Learning about critical phenomena from scribbles and sandpiles.
\newblock In {\em Mod\'{e}lisation {A}l\'{e}atoire et
  {S}tatistique---{J}ourn\'{e}es {MAS} 2014}, volume~51 of {\em ESAIM Proc.
  Surveys}, pages 60--73. EDP Sci., Les Ulis, 2015.
\newblock \href {https://doi.org/10.1051/proc/201551004}
  {\path{doi:10.1051/proc/201551004}}.

\bibitem[Ken11]{Kenyon}
Richard Kenyon.
\newblock Spanning forests and the vector bundle {L}aplacian.
\newblock {\em Ann. Probab.}, 39(5):1983--2017, 2011.
\newblock \href {https://doi.org/10.1214/10-AOP596}
  {\path{doi:10.1214/10-AOP596}}.

\bibitem[Ken19]{Kenyon-massive}
Richard Kenyon.
\newblock Determinantal spanning forests on planar graphs.
\newblock {\em Ann. Probab.}, 47(2):952--988, 2019.
\newblock \href {https://doi.org/10.1214/18-AOP1276}
  {\path{doi:10.1214/18-AOP1276}}.

\bibitem[Kir47]{Kirchhoff}
Gustav Kirchhoff.
\newblock Ueber die {A}ufl\"osung der {G}leichungen, auf welche man bei der
  {U}ntersuchung der linearen {V}ertheilung galvanischer {S}tr\"ome gef\"uhrt
  wird.
\newblock {\em Ann.\ Phys.\ und Chem.}, 72(12):497--508, 1847.
\newblock \href {https://doi.org/10.1002/andp.18471481202}
  {\path{doi:10.1002/andp.18471481202}}.

\bibitem[KK17]{Kassel-Kenyon}
Adrien Kassel and Richard Kenyon.
\newblock Random curves on surfaces induced from the {L}aplacian determinant.
\newblock {\em Ann. Probab.}, 45(2):932--964, 2017.
\newblock \href {https://doi.org/10.1214/15-AOP1078}
  {\path{doi:10.1214/15-AOP1078}}.

\bibitem[KL20]{KL2}
Adrien Kassel and Thierry L\'{e}vy.
\newblock A colourful path to matrix-tree theorems.
\newblock {\em Algebr. Comb.}, 3(2):471--482, 2020.
\newblock \href {https://doi.org/10.5802/alco.100}
  {\path{doi:10.5802/alco.100}}.

\bibitem[KL21]{KL1}
Adrien Kassel and Thierry L\'{e}vy.
\newblock Covariant {S}ymanzik identities.
\newblock {\em Probab. Math. Phys.}, 2(3):419--475, 2021.
\newblock \href {https://doi.org/10.2140/pmp.2021.2.419}
  {\path{doi:10.2140/pmp.2021.2.419}}.

\bibitem[KL22]{KL3}
Adrien Kassel and Thierry L\'evy.
\newblock Determinantal probability measures on {G}rassmannians.
\newblock {\em Ann. Inst. Henri Poincar\'e{} D}, 9(4):659--732, 2022.
\newblock \href {https://doi.org/10.4171/aihpd/152}
  {\path{doi:10.4171/aihpd/152}}.

\bibitem[KL23]{KL4}
Adrien Kassel and Thierry L\'evy.
\newblock On the mean projection theorem for determinantal point processes.
\newblock {\em ALEA, Lat. Am. J. Probab. Math. Stat.}, (20):497--504, (2023).
\newblock \href {https://doi.org/DOI: 10.30757/ALEA.v20-17} {\path{doi:DOI:
  10.30757/ALEA.v20-17}}.

\bibitem[KL26a]{KL2C}
Adrien Kassel and Thierry L{\'e}vy.
\newblock Corrigendum to: ``{A} colourful path to matrix-tree theorems''.
\newblock {\em Algebr. Comb.}, 9(3):595--596, 2026.
\newblock \href {https://doi.org/10.5802/alco.507}
  {\path{doi:10.5802/alco.507}}.

\bibitem[KL26b]{KL6}
Adrien Kassel and Thierry L\'evy.
\newblock Determinantal linear spanning forests.
\newblock 2026.
\newblock In preparation.

\bibitem[KS00]{Kotani-Sunada}
Motoko Kotani and Toshikazu Sunada.
\newblock Jacobian tori associated with a finite graph and its abelian covering
  graphs.
\newblock {\em Adv. in Appl. Math.}, 24(2):89--110, 2000.
\newblock \href {https://doi.org/10.1006/aama.1999.0672}
  {\path{doi:10.1006/aama.1999.0672}}.

\bibitem[KW15]{KaWu}
Adrien Kassel and Wei Wu.
\newblock Transfer current and pattern fields in spanning trees.
\newblock {\em Probab. Theory Relat. Fields}, 163(1-2):89--121, 2015.

\bibitem[KW16]{Kassel-Wilson-sand}
Adrien Kassel and David~B. Wilson.
\newblock The looping rate and sandpile density of planar graphs.
\newblock {\em Amer. Math. Monthly}, 123(1):19--39, 2016.
\newblock \href {https://doi.org/10.4169/amer.math.monthly.123.1.19}
  {\path{doi:10.4169/amer.math.monthly.123.1.19}}.

\bibitem[LC81]{Liu-Chow}
Chyan-Jyue Liu and Yutze Chow.
\newblock Enumeration of forests in a graph.
\newblock {\em Proc. Amer. Math. Soc.}, 83(3):659--662, 1981.
\newblock \href {https://doi.org/10.2307/2044142} {\path{doi:10.2307/2044142}}.

\bibitem[LLY26]{Lam-Lo-Yuen}
Wai~Yeung Lam, On-Hei~Solomon Lo, and Chi~Ho Yuen.
\newblock Period matrices and homological quasi-trees on discrete riemann
  surfaces.
\newblock {\em Journal für die reine und angewandte Mathematik (Crelles
  Journal)}, 2026(837):223--251, 2026.
\newblock \href {https://doi.org/10.1515/crelle-2026-0033}
  {\path{doi:10.1515/crelle-2026-0033}}.

\bibitem[LP16]{Lyons-Peres}
Russell Lyons and Yuval Peres.
\newblock {\em Probability on Trees and Networks}, volume~42 of {\em Cambridge
  Series in Statistical and Probabilistic Mathematics}.
\newblock Cambridge University Press, New York, 2016.
\newblock Available at \url{http://pages.iu.edu/~rdlyons/ }.
\newblock \href {https://doi.org/10.1017/9781316672815}
  {\path{doi:10.1017/9781316672815}}.

\bibitem[LSW04]{LSW}
Gregory~F. Lawler, Oded Schramm, and Wendelin Werner.
\newblock Conformal invariance of planar loop-erased random walks and uniform
  spanning trees.
\newblock {\em Ann. Probab.}, 32(1B):939--995, 2004.
\newblock \href {https://doi.org/10.1214/aop/1079021469}
  {\path{doi:10.1214/aop/1079021469}}.

\bibitem[Lyo03]{Lyons-DPP}
Russell Lyons.
\newblock Determinantal probability measures.
\newblock {\em Publ. Math. Inst. Hautes \'Etudes Sci.}, (98):167--212, 2003.
\newblock \href {https://doi.org/10.1007/s10240-003-0016-0}
  {\path{doi:10.1007/s10240-003-0016-0}}.

\bibitem[Lyo09]{Lyons-Betti}
Russell Lyons.
\newblock Random complexes and {$\ell^2$}-{B}etti numbers.
\newblock {\em J. Topol. Anal.}, 1(2):153--175, 2009.
\newblock \href {https://doi.org/10.1142/S1793525309000072}
  {\path{doi:10.1142/S1793525309000072}}.

\bibitem[Mac75]{Macchi}
Odile Macchi.
\newblock The coincidence approach to stochastic point processes.
\newblock {\em Advances in Appl. Probability}, 7:83--122, 1975.
\newblock \href {https://doi.org/10.2307/1425855} {\path{doi:10.2307/1425855}}.

\bibitem[Mat77]{Matthews}
Laurence~R. Matthews.
\newblock Bicircular matroids.
\newblock {\em Quart. J. Math. Oxford Ser. (2)}, 28(110):213--227, 1977.
\newblock \href {https://doi.org/10.1093/qmath/28.2.213}
  {\path{doi:10.1093/qmath/28.2.213}}.

\bibitem[Mau76]{Maurer}
Stephen~B. Maurer.
\newblock Matrix generalizations of some theorems on trees, cycles and cocycles
  in graphs.
\newblock {\em SIAM J. Appl. Math.}, 30(1):143--148, 1976.
\newblock \href {https://doi.org/10.1137/0130017} {\path{doi:10.1137/0130017}}.

\bibitem[Mer07]{Mercat}
Christian Mercat.
\newblock Discrete {Riemann} surfaces.
\newblock In {\em Handbook of Teichm\"uller theory. Volume I}, pages 541--575.
  Z{\"u}rich: European Mathematical Society (EMS), 2007.
\newblock \href {https://doi.org/10.4171/029-1/14}
  {\path{doi:10.4171/029-1/14}}.

\bibitem[MS18]{Moffatt-Smith}
Iain Moffatt and Ben Smith.
\newblock Matroidal frameworks for topological {T}utte polynomials.
\newblock {\em J. Combin. Theory Ser. B}, 133:1--31, 2018.
\newblock \href {https://doi.org/10.1016/j.jctb.2017.09.009}
  {\path{doi:10.1016/j.jctb.2017.09.009}}.

\bibitem[Myr92]{Myrvold}
Wendy Myrvold.
\newblock Counting {$k$}-component forests of a graph.
\newblock {\em Networks}, 22(7):647--652, 1992.
\newblock \href {https://doi.org/10.1002/net.3230220704}
  {\path{doi:10.1002/net.3230220704}}.

\bibitem[NS61]{Nerode-Shank}
Anil Nerode and Herbert Shank.
\newblock An algebraic proof of {K}irchhoff's network theorem.
\newblock {\em Amer. Math. Monthly}, 68:244--247, 1961.
\newblock \href {https://doi.org/10.2307/2311455} {\path{doi:10.2307/2311455}}.

\bibitem[Oxl11]{Oxley}
James Oxley.
\newblock {\em Matroid theory}, volume~21 of {\em Oxford Graduate Texts in
  Mathematics}.
\newblock Oxford University Press, Oxford, second edition, 2011.
\newblock \href {https://doi.org/10.1093/acprof:oso/9780198566946.001.0001}
  {\path{doi:10.1093/acprof:oso/9780198566946.001.0001}}.

\bibitem[PB83]{Provan-Ball}
J.~Scott Provan and Michael~O. Ball.
\newblock The complexity of counting cuts and of computing the probability that
  a graph is connected.
\newblock {\em SIAM J. Comput.}, 12(4):777--788, 1983.
\newblock \href {https://doi.org/10.1137/0212053} {\path{doi:10.1137/0212053}}.

\bibitem[Pem91]{Pemantle}
Robin Pemantle.
\newblock Choosing a spanning tree for the integer lattice uniformly.
\newblock {\em Ann. Probab.}, 19(4):1559--1574, 1991.
\newblock URL:
  \url{http://links.jstor.org/sici?sici=0091-1798(199110)19:4<1559:CASTFT>2.0.CO;2-E&origin=MSN}.

\bibitem[Piq19]{Piquerez}
Matthieu Piquerez.
\newblock A multidimensional generalization of {S}ymanzik polynomials.
\newblock 2019.
\newblock \arXiv{1901.09797}.

\bibitem[Sch00]{Schramm}
Oded Schramm.
\newblock Scaling limits of loop-erased random walks and uniform spanning
  trees.
\newblock {\em Israel J. Math.}, 118:221--288, 2000.
\newblock \href {https://doi.org/10.1007/BF02803524}
  {\path{doi:10.1007/BF02803524}}.

\bibitem[SP72]{Simoes-Pereira}
Jos{é} Manuel dos~Santos Sim{õ}es-Pereira.
\newblock On subgraphs as matroid cells.
\newblock {\em Math. Z.}, 127:315--322, 1972.
\newblock \href {https://doi.org/10.1007/BF01111390}
  {\path{doi:10.1007/BF01111390}}.

\bibitem[Tut54]{Tutte}
William~T. Tutte.
\newblock A contribution to the theory of chromatic polynomials.
\newblock {\em Canad. J. Math.}, 6:80--91, 1954.
\newblock \href {https://doi.org/10.4153/cjm-1954-010-9}
  {\path{doi:10.4153/cjm-1954-010-9}}.

\bibitem[Val79]{Valiant}
Leslie~G. Valiant.
\newblock The complexity of computing the permanent.
\newblock {\em Theoret. Comput. Sci.}, 8(2):189--201, 1979.
\newblock \href {https://doi.org/10.1016/0304-3975(79)90044-6}
  {\path{doi:10.1016/0304-3975(79)90044-6}}.

\bibitem[Wel93]{Welsh}
Dominic J.~A. Welsh.
\newblock {\em Complexity: knots, colourings and counting}, volume 186 of {\em
  London Mathematical Society Lecture Note Series}.
\newblock Cambridge University Press, Cambridge, 1993.
\newblock \href {https://doi.org/10.1017/CBO9780511752506}
  {\path{doi:10.1017/CBO9780511752506}}.

\bibitem[Whi35]{Whitney}
Hassler Whitney.
\newblock On the {A}bstract {P}roperties of {L}inear {D}ependence.
\newblock {\em Amer. J. Math.}, 57(3):509--533, 1935.
\newblock \href {https://doi.org/10.2307/2371182} {\path{doi:10.2307/2371182}}.

\bibitem[Whi87]{White}
Neil White, editor.
\newblock {\em Combinatorial geometries}, volume~29 of {\em Encyclopedia of
  Mathematics and its Applications}.
\newblock Cambridge University Press, Cambridge, 1987.
\newblock \href {https://doi.org/10.1017/CBO9781107325715}
  {\path{doi:10.1017/CBO9781107325715}}.

\bibitem[Wil96]{Wilson}
David~Bruce Wilson.
\newblock Generating random spanning trees more quickly than the cover time.
\newblock In {\em Proceedings of the {T}wenty-eighth {A}nnual {ACM} {S}ymposium
  on the {T}heory of {C}omputing ({P}hiladelphia, {PA}, 1996)}, pages 296--303.
  ACM, New York, 1996.
\newblock \href {https://doi.org/10.1145/237814.237880}
  {\path{doi:10.1145/237814.237880}}.

\bibitem[Zha93]{Zhang}
Shouwu Zhang.
\newblock Admissible pairing on a curve.
\newblock {\em Invent. Math.}, 112(1):171--193, 1993.
\newblock \href {https://doi.org/10.1007/BF01232429}
  {\path{doi:10.1007/BF01232429}}.

\end{thebibliography}


\end{document}